\documentclass[10pt]{article} 
 
\usepackage{latexsym} 
\usepackage{amsfonts} 
\usepackage{amsmath} 
\usepackage{amsthm} 
\usepackage{mathrsfs,esint,enumitem} 
\usepackage{dsfont} 
\usepackage{bbold} 
\usepackage[english]{babel} 
\usepackage{caption} 
\usepackage{epsfig} 
\usepackage{float} 
\usepackage{subfigure} 
\usepackage{psfrag} 
\usepackage{graphicx} 
\usepackage{epsfig} 
\usepackage{amsbsy} 
\usepackage{hyperref,color}
\usepackage{setspace}

\DeclareMathOperator{\Val}{\matV} 
\DeclareMathOperator{\meas}{meas}

\newtheorem{theorem}{Theorem} 
\newtheorem*{prop*}{Theorem} 
\newtheorem{theo}[theorem]{Theorem} 
\newtheorem{coro}[theorem]{Corollary} 
\newtheorem{defi}[theorem]{Definition} 
\newtheorem{lemma}[theorem]{Lemma} 
\newtheorem{prop}[theorem]{Proposition} 
\newtheorem{rmk}[theorem]{Remark}

\newtheorem{ex}[theorem]{Example}
 
\numberwithin{theorem}{section}

\newcommand{\zerarcounters}{\setcounter{equation}{0}\setcounter{theorem}{0}}

\newcommand{\beq}{\begin{equation}}
\newcommand{\eeq}{\end{equation}}

\newcommand{\ZZZ}{\mathds{Z}} 
\newcommand{\CCC}{\mathds{C}} 
\newcommand{\NNN}{\mathds{N}} 
 
\newcommand{\RRR}{\mathds{R}} 
\newcommand{\TTT}{\mathds{T}}

\newcommand{\BB}{{\mathcal B}} 
 
\newcommand{\DD}{{\mathcal D}} 
 
\newcommand{\calF}{{\mathcal F}} 
\newcommand{\calG}{{\mathcal G}} 
\newcommand{\calH}{{\mathcal H}} 
\newcommand{\calI}{{\mathcal I}}

\newcommand{\LL}{{\mathcal L}} 
\newcommand{\MM}{{\mathcal M}}

\newcommand{\calP}{{\mathcal P}}

\newcommand{\SSSS}{{\mathcal S}}

\newcommand{\calW}{{\mathcal W}}

\DeclareFontFamily{U}{BOONDOX-calo}{\skewchar\font=45 }
\DeclareFontShape{U}{BOONDOX-calo}{m}{n}{
  <-> s*[1.05] BOONDOX-r-calo}{}
\DeclareMathAlphabet{\mathcalboondox}{U}{BOONDOX-calo}{m}{n}
\DeclareMathAlphabet{\mathbcalboondox}{U}{BOONDOX-calo}{b}{n}

\newcommand{\gotn}{{\mathfrak n}}

\newcommand{\gots}{{\mathfrak s}}

\newcommand{\gotB}{{\mathbcalboondox B}} 
\newcommand{\gotC}{{\mathfrak C}}

\newcommand{\gotF}{{\mathfrak F}}

\newcommand{\gotN}{{\mathfrak N}}

\newcommand{\matA}{{\mathscr A}}

\newcommand{\matF}{{\mathscr F}} 
 
\newcommand{\matH}{{\mathscr H}}

\newcommand{\matL}{{\mathscr L}} 
\newcommand{\matM}{{\mathscr M}}

\newcommand{\matR}{{\mathscr R}}

\newcommand{\matU}{{\mathscr U}} 
\newcommand{\matV}{{\mathscr V}}

\newcommand{\matW}{{\mathscr W}}

\newcommand{\und}{\underline}

\newcommand{\ol}{\overline} 
 
\newcommand{\Fullbox}{{\rule{2.0mm}{2.0mm}}} 
 
\newcommand{\EP}{\hfill\Fullbox\vspace{0.2cm}} 
\newcommand{\prova}{\noindent{\it Proof. }} 
\newcommand{\io}{\infty} 
\newcommand{\e}{\varepsilon} 
\newcommand{\al}{\alpha} 
\newcommand{\de}{\delta} 
\newcommand{\be}{\beta} 
 
\newcommand{\m}{\mu} 
\newcommand{\x}{\xi} 
\newcommand{\p}{\pi} 
 
\newcommand{\g}{\gamma} 
\newcommand{\om}{\omega} 
\newcommand{\h}{\eta}

\newcommand{\f}{\varphi} 
\newcommand{\s}{\sigma} 
 
\newcommand{\del}{\partial}

\newcommand{\avg}[1]{\langle #1 \rangle}

\newcommand{\oo}{{\omega}}

\newcommand{\nn}{{\nu}}

\newcommand{\der}{{\rm d}} 
\newcommand{\ii}{{\rm i}}

\newcommand{\jap}[1]{\langle #1 \rangle}

\makeatletter
\newcommand{\oset}[3][0ex]{%
  \mathrel{\mathop{#3}\limits^{
    \vbox to#1{\kern-1\ex@
    \hbox{$\scriptstyle#2$}\vss}}}}
\makeatother

\def\tilde#1{\widetilde{#1}}
 
\def\ins#1#2#3{\vbox to0pt{\kern-#2 \hbox{\kern#1 #3}\vss}\nointerlineskip} 

\def\tilde#1{\widetilde{#1}}
 
\def\ins#1#2#3{\vbox to0pt{\kern-#2 \hbox{\kern#1 #3}\vss}\nointerlineskip} 
 
\begin{document}
%%%%%%%%%%%%%%%%%%%%%%%%%%%%%%%%%%%%%%%%%%%%%%%%%%%%%%%%%%%%%%%%%%%%%%%%%% 
%%%%%%%%%%%%%%%%%%%%%%%%%%%%%%%%%%%%%%%%%%%%%%%%%%%%%%%%%%%%%%%%%%%%%%%%%% 
 
%%%%%%%%%%%%%%%%%%%%%%%%%%%%%%%%%%%%%%%%%%%%%%%%%%%%%%%%%%%%%%%%%%%%%%%%%% 
%%%%%%%%%%%%%%%%%%%%%%%%%%%%%%%%%%%%%%%%%%%%%%%%%%%%%%%%%%%%%%%%%%%%%%%%%% 
\title{\bf Maximal tori in infinite-dimensional Hamiltonian systems:
a Renormalization Group approach}
%%%%%%%%%%%%%%%%%%%%%%%%%%%%%%%%%%%%%%%%%%%%%%%%%%%%%%%%%%%%%%%%%%%%%%%%%% 
%%%%%%%%%%%%%%%%%%%%%%%%%%%%%%%%%%%%%%%%%%%%%%%%%%%%%%%%%%%%%%%%%%%%%%%%%% 
 
\author{
\vspace{.2cm}
\textbf{Livia Corsi, Guido Gentile, Michela Procesi}\\
\small Dipartimento di Matematica e Fisica, Universit\`a Roma Tre, Roma, 
00146, Italy\\
\small e-mail: livia.corsi@uniroma3.it,
guido.gentile@uniroma3.it,  michela.procesi@uniroma3.it}

\date{} 
 
\maketitle 

%%%%%%%%%%%%%%%%%%%%%%%%%%%%%%%%%%%%%%%%%%%%%%%%%%%%%%%%%%%%%%%%%%%%%%%%%% 
\begin{abstract} 
We study the existence of infinite-dimensional invariant tori in a mechanical system of infinitely many rotators weakly interacting with each other.
We consider explicitly interactions depending only on the angles,
with the aim of discussing in a simple case the analyticity properties to be required on the perturbation of the integrable system
in order to ensure the persistence of a large measure set of invariant tori with finite energy. 
The proof we provide of the persistence of the invariant  tori implements the Renormalization Group scheme
based on the tree formalism -- i.e.~the graphical representation of the solutions of the equations of motion in terms of trees --
which has been widely used in finite-dimensional problems. The method is very effectual and flexible:
it naturally extends, once the functional setting has been fixed, to the infinite-dimensional case
with only minor technical-natured adaptations.
\end{abstract} 
%%%%%%%%%%%%%%%%%%%%%%%%%%%%%%%%%%%%%%%%%%%%%%%%%%%%%%%%%%%%%%%%%%%%%%%%%% 
 
\begin{spacing}{0.0}
\tableofcontents
\end{spacing}
 
%%%%%%%%%%%%%%%%%%%%%%%%%%%%%%%%%%%%%%%%%%%%%%%%%%%%%%%%%%%%%%%%%%%%%%%%%% 
%%%%%%%%%%%%%%%%%%%%%%%%%%%%%%%%%%%%%%%%%%%%%%%%%%%%%%%%%%%%%%%%%%%%%%%%%% 
 \zerarcounters 
\section{Introduction} 
\label{intro} 
%%%%%%%%%%%%%%%%%%%%%%%%%%%%%%%%%%%%%%%%%%%%%%%%%%%%%%%%%%%%%%%%%%%%%%%%%% 
%%%%%%%%%%%%%%%%%%%%%%%%%%%%%%%%%%%%%%%%%%%%%%%%%%%%%%%%%%%%%%%%%%%%%%%%%% 

One of the great  achievements of the past century in the context of finite-dimensional Hamiltonian Mechanics is the
celebrated KAM theorem on the persistence of Lagrangian tori in quasi-integrable systems \cite{Kol,Arn,Mos}.
The theorem states that, under an appropriate non-degeneracy condition -- the so-called
\emph{twist condition} -- a large set of Lagrangian tori survive any small enough perturbation. It is therefore natural
to ask whether a similar result holds also in infinite dimension. In fact, by looking at the proof of the classical KAM 
theorem, the smallness condition on the perturbation strongly depends on the dimension and there is no way
to pass directly to the limit. This is not surprising, because in the infinite-dimensional context
most statements heavily depend on the topology, and even the concept of Lagrangian torus needs to be clarified. 

The persistence of a positive measure set (in some topology) of invariant tori has deep dynamical consequences
on the stability of the whole system, and therefore it is clearly of interest.
Note that in  the literature on the topic all results are in quite strong topology; see for instance
refs.~\cite{Po,Bjfa,Yuan,BMP1,cong2} for the existence of almost periodic solutions and refs.~\cite{BaGr,FaGr} for results on stability.
This might suggest that only very special invariant tori survive if the topology is not strong enough;
see for instance ref.~\cite{BMP2} for the construction of non-maximal tori with only Sobolev regularity.

Let us consider, as a toy model, the Hamiltonian system consisting in a set of rotators weakly interacting with each other.
The system is described, in the standard Darboux coordinates, by the Hamiltonian
\begin{equation}\label{ham0}
H(\theta,I)= \frac{1}{2}I\cdot I + \e f(\theta),\qquad (\theta,I)\in\TTT^d\times\RRR^d,
\end{equation}  
with the perturbation (i.e.~the interaction function $f$) assumed to be real analytic.
The Hamiltonian \eqref{ham0} has been explicitly
studied in the literature as a paradigmatic example to illustrate the main features of KAM theory
avoiding as much as possible any technical intricacies
\cite{Th,Galla,GG,GM1,CF,BC1,BC2,BC3}.
We mention that last three references focus on secondary tori, i.e.~non-contractible tori which appear only for $\e\neq0$, while
here we concentrate on primary tori, i.e.~tori winding around $\TTT^d$.

To illustrate the result in such a finite-dimensional case, we first recall that a function $f$ is real analytic
if it admits a holomorphic extension to a complex neighbourhood of $\TTT^d$;
since all norms are equivalent in finite dimension, this is tantamount to requiring that,
for any choice of the norm $|\cdot|_\star$ and for some positive constant $s$, one has
\begin{equation}\label{fouriero}
f(\theta)=\sum_{\nu\in\ZZZ^d} f_\nu e^{\ii\theta\cdot\nu},
\qquad \ol{f_\nu}= f_{-\nu} ,
\qquad M_f :=\sum_{\nu\in\ZZZ^d} e^{2s|\nu|_\star}|f_\nu| < \io, 
\end{equation}
where $\ol z$ denotes the complex conjugate of $z$.

Then, any torus with frequency vector $\om$ satisfying a suitable non-resonance condition survives any perturbation $f$
if $\e$ is smaller than a threshold value depending on $\om$. Specifically, if $\om\in\RRR^d$ satisfies, say,
the standard Diophantine condition\footnote{Here and in the following $\|\cdot\|_p$ denotes as usual the $\ell^p$ norm,
while $\|\cdot\|$ is the sup-norm.}
\begin{equation}\label{diop}
|\om\cdot\nu|>\frac{\g}{\|\nu\|_1^\tau}\qquad \forall \,\nu\in\ZZZ^d\setminus\{0\} ,
\end{equation}
for some $\g>0$ and $\tau>d-1$, then, for any real analytic
perturbation $f$ satisfying \eqref{fouriero} with $|\cdot|_\star=\|\cdot\|_1$,
the perturbed system admits an invariant torus with frequency $\om$, provided that one has
\begin{equation}\label{ddd}
|\e|< \e_0(d,\g) , \qquad \e_0(d,\g) \le \frac{C \g^2}{(d!)^{\al}}\ ,
\end{equation}
for suitable constants $C=C(s,M_f)$ and $\al$. 
Moreover, using that the complement of the set of frequencies satisfying \eqref{diop} 
has probability measure of order $O(\g)$, one finds that the complement of the set of invariant  primary tori is of order $O({\sqrt{\e}})$;
see for instance refs.~\cite{Kou,CK2} for a recent and quantitative presentation.
The estimate \eqref{ddd} can be improved, especially if one assumes the interaction to be short range~\cite{Way},
but still, in general, the best bounds for $\e_0(d,\g)$ one is able to obtain
go to zero at least exponentially in $d$~\cite{CK1}.
Thus, it is clear that the case $d\to+\io$ requires some refinement.

To deal with the infinite-dimensional limit one has to make some strong assumptions either on $f$ or on $\om$.
To envision this,  consider an interaction of the form
\begin{equation}\label{poschel}
f(\theta) = \sum_{j=1}^D \e_j f_j(\theta_1,\ldots,\theta_{d_j}),\qquad d_j<d_{j+1} \quad \forall \, j=1,\ldots,D-1 .
\end{equation}
Then, if, for each $j=1,\ldots,D$, one assumes a Diophantine condition
\begin{equation}\label{diopj}
|\varpi_j\cdot\nu|>\frac{\g_j}{\|\nu\|_1^{\tau_j}}\qquad \forall \,\nu\in\ZZZ^{d_j}\setminus\{0\} ,
\qquad \varpi_j=(\om_1,\ldots,\om_{d_j}) ,
\end{equation}
and conditions like \eqref{ddd} are required to hold with
$\e$, $\g$ and $d$ replaced with $\e_j$, $\g_j$ and $d_j$, respectively,
for $j=1,\ldots,D$,
one can expect to be able to prove  persistence of tori also in the limit $D\to+\io$, by a minor modification
of any classical KAM algorithm.\footnote{
We mention also that, in the infinite-dimensional case, Berti and Biasco provided a super-exponential bound
of the maximum size of the perturbation allowed for the persistence of finite-dimensional tori in terms of the dimensions of the tori~\cite{BeBi}.
The bound was slightly improved subsequently by Li and Liu \cite{LiLi}.}
For instance, one may assume that $\e_j$ goes to zero fast enough, e.g.~$\e_j\sim (d_j!)^{-\al}$;
this is in essence the idea exploited by Poschel  \cite{Po} to prove the existence of almost-periodic solutions
to a nonlinear Schr\"odinger-type equation.
Another possibility consists in taking $\g_j$ increasingly larger when increasing $j$, for instance of order $\sim (d_j!)^{\al/2}$;
in this case, for $\varpi_j$ to satisfy \eqref{diopj}, one needs $\om_j$ to grow at least super-exponentially.
An approach of this type was proposed by Chierchia and Perfetti \cite{CP}.
 
However it is worth noticing that, although in finite dimension all norms are equivalent,
the equivalence constant depends in fact on the dimension.
This means that we have some freedom in requiring the regularity condition on $f$ and the non-resonance condition on $\om$,
and it may be possible that, by choosing carefully the norm $|\cdot|_\star$ in \eqref{fouriero}
and replacing \eqref{diop} with a different condition involving such a norm, 
one gets a bound on $\e$ uniform in $d$. This approach is at the basis of Bourgain's result \cite{Bjfa} and of the present paper.

To highlight the dependence on the topology of \eqref{diop} it is actually more convenient to consider Bryuno vectors
instead of Diophantine ones; this allows us to introduce a notion of non-resonance which turns out to be more natural and,
at the same time, provides stronger results.
Introducing the non-increasing sequence $\be^\star_\om:=\{\be^\star_\om (m)\}_{m\ge0}$ and the 
function $\BB^\star(\om)$ by setting
\begin{equation} \nonumber %\label{beta0}
\beta^\star_\om(m):= \inf_{\substack{\nu\in\ZZZ^{d} \vspace{-.00cm}\\ 0 <|\nu|_{\star}\le 2^m}} \!\!\!\!\!\!  |\om\cdot\nu| , 
\qquad\qquad 
\BB^\star(\om) := \sum_{m\ge1}\frac{1}{2^m}\log \left(\frac{1}{\be^\star_\om(m)}\right) ,
\end{equation}
we say that $\om$ is a Bryuno vector if $\BB^\star(\om)<\io$. If $|\cdot|_\star=\|\cdot\|_1$ this is the original definition used by Bryuno,
who proved that such a condition is sufficient for the persistence of invariant tori for $\e$ small enough depending on $\om$ and $d$,
and moreover that Diophantine vectors are also Bryuno \cite{Bry}.

The aim of the present paper is to provide sufficient conditions on the norm $|\cdot|_\star$ in such a way that
\begin{enumerate}[topsep=0.5ex]
	\itemsep-0.1em
	\item if $\om$ is a Bryuno frequency vector one has a bound on the maximum value of $\e$ uniform in $d$;
	\item the corresponding torus has uniformly bounded energy in $d$;
	\item the set of the Bryuno frequency vectors  has positive measure.
\end{enumerate}
With this in mind, we may consider directly the case $d=+\io$, provided one introduces an appropriate functional setting 
to deal with functions of infinitely many angles.
In order for the Hamiltonian \eqref{ham0} to make sense in the infinite-dimensional setting,
we need, at least, both $I$ to be a sequence in $\ell^2(\RRR)$ and $f$ to be the uniform limit of functions of the form \eqref{poschel}.
For our result however we need further conditions, namely
\begin{enumerate}[topsep=0.5ex]
	\itemsep-0.1em
	\item for item 1 above we need that $|\nu|_\star:= \sum_j h_j |\nu_j|$ with $h_j\to +\io$;
	\item for item 2 above we need that $h_j$ goes to $+\io$ faster than $\log j$;
	\item for item 3 above we need that $h_j$ goes to $+\io$ faster than $(\log j)^2$.
\end{enumerate}

The technique we rely upon is the tree formalism originally introduced by Eliasson \cite{E}
and, more systematically, by Gallavotti \cite{Galla} to discuss the finite-dimensional KAM theorem.
Such a technique, based on Renormalization Group ideas usually implemented in Quantum Field Theory \cite{GM2},
turns out to be very flexible and can be applied to study quasi-integrable systems, whether of finite or infinite dimensions,
essentially with minimal differences, which are mostly of a technical nature and involve only the functional setting.
Recently, the tree formalism has been applied to prove the existence of almost-periodic solutions
to the one-dimensional NLS equation with a convolution potential of class $C^p$, for any integer $p$ \cite{CGP}.
The infinite-dimensional system we study in this paper provides a simpler model, where most of the basic ingredients of the tree formalism
can be exposed more easily. In particular, the twist condition of the unperturbed Hamiltonian allows to fix the frequency vector since the beginning,
without the need to introduce any counterterm and to deal with the resulting, quite non-trivial implicit function problem.
Therefore, the present paper can also be considered as an introduction to the tree formalism and to the method
used in the quoted reference to attack a technically more involved  problem.

In this paper we focus on the case of analytic Hamiltonian systems.
Recently, persistence of infinite-dimensional tori has been obtained
with KAM techniques, in perturbations of linear
first order equations, also in the $C^\io$ case \cite{TL}.

%%%%%%%%%%%%%%%%%%%%%%%%%%%%%%%%%%%%%%%%%%%%%%%%%%%%%%%%%%%%%%%%%%%%%%%%%% 
%%%%%%%%%%%%%%%%%%%%%%%%%%%%%%%%%%%%%%%%%%%%%%%%%%%%%%%%%%%%%%%%%%%%%%%%%% 
\zerarcounters
\section{Functional setting and main result}
\label{functional}
%%%%%%%%%%%%%%%%%%%%%%%%%%%%%%%%%%%%%%%%%%%%%%%%%%%%%%%%%%%%%%%%%%%%%%%%%% 
%%%%%%%%%%%%%%%%%%%%%%%%%%%%%%%%%%%%%%%%%%%%%%%%%%%%%%%%%%%%%%%%%%%%%%%%%% 

As outlined in Section \ref{intro}, if we want to consider an infinite-dimensional version of  the Hamiltonian \eqref{ham0},
we need to define carefully the functional setting. We start be introducing the proper Banach space.

%%%%%%%%%%%%%%%%%%%%%%%%%%%%%%%%%%%%%%%%%%%%%%%%%%%%%%%%%%%%%%%%%%%%%%%%% 
\begin{defi} 
\label{normstar}
Set\footnote{The subscript $f$ stands for ``finite support'', and it has nothing to do with the perturbation $f$ we introduce in \eqref{ham}.}
\[
\ZZZ^{\ZZZ}_{f}:=  \Bigl\{ \nu \in \ZZZ^\ZZZ : \|\nu\|_1 < \infty \Bigr\}
\]
and define in $\ZZZ^\ZZZ_f$ the \emph{$\star$-norm} as
\begin{equation}\label{normalfa}
|\nu|_\star:= \sum_{j\in\ZZZ} h_j |\nu_j|,\quad  h_j= h_{-j} \in \RRR_+ \; \forall\,j\in\ZZZ , 
\quad h_{j+1} \ge h_j \; \forall j \in \ZZZ_+ ,
\quad 
\lim_{j\to +\io} h_j= +\io .
\end{equation}
\end{defi}
%%%%%%%%%%%%%%%%%%%%%%%%%%%%%%%%%%%%%%%%%%%%%%%%%%%%%%%%%%%%%%%%%%%%%%%%%% 

%%%%%%%%%%%%%%%%%%%%%%%%%%%%%%%%%%%%%%%%%%%%%%%%%%%%%%%%%%%%%%%%%%%%%%%%% 
\begin{rmk} \label{hj-infty1}
\emph{
If we allow $h_j \in \RRR_+\cup\{+\io\}$ in \eqref{normalfa} and interpret $h_j |\nu_j|=0$ when $\nu_j=0$
for any value of $h_j$ ($+\io$ being included), then the persistence of invariant tori in the case of the finite-dimensional system with Hamiltonian
\eqref{ham0} can be seen as a special case of the general result we discuss in this paper (see Remark \ref{classico}).
As explicit and paradigmatic examples of $\star$-norms, one may consider the cases\footnote{Here and in the following $\jap{j}=\max\{1,|j|\}$.}
$h_j=(\log(1+\jap{j}))^{\s}$, with $\s>1$, and $h_j=\jap{j}^{\al}$, with $\al>0$.
}
\end{rmk}
%%%%%%%%%%%%%%%%%%%%%%%%%%%%%%%%%%%%%%%%%%%%%%%%%%%%%%%%%%%%%%%%%%%%%%%%%%% 

Next, following Montalto and Procesi \cite{MP}, we identify the space where to define the dynamical system
we want to study and the solutions we want to investigate.

%%%%%%%%%%%%%%%%%%%%%%%%%%%%%%%%%%%%%%%%%%%%%%%%%%%%%%%%%%%%%%%%%%%%%%%%% 
\begin{defi} \label{spaceH}
Given the infinite-dimensional torus 
\begin{equation} \nonumber 
\TTT^\ZZZ= \Bigl\{ \theta=\{\theta_j\}_{j\in\ZZZ}\ \colon \theta_j\in \RRR/2\pi\ZZZ \Bigr\} , 
\end{equation}
endowed with the metric\footnote{In fact,
the topology induced by the metric \eqref{dist}
makes $\TTT^\ZZZ$ a
Banach manifold modeled on $\ell^\io(\RRR)$.} 
\begin{equation} \label{dist}
{\rm dist}(\theta,\theta')= \sup_{j\in\ZZZ} \inf_{k\in \ZZZ} |\theta_j-\theta'_j- 2\pi k|
\end{equation}
inherited form $\ell^\io(\RRR)$, and given a Banach space $(X, |\cdot|_X)$, for any $\gots>0$ we set
\begin{equation} \label{Hs}
\matH_\star^{\gots}(\TTT^\ZZZ,X) :=\Biggl\{ u\!:\TTT^\ZZZ \to X\; :
\;  u(\f) \!=\!\!\!\sum_{ \nu \in \ZZZ^\ZZZ_f} e^{\ii \nu\cdot\f} u_\nu , \quad
\|u\|_{\gots,\star,X}\!:=\!\!\!\sum_{\nu\in\ZZZ^\ZZZ_f}|u_\nu|_Xe^{\gots|\nu|_\star} < + \io \Biggr\} .
\end{equation}
\end{defi}
%%%%%%%%%%%%%%%%%%%%%%%%%%%%%%%%%%%%%%%%%%%%%%%%%%%%%%%%%%%%%%%%%%%%%%%%% 

Then, we consider the infinite-dimensional Hamiltonian
\begin{equation}\label{ham}
	\calH(\theta,I)= \frac{1}{2}I\cdot I + \e f(\theta) ,\qquad a\cdot b:= \sum_{j\in\ZZZ} a_j b_j \,,
\end{equation}  
where $I=\{I_j\}_{j\in \ZZZ}$ is a sequence in
\[
\ell^{q,\io}(\RRR)=\ell^{q,\io}(\ZZZ,\RRR):= \biggl\{ I=\{I_j\}_{j\in \ZZZ}\ :\ I_j\in \RRR , \ \sup_{j\in\ZZZ} |I_j|\jap{j}^q< + \io\biggr\}\,,\
\] 
with $q>1/2$, so that $\ell^{q,\io}(\ZZZ,\RRR)\subset \ell^2(\RRR)$, while 
the \emph{perturbation} $f$ is in $\matH_\star^{2s}(\TTT^\ZZZ,\RRR)$ for some $s>0$.
In particular, $f$ admits a Fourier series representation such that
\begin{equation}\label{fourierf}
f(\theta) = \sum_{\nu\in\ZZZ^\ZZZ_f} e^{\ii \nu\cdot \theta} f_\nu, \qquad
f_{-\nu}=\ol{f_{\nu}} .
\qquad
\| f \|_{2s,\star,\RRR} < + \io .
\end{equation}

We want to show that the last condition in \eqref{normalfa} ensures that the series \eqref{fourierf} defines a real analytic function.
Since we are dealing with a function of infinitely many angles, the notion of analyticity needs some clarification.
Analytic functions on complex Banach spaces are widely studied; following Mujica \cite{Muj},
we say that a function $f$ is analytic (or holomorphic) in a complex domain
$\matA$ if, restricted to any ball $\matU\subset \matA$,
$f$ is the uniform limit of polynomials. Functions of infinitely many angles,
with values in some Banach space $X$,  are defined as the uniform limits of functions of finitely many angles,
so, by definition, $f$ is a real analytic function of infinitely many angles if it admits a
holomorphic extension which is the uniform limit of trigonometric polynomials.
Such functions are studied in details in ref.~\cite{MP}; we state here a few results, 
by referring to Section \ref{sellettina} for the proofs.

If the $\star$-norm is chosen as in Definition \ref{normstar},
requiring $f$ to belong to $\matH_\star^{2s}(\TTT^\ZZZ,\RRR)$ provides a rather strong regularity condition.
This, together the strong topology on the torus $\TTT^\ZZZ$, yields (Lemma \ref{cau})
that the functions satisfying \eqref{fourierf} are $C^\io$ and admit a
holomorphic and uniformly bounded extension to the thickened torus
\begin{equation} \nonumber 
\TTT^\ZZZ_{2s}:=\biggl\{ z\in \CCC^\ZZZ :  {\rm Re}(z)\in \TTT^\ZZZ , \;  \sup_{j\in \ZZZ} | h_j \,{\rm Im}z_j|<2s\biggr\} ,
\end{equation}
and allows us to define, for any $p\ge 0$, the $p$-th order differential
\begin{equation} \nonumber 
{\mathtt d}^p_\theta f[z_1,\dots,z_p] := \sum_{ \nu \in \ZZZ^\ZZZ_f} \ii^p 
\prod_{r=1}^p (\nu\cdot z_r) e^{\ii \nu\cdot\theta} f_\nu  \qquad \forall \,  z_1,\dots,z_p\in \ell^\io(\ZZZ) .
\end{equation}
Next, we prove (again see Lemma \ref{cau}) that the $\ell^2$-gradient of $f$ equals $\del_\theta f= \{\del_{\theta_j} f\}_{j\in\ZZZ}$
and satisfies the Cauchy estimates in $\ell^\io(\RRR)$. Finally we show
that the average
\[
\fint_{\TTT^\ZZZ} f(\theta) \, \der\theta:= \lim_{d\to +\infty} \frac{1}{(2\p)^{2d+1}}
\int_{\TTT^{2d+1}} f(\theta) \, \der \theta_{-d}\dots \der\theta_0\dots \der \theta_d  
\]
is well defined and equals $f_0$. The advantage, in using the definitions above, is that
all the usual properties of analytic functions on $\TTT^d$ still hold and one can essentially ignore the dimension;
we refer to  ref.~\cite{loom} for more abstract formulations.

%%%%%%%%%%%%%%%%%%%%%%%%%%%%%%%%%%%%%%%%%%%%%%%%%%%%%%%%%%%%%%%%%%%%%%%%%% 
\begin{rmk}
\emph{
The functional setting introduced above should be  compared with the one used by Chierchia and Perfetti \cite{CP}:
when discussing a Hamiltonian like \eqref{ham} the authors impose much weaker conditions
on both the topology of the torus (product topology) and the analyticity of $f$
(uniformity of $h_j$ in $j$ or even $h_j\to 0$),
but this is compensated by a very strong request on the frequency vector $\omega$
(super-exponential growth of $\om_j$).
}
\end{rmk}
%%%%%%%%%%%%%%%%%%%%%%%%%%%%%%%%%%%%%%%%%%%%%%%%%%%%%%%%%%%%%%%%%%%%%%%%%% 

Finally, we have to define a probability measure on $\ell^{q,\io}(\RRR)$.

%%%%%%%%%%%%%%%%%%%%%%%%%%%%%%%%%%%%%%%%%%%%%%%%%%%%%%%%%%%%%%%%%%%%%%%%% 
\begin{defi} \label{probmeas}
If $\matU_\rho(\ell^{q,\io}(\RRR))$ denotes the open ball of radius $\rho$ centered at zero in $\ell^{q,\io}(\RRR)$, then,
for any set 
\[
 A = \prod_{j\in\ZZZ} A_j\subseteq  \ol{\matU}_\rho(\ell^{q,\io}(\RRR)) = \prod_{j\in\ZZZ} [-\jap{j}^{-q}\rho,\jap{j}^{-q}\rho] ,
\]
we define the probability measure of $A$ as
\begin{equation} \nonumber 
\meas(A) := \lim_{d\to+\io} \prod_{j=-d}^{d} \meas (A_j) ,\qquad \meas(A_j)=\frac{m(A_j)}{2\jap{j}^{-q}\rho},
\end{equation}
with $m(\cdot)$ denoting the Lebesgue measure.
\end{defi}
%%%%%%%%%%%%%%%%%%%%%%%%%%%%%%%%%%%%%%%%%%%%%%%%%%%%%%%%%%%%%%%%%%%%%%%%% 

With the definitions above in hand, given the Hamiltonian \eqref{ham}
with $f$ satisfying \eqref{fourierf}, the associated dynamical system is described by 
\begin{equation}\label{sistema}
	\left\{
	\begin{aligned}
		&\dot \theta=  I , \\
		&\dot I = -\e\del_\theta f(\theta),
	\end{aligned}
	\right.
\end{equation}
or, in the Lagrangian form, by\footnote{This is the model which is studied by Chierchia and Perfetti \cite{CP}.}
\begin{equation}\label{lagra}
	\ddot \theta = -\e \del_\theta f(\theta).
\end{equation}

A solution to \eqref{sistema} is a $C^1$ map $[-T,T] \ni t \to (\theta(t),I(t)) \TTT^\ZZZ\times \ell^{q,\io}(\RRR)$,
and one may readily verify that \eqref{sistema} 
is at least locally well posed. An invariant torus of frequency $\oo\in \ell^{q,\io}(\RRR)$ for the Hamiltonian system \eqref{sistema}
is a map at least of class $C^1$ of the form
\[
\TTT^\ZZZ\ni\f \mapsto (\f + u(\f),\calI(\f))\in \TTT^\ZZZ\times \ell^{q,\io}(\RRR) \,,
\]
such that $(\oo t+ u(\oo t),\calI(\oo t))$ is a global solution to \eqref{sistema}.
Note that, since the map
$
\RRR \ni t \mapsto \om t \in \TTT^\ZZZ
$
is smooth, if $(\f+u(\f),\calI(\f))$ solves \eqref{lagrapdetot}, then $(\theta(t),I(t))=(\oo t + u(\oo t),\calI(\oo t))$
is automatically a solution to \eqref{sistema}.
We look for an analytic map $u\!:\TTT^\ZZZ\to\ell^{\io}(\RRR)$ in $\matH_\star^{s}(\TTT^\ZZZ, \ell^\io(\RRR))$ which solves
\begin{subequations}\label{lagrapdetot}
\begin{align}
	&(\om\cdot\del_\f)^2 u = -\e\del_\theta f(\f+u), \phantom{\sum)}
	\label{lagrapde} \\
	&\calI(\f):= \omega + \omega\cdot\del_\f u(\f) \in\ell^{q,\io}(\RRR). \phantom{sum)}
	\label{condi2}
\end{align}
\end{subequations}
We call an invariant torus with frequency vector $\om$ a solution $(U(\f),\calI(\f))$ to \eqref{lagrapdetot}.

Our main result is the following.

%%%%%%%%%%%%%%%%%%%%%%%%%%%%%%%%%%%%%%%%%%%%%%%%%%%%%%%%%%%%%%%%%%%%%%%%% 
\begin{theo}\label{main}
Let the $\star$-norm satisfy the further condition
\begin{equation} \label{normalfabis}
\limsup_{j\to+\io} \frac{(\log(1+\jap{j}))^{\s}}{h_j}< + \io \qquad \mbox{for some}\  \s>2 .
\end{equation}
Then, for any $s>0$ and any $f\in \matH_\star^{2s}(\TTT^\ZZZ, \RRR)$,
there exists $\e^\star=\e^\star(\s,s,\|f\|_{2s,\star,\RRR})$ such that,
for any $\e\in(-\e^\star,\e^\star)$ and any ball $\matU\subset \ell^{q,\io}(\RRR)$, the Hamiltonian \eqref{ham}
admits an invariant torus in $\matH_\star^{s}(\TTT^\ZZZ, \ell^\io(\RRR))$ 
with frequency vector $\om$ for any $\om$ in a set $\Omega_{\e}\subset \matU$
such that
\[
\lim_{\e\to 0} \meas(\matU\setminus\Omega_{\e}) = 0 .
\]
\end{theo}
%%%%%%%%%%%%%%%%%%%%%%%%%%%%%%%%%%%%%%%%%%%%%%%%%%%%%%%%%%%%%%%%%%%%%%%%% 

%%%%%%%%%%%%%%%%%%%%%%%%%%%%%%%%%%%%%%%%%%%%%%%%%%%%%%%%%%%%%%%%%%%%%%%%%% 
\begin{rmk} \label{classico}
\emph{
If there is $j_0>0$ such that $h_{j}=+\io$ for all $|j|>j_0$ in \eqref{normalfa} (see Remark \ref{hj-infty1}),
Theorem \ref{main} above implies the classical KAM theorem in dimension $d=2j_0+1$.
}
\end{rmk}
%%%%%%%%%%%%%%%%%%%%%%%%%%%%%%%%%%%%%%%%%%%%%%%%%%%%%%%%%%%%%%%%%%%%%%%%% 

To characterize the set $\Omega_{\e}$ containing the frequency vectors for which the invariant tori persist,
we introduce the Bryuno vectors associated with the $\star$-norm.

%%%%%%%%%%%%%%%%%%%%%%%%%%%%%%%%%%%%%%%%%%%%%%%%%%%%%%%%%%%%%%%%%%%%%%%%% 
\begin{defi} 
\label{br}
For $\om\in \RRR^\ZZZ$, consider the non-increasing sequence $\be^\star_\om=\{\be^\star_\om (m)\}_{m\ge0}$ given by
\begin{equation}\label{beta}
\beta^\star_\om(m):= \inf_{\substack{\nu\in\ZZZ^{\ZZZ}_{f} \vspace{-.05cm}\\ 0 <|\nu|_{\star}\le 2^m}}  |\om\cdot\nu| , 
\qquad\qquad
\om\cdot\nu := \sum_{i\in\ZZZ} \om_i \nu_i .
\end{equation}
Define the \emph{$\star$-Bryuno function} as
\begin{equation} \nonumber 
\BB^\star(\om) := \sum_{m\ge1}\frac{1}{2^m}\log \left(\frac{1}{\be^\star_\om(m)}\right)
\end{equation}
and set
\[
\gotB^\star := \{\om\in \RRR^\ZZZ : \BB^\star(\om) < + \io\} .
\]
We say that $\om$ satisfies the
\emph{Bryuno condition w.r.t.~the $\star$-norm}
or, equivalently, that $\om$ is a 
\emph{$\star$-Bryuno vector}, if  $\om\in\gotB^\star$.
\end{defi}
%%%%%%%%%%%%%%%%%%%%%%%%%%%%%%%%%%%%%%%%%%%%%%%%%%%%%%%%%%%%%%%%%%%%%%%%% 

The main step to prove Theorem \ref{main} is the following result, proved in Sections \ref{brutale} to \ref{R&C}.

%%%%%%%%%%%%%%%%%%%%%%%%%%%%%%%%%%%%%%%%%%%%%%%%%%%%%%%%%%%%%%%%%%%%%%%%% 
\begin{prop}\label{esiste1}
For any $\star$-Bryuno vector $\oo\in\ell^{q,\io}(\RRR)$ and any $f\in\matH_\star^{2s}(\TTT^\ZZZ,\RRR)$,
there exists $\e_1=\e_1(s,\om,\|f\|_{s,\star,\RRR})>0$
such that for all $\e\in(-\e_1,\e_1)$ there exists a solution
$u\in \matH_\star^{s}(\TTT^\ZZZ,\ell^\io(\RRR))$ to \eqref{lagrapde}
such that $\|u\|_{s,\star,\RRR} \to 0$ as $\e\to 0$.
\end{prop}
%%%%%%%%%%%%%%%%%%%%%%%%%%%%%%%%%%%%%%%%%%%%%%%%%%%%%%%%%%%%%%%%%%%%%%%%% 

The next step consists in checking that, if we strengthen the assumptions on the $\star$-norm,
then, once \eqref{lagrapde} is solved, \eqref{condi2} is satisfied as well.
This is the content of the following result, proved in Section \ref{sellettina}.

%%%%%%%%%%%%%%%%%%%%%%%%%%%%%%%%%%%%%%%%%%%%%%%%%%%%%%%%%%%%%%%%%%%%%%%%% 
\begin{prop}\label{esiste2}
Assume that the $\star$-norm satisfies the further condition
\begin{equation} \label{normalfater}
\lim_{j\to+\io} \frac{(\log(1+\jap{j}))}{h_j}=0. 
\end{equation}
For any $\star$-Bryuno vector $\oo\in\ell^{q,\io}(\RRR)$ and any $f\in\matH_\star^{2s}(\TTT^\ZZZ,\RRR)$
the following holds.
Let $\e_1$ and the solution $u$ to \eqref{lagrapde} be defined as in Proposition \ref{esiste1}.
Then, for any $s_2\in(0,s)$,
there exists $\e_2=\e_2(s,s_2,\om,\|f\|_{s,\star,\RRR}) \in (0,\e_1)$
such that, for all $\e\in(-\e_2,\e_2)$, the function $u$ satisfies \eqref{condi2} as well.
The value of $\e_2$ is fixed so that $\|u\|_{s,\star,\ell^\io(\RRR)} \le (s-s_2)/2$ for all $\e\in(0,\e_2)$.
\end{prop}
%%%%%%%%%%%%%%%%%%%%%%%%%%%%%%%%%%%%%%%%%%%%%%%%%%%%%%%%%%%%%%%%%%%%%%%%% 

In the light of Propositions \ref{esiste1} and \ref{esiste2},
it is natural to wonder whether the set of  $\star$-Bryuno vectors is of positive probability measure in 
any ball in $\ell^{q,\io}(\RRR)$.  
A possibile way to secure such a property is to show that $\gotB^\star$ contains vectors satisfying a suitable Diophantine condition
which ensure that they form a set of positive probability measure.
In this paper we consider the set
\begin{equation} \label{diofa}
\DD_{\g,\mu_1,\mu_2}
:=\bigg\{  \omega \in \RRR^\ZZZ: \, | \omega\cdot\nu|\ge \g \prod_{j\in \ZZZ}(1+\jap{j}^{\mu_1} |\nu_j|^{\mu_2})^{-1}
\quad \forall\, \nu\in\ZZZ^\ZZZ_f\setminus\{0\}\bigg\} ,
\end{equation}
for some positive constants $\mu_1>1+q$ and $\mu_2>1$, which was proved by Bourgain \cite{Bjfa}
to have probability measure $1-O(\g)$ in any ball in  $\ell^{q,\io}(\RRR)$.

At the end of Section \ref{sellettina} we prove the following result.

%%%%%%%%%%%%%%%%%%%%%%%%%%%%%%%%%%%%%%%%%%%%%%%%%%%%%%%%%%%%%%%%%%%%%%%%% 
\begin{prop}\label{brubru}
Let the $\star$-norm satisfy \eqref{normalfabis}.
Then, for any $\mu_1>1+q$ and $\mu_2>1$, one has  $\DD_{\g,\mu_1, \mu_2}\subseteq\gotB^\star$.
Moreover there exists a sequence $\be_*=\{\be_*(m)\}_{m\in\NNN}$, depending only on $\g,$ $\m_1$, $\m_2$ and $\s$, such that,
for all $\om\in \DD_{\g,\mu_1,\mu_2}$,
one has $\be_\om^\star(m) \ge \be_*(m)$ for all $m\in\NNN$; in particular,
the $\star$-Bryuno function $\BB^\star(\om)$ is bounded uniformly, for $\om\in \DD_{\g,\mu_1,\mu_2}$, as
\[
\BB^\star(\om) \le \BB_* := \sum_{m\ge 1} \frac{1}{2^m} \log \frac{1}{\be_*(m)} < + \io .
\]
More generally, for any sequence $\be = \{\be(m)\}_{m\in\NNN}$ such that
\[
\be(m) \le \be_*(m) \quad \forall m \ge 1 , \qquad 
\BB := \sum_{m\ge 1} \frac{1}{2^m} \log \frac{1}{\be(m)} < + \io ,
\]
there exists a set $\DD \subseteq \gotB^\star$ such that $\DD_{\g,\mu_1,\mu_2}\subseteq \DD$ and $\BB^\star(\om) \le \BB$ for all $\om\in\DD$.
\end{prop}
%%%%%%%%%%%%%%%%%%%%%%%%%%%%%%%%%%%%%%%%%%%%%%%%%%%%%%%%%%%%%%%%%%%%%%%%% 

The last ingredient we need to obtain Theorem \ref{main} is the following result, proved in Section \ref{10}.

%%%%%%%%%%%%%%%%%%%%%%%%%%%%%%%%%%%%%%%%%%%%%%%%%%%%%%%%%%%%%%%%%%%%%%%%% 
\begin{prop}\label{bieco}
Under the hypotheses of Proposition \ref{brubru}, for $\m_1>1+q$ and $\m_2>1$, one has,
for a suitable $n_0\in\NNN$,
\begin{equation} \label{e*}
\e_* :=	\inf_{\om\in \DD_{\g,\m_1, \m_2}}
\e_2(s,s_2,\om,\|f\|_{s,\star,\RRR}) \ge
c_*\biggl( \frac{s \be_*(m_{n_0})}{\|f\|_{2s,\star,\RRR}} \biggr)^4 ,
\end{equation}
where $\e_2$ is defined in Proposition \ref{esiste2},
$\be_*(m_{n_0})$ is proportional to $\g$, and
$c_*$ is a positive constant depending only on $\s$, $s$ and $s_2$.
Analogously, one has, for a suitable $n_0^\star\ge n_0$,
\begin{equation} \label{estar}
\e_\star :=	\inf_{\om\in \DD}
\e_2(s,s_2,\om,\|f\|_{s,\star,\RRR}) \ge
c_* \biggl( \frac{s \be(m_{n_0^\star})}{\|f\|_{2s,\star,\RRR}} \biggr)^4 ,
\end{equation}
with $c_*$ as in \eqref{e*}. In particular, one has $\e_*\ge \e_\star >0$.
\end{prop}
%%%%%%%%%%%%%%%%%%%%%%%%%%%%%%%%%%%%%%%%%%%%%%%%%%%%%%%%%%%%%%%%%%%%%%%%% 

The proof of Theorem \ref{main} follows from Propositions \ref{esiste1} to \ref{bieco}. 
Choosing the $\star$-norm so as to satisfy \eqref{normalfabis}, and combining Propositions \ref{esiste1} and \ref{esiste2}
we obtain, for all  $\star$-Bryuno vectors $\omega\in \ell^{q,\io}(\RRR)$ and all $\e\in(-\e_2,\e_2)$,
the existence of an invariant torus with frequency vector $\om$.
We can fix arbitrarily $s_2$, say $s_2=s/2$, to avoid using too many parameters. 
By Proposition \ref{bieco}, for all $\om\in \DD_{\g,\m_1, \m_2}$ one has
$\e_2 \ge \e_*$, with $\e_*$ depending on $\om$ only through the sequence $\{\be_*(m)\}_{m\in\NNN}$,
which in turn depends only on the constants $\g$, $\mu_1$, $\mu_2$ and $\s$.
In particular, as $\g\to 0$ one has $\be^*(m_{n_0})\to 0^+$, and hence
$\e_*\to 0^+$ as well.
Let us write $\e_*=\e_*(\g)$ to stress the dependence of $\e_*$ on $\g$.
For any $\e>0$ small enough we can fix $\g=\g(\e)$ such that $\g(\e)$ goes to zero as $\e\to 0$
and the corresponding $\e_*(\g(\e))$ is larger than $\e$.
Next, define $\Omega_{\e}$ as the set of frequency vectors $\om\in \gotB^\star \cap \ell^{q,\io}(\RRR)$ for which $\e_2 \ge \e$.
Since, by Propositions \ref{brubru} and \ref{bieco}, we have
\[
\min_{\om \in \Omega_\e} \e_2 \le \e_\star \le \e_* ,
\] 
the set $\Omega_{\e}$ contains the set of Diophantine vectors  \eqref{diofa}
and hence has probability measure $1-O(\g(\e))$.
The fact that $\Omega_{\e}$ has asymptotically full measure follows from the fact that $\g(\e)$ goes to zero as $\e\to 0$.
The bound provided by Proposition \ref{bieco} gives $\e_*(\g) \sim \g^4$; however, by simple scaling arguments
one can show that $\e_*(\g)\sim \g^2$ (see Remark \ref{finale} at the end of Section \ref{10}).

%%%%%%%%%%%%%%%%%%%%%%%%%%%%%%%%%%%%%%%%%%%%%%%%%%%%%%%%%%%%%%%%%%%%%%%%% 
%%%%%%%%%%%%%%%%%%%%%%%%%%%%%%%%%%%%%%%%%%%%%%%%%%%%%%%%%%%%%%%%%%%%%%%%% 
\zerarcounters
%%%%%%%%%%%%%%%%%%%%%%%%%%%%%%%%%%%%%%%%%%%%%%%%%%%%%%%%%%%%%%%%%%%%%%%%% 
%%%%%%%%%%%%%%%%%%%%%%%%%%%%%%%%%%%%%%%%%%%%%%%%%%%%%%%%%%%%%%%%%%%%%%%%% 
\zerarcounters
\section{Remarks and open problems}
%%%%%%%%%%%%%%%%%%%%%%%%%%%%%%%%%%%%%%%%%%%%%%%%%%%%%%%%%%%%%%%%%%%%%%%%% 
%%%%%%%%%%%%%%%%%%%%%%%%%%%%%%%%%%%%%%%%%%%%%%%%%%%%%%%%%%%%%%%%%%%%%%%%% 

We have stated our results for Hamiltonians of the form \eqref{ham}, because we wanted to highlight the key ideas behind our
\emph{dimension-independent} approach, while minimizing the technical issues. Let us give a brief overview of possible generalizations.

In view of applications to PDEs it is very natural to study a Hamiltonian of the form
\[
\calH(\theta,I)=\omega\cdot I +\frac12 I\cdot I + \e f(\theta) , 
\]
where $\omega=\{\omega_j\}_{j\in\ZZZ}\in \RRR^\ZZZ$ does not belong to $\ell^{q,\io}(\RRR)$ and possibly  $\omega_j\to + \infty$
 {as $j\to\pm\io$.} If one works at the level of the second order equation \eqref{lagra},  and looks for an invariant torus in 
$\matH_\star^{s}(\TTT^\ZZZ,\ell^\io(\RRR))$  with frequency vector $\omega$, then our algorithm can be applied word by word. 
In particular, 
in the construction we describe in the forthcoming sections,
we never use any upper bound on $\omega_j$,
and in fact our definition of $\star$-Bryuno vectors is in $\RRR^\ZZZ$.

In this context, one may consider a frequency vector $\om$ with a strong (super-exponential) divergence (say $\om_j\sim e^{j^2}$) 
in order to impose the very strong Diophantine conditions used by Chierchia and Perfetti~\cite{CP}, and thus weaken the request on $\{h_j\}_{j\in\ZZZ}$. 
Pushing this argument to full extent is beyond the scope of the present paper, but we mention that obtaining a condition weaker than \eqref{normalfabis}
in this case is straightforward.  However, as pointed out in ref.~\cite{CP},
the solutions to \eqref{lagra} found in this case
are not classical solutions because the map $t\mapsto \omega t$ is not continuous
in the topology we consider in Section \ref{functional}
if the sequence $\{\omega_j\}_{j\in\ZZZ}$ is not bounded.
Furthermore it is very well possible that such solutions have infinite energy, if one does not make
very strong assumptions on the regularity of $f$.

Another possible extension is to study whether our results hold also for the more general case
\begin{equation}
	\label{Hamaz}
	\calH(\theta,I)=\om\cdot I + \frac 12 I\cdot I + \e f(\theta,I) \, ,
\end{equation}  
where $f\!:\TTT^\ZZZ\times \matU_\rho(\ell^{q,\io}(\RRR))\to\RRR$ is real analytic for some $\rho>0$.
Up to technical difficulties making the proof more involved,
we expect a result like Theorem \ref{main} to hold 
if one further requires that  $\del_I f:\TTT^\ZZZ\times \matU_{\rho}(\ell^{q,\io}(\RRR)) \to \ell^{q,\io}(\RRR)$ is real analytic as well.
Note that the latter condition is in general not a consequence of the analitycity of $f$. 
In particular, as it happens in our model (see Lemma \ref{micnondorme}), it turns out that the regularity of $I$ is stronger that what is needed,
since we find $\sup_j |I_j e^{s_1 h_j}|<+\io$ for some\footnote{In Lemma \ref{micnondorme} we prove only $s_1<s$
due to the fact that, merely to simplify the exposition, we arbitrarily fix the loss of regularity to be $s$.} $s_1<2s$,
if one assumes $f(I,\cdot)\in \matH_\star^{2s}(\TTT^\ZZZ,\RRR)$.

It is natural to wonder whether a model such as \eqref{Hamaz} has applications to Hamiltonian PDEs. 
Unfortunately, there are some serious difficulties, which can be illustrated as follows.
Consider for instance a NLS equation on the circle:
\begin{equation} \label{nls}
\ii u_t-u_{xx}  - |u|^2u + \mbox{higher order terms}=0 \,, \qquad x\in \TTT\,.
\end{equation}
This is an infinite-dimensional Hamiltonian  system and, if one passes to the Fourier representation
$u(t,x)=\sum_{j\in \ZZZ} u_j(t)e^{\ii j x}$, performs one step of Birkhoff Normal Form, and a phase shift~\cite{KP,PP},
then the corresponding Hamiltonian has the form
\begin{equation} \label{ricompare}
\sum_{j\in\ZZZ} j^2 |u_j|^2 + \frac12 \sum_{j\in\ZZZ} |u_j|^4 \quad - \!\!\!\!\!\!\!\!\!\!
\sum_{\substack{j_1,j_2,j_3,j_4,j_5,j_6 \in \ZZZ \\j_1-j_2+j_3-j_4+j_5-j_6=0}}
h(j_1,\dots,j_6) \, u_{j_1}\bar u_{j_2} u_{j_3}\bar u_{j_4}u_{j_5}\bar u_{j_6} ,
\end{equation}
for a suitable weight $h(j_1,\dots,j_6)$, plus higher order terms. For simplicity, let us set $h(j_1,\dots,j_6)=1$ and ignore the higher order terms. 
Passing to\footnote{This is an analytic change of variables provided one imposes some lower bounds on the sequence $\{I_j\}_{j\in\ZZZ}$
and assumes that $u=\{u_j\}_{j\in\ZZZ}$ lives in an $\ell^\io$-based Banach space.}
the linear action-angle variables  $u_j= \sqrt{I_j}e^{\ii \theta_j}$ one gets the Hamiltonian \eqref{Hamaz} 
with $\omega_j= j^2$ and 
\begin{equation} \label{fItheta}
f(\theta,I)=  - \!\!\!\!\!\!\!\!\!\! \sum_{\substack{j_1,j_2,j_3,j_4,j_5,j_6 \in \ZZZ \\ j_1-j_2+j_3-j_4+j_5-j_6=0}} \!\!\!\!\!\!\!\!\!\!
e^{\theta_{j_1}-\theta_{j_2}+\theta_{j_3}-\theta_{j_4}+\theta_{j_5}-\theta_{j_6}} \, 
\sqrt{I_{j_1}I_{j_2}I_{j_3}I_{j_4}I_{j_5}I_{j_6}} .
\end{equation}
The analiticity in $\theta$ of \eqref{fItheta} is {driven} by the regularity in $I$; for instance if we only require that $I\in \ell^{q,\io}(\RRR)$ then we obtain
$h_j< \log(1+\jap{j})$, which does not satisfy \eqref{normalfater}. One might assume $I$ to be more regular, but even in such a case,
unavoidably, the regularity in the angles would turn out to be
weaker than the one in the actions and our method would not work.
We believe that this is not a purely technical problem,
also because there are no results on almost-periodic solutions for PDEs
which use action-angle variables.

Indeed, it is a general fact that,
as  pointed out by Bourgain \cite{Bjfa}, when dealing with  Hamiltonian PDEs, it is better not to resort to action-angle variables.
In the finite-dimensional setting,  there are KAM results without action-angle variables 
on the existence of maximal tori in Hamiltonian systems close to an elliptic fixed point; see for instance refs.~\cite{SZ,llave,CGP0,sto}.
For a Hamiltonian like \eqref{ricompare},with $j^2\rightsquigarrow j^2+V_j$, 
where $V=\{V_j\}_{j\in \ZZZ}\in \ell^\io(\RRR)$ is a set of external parameters, there are a number of results;
we mention, among others, refs.~\cite{Bjfa,BMP1,Yuan, Cong}.
In this context, the r\^ole of the analiticity in $\theta$ is taken by the regularity of $u$, 
namely the decay of the sequence $\{u_j\}_{j\in\ZZZ}$;
the natural condition is something like $\sup_j |u_j e^{s h_j}|<+\io$.
To study the existence of maximal tori for such a system,
one has essentially the same difficulties that we have to face in this paper.
In fact, to the best of our knowledge,
the slowest growth for the sequence $\{h_j\}_{j\in\ZZZ}$ in the literature is \eqref{normalfabis},
obtained by Cong \cite{cong2}
for the quintic NLS with external parameters,
and the only result on existence of almost-periodic solutions with only Sobolev regularity
(i.e.~with $h_j\sim \log(1+\jap{j})$)
is the result by Biasco, Massetti and Procesi
for special, non-maximal infinite-dimensional tori \cite{BMP2}.
Whether maximal tori of Sobolev regularity persist is a completely open problem.

A weak point of the results for PDEs described above is that, in order to impose
non-resonance conditions, whether Diophantine or Bryuno,
on the frequency,  one needs the set of external parameters $V\in \ell^\io(\RRR)$.
This means that the existence of almost-periodic solutions 
can be proved  only for ``many'' $V$. This  amounts to  considering families of NLS equations  \eqref{nls}, 
with $-\del_{xx}\rightsquigarrow  -\del_{xx} + T_V$, where $T_V$ is a Fourier multiplier defined as $T_Ve^{\ii j x}:= V_j e^{\ii j x}$.
The physically more relevant case $V\in \ell^{N,\io}(\RRR)$, with $N\in \NNN$,
was studied by Corsi, Gentile and Procesi~\cite{CGP}.
However, the case of $V$ fixed -- for instance to zero -- is still wide open.

Since the Hamiltonian \eqref{ricompare} has the twist term 
\[
\sum_{j\in\ZZZ} |u_j|^4 =  I\cdot I ,
\]
it would be natural to use the internal actions as parameters.
This is obtained by making a translation $I_j= V_j + y_j$, where $V_j$ are parameters and $y_j$  are the new actions.
Unfortunately there is a serious gap in this reasoning:
to obtain meaningful measure estimates on the set of parameters
for which the existence of almost-periodic solutions can be proved,
one needs the sequence of $V_j\sim I_j$ not to be too regular.
Therefore, if the one hand the existence of the maximal tori can be proved to hold for most values of the parameters
only if one assumes at most $V\in \ell^{N,\io}(\RRR)$, for some finite $N$,
one the other hand, for the perturbation to be as regular as needed,
one has to require the regularity of $|u_j|\sim \sqrt{I_j}$ to be strictly stronger than Sobolev.
In conclusion, in order to modulate both $u$ and $V$ through the same parameters one should either lower the regularity of $u$
or increase the one of $V$. This does not seem at all a simple question to handle.

%%%%%%%%%%%%%%%%%%%%%%%%%%%%%%%%%%%%%%%%%%%%%%%%%%%%%%%%%%%%%%%%%%%%%%%%%% 
%%%%%%%%%%%%%%%%%%%%%%%%%%%%%%%%%%%%%%%%%%%%%%%%%%%%%%%%%%%%%%%%%%%%%%%%%% 
\zerarcounters 
\section{On the $\star$-norm} 
\label{sellettina} 
%%%%%%%%%%%%%%%%%%%%%%%%%%%%%%%%%%%%%%%%%%%%%%%%%%%%%%%%%%%%%%%%%%%%%%%%%% 
%%%%%%%%%%%%%%%%%%%%%%%%%%%%%%%%%%%%%%%%%%%%%%%%%%%%%%%%%%%%%%%%%%%%%%%%%% 

In this section we prove some properties of the spaces $\matH_\star^{\gots}(\TTT^\ZZZ,X)$, in particular
that the conditions \eqref{normalfa} on the $\star$-norm ensure that $f$ is real-analytic.
Then, we use such results to prove Propositions \ref{esiste2} and \ref{brubru},
by showing that \eqref{condi2} is automatically satisfied if \eqref{normalfater} holds and $u$ solves \eqref{lagrapde},
and that Diophantine vectors of the form \eqref{diofa} are $\star$-Bryuno if \eqref{normalfabis} holds.

As anticipated in Section \ref{functional}, we want to prove that $f$ admits an extension to a complex thickened torus,
where it is uniform limit of trigonometric polynomials,and  that its Fourier representation \eqref{fourierf} of its gradient
is totally convergent and, together with its $p$-differentials, satisfies the Cauchy estimates.

For any $\nu\in\ZZZ^\ZZZ_f$, we define $\nu^{p+1}$ as the symmetric $p$-linear tensor
\[
\nu^{p+1}:\CCC^\ZZZ\times\ldots\times\CCC^\ZZZ\longrightarrow \CCC^\ZZZ ,
\]
whose action is given by
\begin{equation}\label{tensore0}
\begin{aligned}
\nu^{p+1}[z_1,\ldots,z_p] &:= \nu \prod_{r=1}^p (\nu\cdot z_r) ,
\end{aligned}
\end{equation}
and set
\begin{equation}\label{normaop}
\|\nu^{p+1}\|_{\rm op} := \sup_{\|z_1\|,\ldots,\|z_p\|=1} \|\nu^{p+1}[z_1,\ldots,z_p]\|= \|\nu\| \, \|\nu\|_1^p.
\end{equation}
%

%%%%%%%%%%%%%%%%%%%%%%%%%%%%%%%%%%%%%%%%%%%%%%%%%%%%%%%%%%%%%%%%%%%%%%%%%% 
\begin{rmk}\label{tensore} 
\emph{
Since $\nu\in\ZZZ^\ZZZ_f$, the sum in \eqref{tensore0} is finite, so that the tensor $\nu^{p+1}$ is well defined and it is 
continuous when acting on $\ell^\io(\CCC)\times\ldots\times\ell^{\io}(\CCC)$.
If $f$ is the perturbation in \eqref{ham}, then $\del_\theta f$ represents the action-component of the Hamiltonian vector field.
Since in \eqref{sistema} we work at the level of the Hamiltonian vector field, in fact we need
$\del_\theta f$ to be a real analytic map from $\TTT^\ZZZ\to \ell^\io(\RRR)$. Thus,
we systematically consider the $(p+1)$-differential in $\CCC$ as a $p$-linear tensor with values in $\ell^{\io}(\CCC)$. 
}
\end{rmk}
%%%%%%%%%%%%%%%%%%%%%%%%%%%%%%%%%%%%%%%%%%%%%%%%%%%%%%%%%%%%%%%%%%%%%%%%%% 

%%%%%%%%%%%%%%%%%%%%%%%%%%%%%%%%%%%%%%%%%%%%%%%%%%%%%%%%%%%%%%%%%%%%%%%%%% 
\begin{lemma}\label{cau}
Let the $\star$-norm and $\gots>0$ be given.
Then, any function $f\in \matH_\star^{\gots}(\TTT^\ZZZ,\RRR)$ is real analytic with a uniformly bounded 
holomorphic extension to $\TTT^\ZZZ_{\gots}$.
Moreover, for all $s'\in(0,\gots)$, one has $\del_\theta f\in \matH_\star^{s'}(\TTT^\ZZZ, \ell^\io(\RRR))$, and the following Cauchy estimates hold: 
\[
\sup_{\|z_1\|,\ldots,\|z_{p}\|=1} \|{\mathtt d}^{p}_\theta \del_\theta f[z_1,\dots,z_{p}]\|_{s',\star,\ell^\io(\RRR)} 
\le p! \|f\|_{\gots,\star,\RRR} \biggl(\frac{C_2^*}{\gots-s'}\biggr)^{p+1}  ,
\]
with $C^*_2$ a constant depending only on the $\star$-norm.
\end{lemma}
%%%%%%%%%%%%%%%%%%%%%%%%%%%%%%%%%%%%%%%%%%%%%%%%%%%%%%%%%%%%%%%%%%%%%%%%%% 

%%%%%%%%%%%%%%%%%%%%%%%%%%%%%%%%%%%%%%%%%%%%%%%%%%%%%%%%%%%%%%%%%%%%%%%%%% 
\prova
The fact that the Fourier series defining $f$ is totally convergent w.r.t~the sup norm in  $\TTT^\ZZZ_{\gots}$ comes directly from  
the assumption $f\in \matH_\star^{\gots}(\TTT^\ZZZ,\RRR)$
 and the definition \eqref{normalfa}, since
\[
\sum_{\nu\in\ZZZ^\ZZZ_{f}} \sup_{\theta\in \TTT^\ZZZ_{\gots}} |f_\nu e^{\ii \nu\cdot\theta}| \le
\sum_{\nu\in\ZZZ^\ZZZ_{f}} |f_\nu| e^{\gots|\nu|_\star} =\|f\|_{\gots,\star,\RRR}<\io\,.
\]
In turn, the total convergence of the serie implies that $f$ can be obtained as the uniform limit of trigonometric polynomials.

If we define
\[
\tilde {\mathtt d}^{p}_\theta f[z_1,\dots,z_{p}](\theta) := \sum_{\nu\in\ZZZ^\ZZZ_{f}} f_\nu \,
{\mathtt d}^{p}_\theta \, e^{\ii \theta\cdot \nu} [z_1,\dots,z_{p}]= \ii^p \sum_{\nu\in\ZZZ^\ZZZ_{f}} f_\nu \, e^{\ii \theta\cdot \nu} 
z_{p}\cdot \nu^{p}[z_1,\ldots,z_{p-1}] ,
\]
then, for any $s'\in(0,\gots)$ and all $p\ge 0$, we find
\[
\sum_{ \nu \in \ZZZ^\ZZZ_f}  \sup_{\|z_1\|,\ldots,\|z_{p}\|=1}\sup_{\theta\in \TTT^\ZZZ_{s'}}
|\prod_{r=1}^p (\nu\cdot z_r) e^{\ii \nu\cdot\theta} f_\nu |\le  \sup_{\nu\in\ZZZ^\ZZZ_{f}} \|\nu\|^p_1 e^{-(\gots-s')|\nu|_\star}  \sum_{ \nu \in \ZZZ^\ZZZ_f} 
e^{\gots |\nu|_\star}|f_\nu | .
\]
where, for any $\de>0$, we may bound
\begin{equation} \label{jedoernome}
\sup_{\nu\in\ZZZ^\ZZZ_{f}} \|\nu\|^p_1 e^{-\delta|\nu|_\star}  \le 
c^p_1 \sup_{\nu\in\ZZZ^\ZZZ_{f}} |\nu|^p_\star e^{-\delta|\nu|_\star}\le p! c_1^p \delta^{-p} ,
\end{equation}
with the constant $c_1>0$ depending on the sequence $\{h_j\}_{j\in\ZZZ}$. 
The bound \eqref{jedoernome} implies the total convergence of $\tilde{\mathtt d}^{p}_\theta f[z_1,\dots,z_{p}]$
and allows us to set ${\mathtt d}^{p}_\theta f[z_1,\dots,z_{p}] = \tilde{\mathtt d}^{p}_\theta f[z_1,\dots,z_{p}]$.

To show that $\ell^2$-gradient of $f$ is in $\matH_\star^{s'}(\TTT^\ZZZ,\ell^\io(\RRR))$, for any $s' \in(0,\gots)$, we note that
\[
\|\del_\theta f\|_{s',\star,\ell^\io(\RRR)} = \sum_{ \nu \in \ZZZ^\ZZZ_f}\|\nu\| | f_\nu| e^{s'|\nu|_\star} \le c_1\sum_{\nu\in\ZZZ^\ZZZ_{f}} 
| f_\nu| e^{\gots |\nu|_\star}  \sup_{ \nu \in \ZZZ^\ZZZ_f}|\nu|_{\star} e^{-(\gots-s') |\nu|_\star}<\io ,
\]
so that we deduce that the $\ell^2$-gradient of $f$ is well defined and equal $\partial_\theta f$.
Finally, as to the Cauchy estimates, we have 
\[
\begin{aligned} 
\sup_{\|z_1\|,\ldots,\|z_{p}\|=1} \|{\mathtt d}^{p}_\theta \del_\theta f[z_1,\dots,z_{p}]\|_{s',\star,\ell^\io}
%&
&\le \sum_{\nu\in\ZZZ^\ZZZ_{f}} | f_\nu| e^{\gots |\nu|_\star} 
\sup_{ \nu \in \ZZZ^\ZZZ_f}  \|\nu\|\|\nu^{p}\|_{1} e^{-(\gots-s') |\nu|_\star} \\
& 
\le (p+1)! \|f\|_{\gots,\star,\RRR} c_1^{p+1} (\gots-s')^{-(p+1)} .
\end{aligned}
\]
This completes the proof and provides the value $C_2^*=2c_1$.
\EP	
%%%%%%%%%%%%%%%%%%%%%%%%%%%%%%%%%%%%%%%%%%%%%%%%%%%%%%%%%%%%%%%%%%%%%%%%%% 

%%%%%%%%%%%%%%%%%%%%%%%%%%%%%%%%%%%%%%%%%%%%%%%%%%%%%%%%%%%%%%%%%%%%%%%%%% 
\begin{rmk}
\emph{
In the proof of Lemma \ref{cau}, we only used that $\| \cdot \| \le \| \cdot \|_1 \le c_1  | \cdot |_\star$ for some positive constant $c_1$.
}
\end{rmk}
%%%%%%%%%%%%%%%%%%%%%%%%%%%%%%%%%%%%%%%%%%%%%%%%%%%%%%%%%%%%%%%%%%%%%%%%%% 

Now, we prove two abstract results, which are discussed in ref.~\cite{MP} in the special case $h_j=\jap{j}^\al$, and
which are fundamental in the iterative procedure discussed in Section \ref{brutale}.

%%%%%%%%%%%%%%%%%%%%%%%%%%%%%%%%%%%%%%%%%%%%%%%%%%%%%%%%%%%%%%%%%%%%%%%%%% 
\begin{lemma}\label{facciofou}
Consider any function $u\in {\matH}^\gots(\TTT^\ZZZ, X)$, with $\gots>0$.
Then, for all $\nu \in\ZZZ^\ZZZ_f$ one has
\begin{equation}\label{ansia1}
u_\nu= \fint_{\TTT^\ZZZ} u(\f) e^{- \ii \nu \cdot \f}\, \der \f := \lim_{d \to \infty} \frac{1}{(2\pi)^{2d+1}}\int_{\TTT^{2d+1}} u(\f) e^{-\ii \nu\cdot\f} \,
\der \f_{-d} \ldots \der \f_{d} .
\end{equation}
\end{lemma}
%%%%%%%%%%%%%%%%%%%%%%%%%%%%%%%%%%%%%%%%%%%%%%%%%%%%%%%%%%%%%%%%%%%%%%%%%% 
%%%%%%%%%%%%%%%%%%%%%%%%%%%%%%%%%%%%%%%%%%%%%%%%%%%%%%%%%%%%%%%%%%%%%%%%%% 
\prova
Since $h_j\to \io$ monotonically, for any $d \ge1$ there exists  $N(d)$ such that if
$\nu \in \ZZZ^\ZZZ_f$ satisfies  $|\nu|_\star <N(d)$ then $\nu_j=0$ for all $j$ such that $|j|>d$. 
Thus, for such $\nu$'s one has
\[
e^{\ii \nu \cdot \f} = e^{\ii \nu_{-d} \f_{-d}} \ldots e^{\ii \nu_{d} \f_{d}} ,
\]
implying that 
\[
\frac{1}{(2 \pi)^{2d+1}} \int_{\TTT^{2d+1}} e^{\ii \nu \cdot \f}\, \der \f_{-d} \ldots \der \f_{d} = 0\,.
\]
Hence
\[
\begin{aligned}
& \frac{1}{(2 \pi)^{2d +1}}\int_{\TTT^{2d +1}}u(\f) \, \der \f_{-d} \ldots \der \f_d \\
& \qquad \qquad = \frac{1}{(2 \pi)^{2d +1}}\int_{\TTT^{2d +1}} \bigg(u_0 + \sum_{0 < |\nu|_\star \leq N(d)} u_\nu e^{\ii \nu \cdot \f} 
+ \sum_{|\nu|_\star > N(d)} u_\nu e^{\ii \nu \cdot \f} \bigg) \der \f_{-d} \ldots \der \f_d \\
& \qquad \qquad = u_0 + \frac{1}{(2 \pi)^{2d +1}}\int_{\TTT^{2d +1}} \sum_{|\nu|_\star > N(d)} u_\nu e^{\ii \nu \cdot \f} \, \der \f_{-d} \ldots \der \f_d\,.
\end{aligned}
\]
Since $u \in {\matH}^\gots(\TTT^\ZZZ, X)$, the tail of the series wih $|\nu|_\star > N(d)$ goes to zero as $d \to \infty$.
This proves \eqref{ansia1} for $\nu=0$.

Now, take $\nu' \in \ZZZ^\ZZZ_f \setminus \{ 0 \}$ and set
\[
u(\nu',\f) := u(\f) e^{- \ii \nu' \cdot \f} = \sum_{k \in \ZZZ^\ZZZ_f}u_k e^{\ii (k - \nu') \cdot \f} = \sum_{h \in \ZZZ^\ZZZ_f} u_{h+\nu'} e^{\ii h \cdot \f}\,.
\]
By applying \eqref{ansia1} with $\nu=0$ to the function $u(\nu',\f)$ and using that $\nu'$ is arbitrary,
the equality \eqref{ansia1} follows for $\nu\neq0$ as well.
\EP
%%%%%%%%%%%%%%%%%%%%%%%%%%%%%%%%%%%%%%%%%%%%%%%%%%%%%%%%%%%%%%%%%%%%%%%%%% 

%%%%%%%%%%%%%%%%%%%%%%%%%%%%%%%%%%%%%%%%%%%%%%%%%%%%%%%%%%%%%%%%%%%%%%%%%% 
\begin{lemma}\label{lemmasopra}
Consider any function $f\in \matH_\star^\gots(\TTT^\ZZZ, \RRR)$, with $\gots>0$.
Then, for any $s'\in(0,\gots)$ there exists $\e_3=\e_3(\gots,s')$ such that, for any $h\in \matH_\star^{s'}(\TTT^\ZZZ,\ell^\io(\RRR))$
with $\| h\|_{s',\star,\ell^\io(\RRR)} \le \e_3$, the Taylor expansion
\begin{equation}\label{sarto}
\partial_\theta f (\f+ h(\f)) = \sum_{p\ge0} \frac{1}{p!} \mathtt{d}^p_\theta \partial_\theta f (\f)[h(\f),\dots,h(\f)]=\sum_{p\ge0} \frac{1}{p!}
\sum_{ \nu_0 \in \ZZZ^\ZZZ_f } 
\ii \nu_0 f_{\nu_0}e^{\ii \nu_0\cdot\f} \prod_{i=1}^{p} (\nu_0\cdot h(\f))
\end{equation}
is totally convergent in $\TTT^\ZZZ_{s'}$, and $\partial_\theta f(\cdot+ h(\cdot)) \in \matH_\star^{s'}(\TTT^\ZZZ,\ell^{\io}(\RRR))$.
\end{lemma}
%%%%%%%%%%%%%%%%%%%%%%%%%%%%%%%%%%%%%%%%%%%%%%%%%%%%%%%%%%%%%%%%%%%%%%%%%% 
	
%%%%%%%%%%%%%%%%%%%%%%%%%%%%%%%%%%%%%%%%%%%%%%%%%%%%%%%%%%%%%%%%%%%%%%%%%% 
\prova
The fact that the Taylor expansion \eqref{sarto} is totally convergent w.r.t.~the sup norm in $\TTT^\ZZZ_{s'}$
follows, using the Cauchy estimates and noting that 
\[
\sup_{\theta\in \TTT^\ZZZ_{s'}}\|h(\theta)\|\le \|h\|_{s',\star,\ell^\io(\RRR)} ,
\]
as soon as we require $2 C_2^* \| h\|_{\gots,\star,\ell^\io}(\gots-s')^{-1}<1$, with $C_2^*$ as in Lemma \ref{cau}.
For $h\in \matH_\star^{s'}(\TTT^\ZZZ,\ell^\io(\RRR))$, indeed, we may consider its Fourier expansion and
insert it into the r.h.s.~of \eqref{sarto}, so as to write
\begin{equation}\label{sarto1}
\partial_\theta f (\f+ h(\f)) 
:=  \sum_{ \nu \in \ZZZ^\ZZZ_f}\sum_{p\ge0} \frac{1}{p!}  e^{\ii \f\cdot\nu}
\sum_{\substack{\nu_0,\dots,\nu_p \in \ZZZ^\ZZZ_f \\  \nu_0+\ldots+\nu_r=\nu}} \ii \nu_0 f_{\nu_0} \prod_{i=1}^p (\nu_0\cdot h_{\nu_r})
\end{equation}
In order to show that $\partial_\theta f(\cdot + h(\cdot)) \in\matH_\star^{s'}(\TTT^\ZZZ,\ell^\io(\RRR))$ for $s'\in(0,\gots)$, 
we bound the r.h.s.~of \eqref{sarto1} as
\[
\sum_{p\ge0} \frac{1}{p!} \sum_{\nu_0\in\ZZZ^\ZZZ_f} e^{s'|\nu_0|_\star}  \|\nu_0 \|^{p+1} |f_{\nu_0}| \| h \|_{ s',\star,\ell^\io(\RRR)}^p
\le \sum_{p\ge0} \| f \|_{\gots} \left( 2 \| h \|_{ s',\star,\ell^\io(\RRR)} \right)^p \left(\gots - s' \right)^{-(p+1)} ,
\]
which converges for $\| h \|_{ s',\star,\ell^\io(\RRR)}$ small enough.
\EP
%%%%%%%%%%%%%%%%%%%%%%%%%%%%%%%%%%%%%%%%%%%%%%%%%%%%%%%%%%%%%%%%%%%%%%%%%% 

Let us now show that the existence of a solution $u$ to \eqref{lagrapde} implies that $I$ in \eqref{condi2}
is defined as well and is in $\ell^{q,\io}(\RRR)$, and hence ensures the existence of an invariant torus with frequency vector $\om$.

%%%%%%%%%%%%%%%%%%%%%%%%%%%%%%%%%%%%%%%%%%%%%%%%%%%%%%%%%%%%%%%%%%%%%%%%%% 
\begin{lemma}\label{micnondorme}
For any $\gots>0$ define the scale of Banach spaces
\begin{equation} \nonumber 
\MM^\gots:= \biggl\{v\in\CCC^\ZZZ\, :\, |v|_{\MM^\gots}:=\sup_{j\in\ZZZ} |e^{\gots h_j} v_j| < + \io \biggr\}.
\end{equation}
Let $\e_1$ and the solution $u$ to \eqref{lagrapde} be as in Proposition \ref{esiste1}.
Then, for any $s_1,s_2\in(0,s)$
 {
there exists $\e_2 \in (0,\e_1)$ such that for all $\e\in(0,\e_2)$
one has $(\om\cdot\del_\f)u \in \matH_\star^{s_1}(\TTT^\ZZZ,\MM^{s_2})$.
The value of $\e_2$ is fixed by requiring that
$\|u\|_{s,\star,\ell^\io(\RRR)} \le (s-s_2)/2$ for all $\e\in(0,\e_2)$.
}
\end{lemma}
%%%%%%%%%%%%%%%%%%%%%%%%%%%%%%%%%%%%%%%%%%%%%%%%%%%%%%%%%%%%%%%%%%%%%%%%%% 

%%%%%%%%%%%%%%%%%%%%%%%%%%%%%%%%%%%%%%%%%%%%%%%%%%%%%%%%%%%%%%%%%%%%%%%%%% 
\prova
By \eqref{lagrapde} and Lemmata \ref{facciofou} and \ref{lemmasopra}, for any $\nu\in\ZZZ^\ZZZ_f\setminus\{0\}$ we have
\[
(\om\cdot\nu) \, u_\nu = \frac{\e}{\om\cdot\nu} \sum_{p\ge0} \frac{1}{p!} 
\sum_{\substack{\nu_0,\ldots, \nu_p\in \ZZZ^\ZZZ_f \setminus\{0\} \\ 
\nu_0+\ldots +\nu_p=\nu}}  \ii \nu_0 f_{\nu_0}  \prod_{r=1}^p (\ii\nu_0 \cdot u_{\nu_r}) ,
\]
so that, setting $\de_1:=s-s_1$ and
\[
C(\de_1):=\sup_{\nu\in\ZZZ^\ZZZ_f\setminus\{0\}}\frac{1}{|\om\cdot\nu|}e^{-\de_1|\nu|_\star},
\]
which is finite if $\om$ is $\star$-Bryuno, we obtain
\[
\|(\om\cdot\del_\f)u\|_{s_1,\star,\MM^{s_2}}
\le |\e| C(\de_1) \sum_{p\ge0} \frac{1}{p!} 
\sum_{\nu\in\ZZZ^\ZZZ_f} e^{s|\nu|_\star} \!\!\!\!\!\!\!\!
 \sum_{\substack{\nu_0,\ldots, \nu_p\in \ZZZ^\ZZZ_f \setminus\{0\} \\ 
 \nu_0+\ldots +\nu_p=\nu}} |f_{\nu_0} | \sup_{j\in\ZZZ}|e^{s_2 h_j} (\nu_0)_j| \prod_{r=1}^p |\nu_0 \cdot u_{\nu_r}| .
\vspace{-.2cm}
\]

Thus, using that
\[
\sup_{j\in\ZZZ}|e^{s_2h_j} (\nu_0)_j| \le \|\nu_0\|e^{s_2|\nu_0|_\star},\qquad |\nu_0 \cdot u_{\nu_r}| \le \|\nu_0\|_1 \|u_{\nu_r}\| ,
\vspace{-.2cm}
\]
we obtain
\[
\begin{aligned}
\|(\om\cdot\del_\f) u\|_{s_1,\star,\MM^{s_2}}
& \le |\e| C(\de_1) \sum_{p\ge0} \frac{1}{p!} 
\sum_{\nu\in\ZZZ^\ZZZ_f} e^{s|\nu|_\star}  \!\!\!\!\!\!\!\!
 \sum_{\substack{\nu_0,\ldots, \nu_p\in \ZZZ^\ZZZ_f \setminus\{0\} \\ 
 \nu_0+\ldots +\nu_p=\nu}} |f_{\nu_0} | \|\nu_0\|e^{s_2|\nu_0|_\star}\prod_{r=1}^p  \|\nu_0\|_1 \|u_{\nu_r}\| \\
&\le |\e| C(\de_1) \sum_{p\ge0} \frac{1}{p!} 
\sum_{\nu\in\ZZZ^\ZZZ_f} 
\sum_{\substack{\nu_0,\ldots, \nu_p\in \ZZZ^\ZZZ_f \setminus\{0\} \\ 
\nu_0+\ldots +\nu_p=\nu}} |f_{\nu_0} \|\nu_0\|^p e^{(s+s_2)|\nu_0|_\star}\prod_{r=1}^p  
e^{s |\nu_r|_\star} \|u_{\nu_r}\|  \\
&\le  |\e| C(\de_1) \sum_{p\ge0} \frac{1}{p!} 
\sum_{\nu_0\in\ZZZ^\ZZZ_f} e^{(s+s_2)|\nu_0|_\star} |f_{\nu_0}| 
|\nu_0|_\star^{p+1}
\|u\|_{s,\star,\ell^\io(\RRR)}^p  ,  
\end{aligned}
\]
where, for any $\de_2>0$, we can bound
\[
\sum_{p\ge0} \frac{1}{p!}  |\nu_0|_\star^{p+1} \|u\|_{s,\star,\ell^\io(\RRR)}^p 
\le e^{\de_2 |\nu_0|_\star}
\sum_{p \ge 0} \biggl( \frac{2\|u\|_{s,\star,\ell^\io(\RRR)}}{\de_2} \biggr)^{p+1} .
\]
Finally, requiring $s_2 < s$ and choosing $\de_2= s-s_2$, we obtain
\begin{equation} \label{sump}
\|(\om\cdot\del_\f) u\|_{s_1,\star,\MM^{s_2}} \le
|\e| C(s-s_1) \, \| f \|_{2s,\star,\RRR} \sum_{p\ge 0} \biggl( \frac{2\|u\|_{s,\star,\ell^\io(\RRR)}}{s-s_2} \biggr)^{p+1} ,
\end{equation}
so that the assertion follows provided $\e$ is small enough for the series in \eqref{sump} to converge,
 {
i.e.~for $2\|u\|_{s,\star,\ell^\io(\RRR)}/(s-s_2)$ to be smaller than 1.
Since $u$ is continuous in $\e$ and goes to 0 as $\e$ tend to 0,
there exists $\e_2\in(0,\e_1)$ such that $\|u\|_{s,\star,\ell^\io(\RRR)} \le (s-s_2)/2$ for all $\e\in(0,\e_2)$.
}
\EP
%%%%%%%%%%%%%%%%%%%%%%%%%%%%%%%%%%%%%%%%%%%%%%%%%%%%%%%%%%%%%%%%%%%%%%%%%% 

\medskip

%%%%%%%%%%%%%%%%%%%%%%%%%%%%%%%%%%%%%%%%%%%%%%%%%%%%%%%%%%%%%%%%%%%%%%%%%% 
\noindent\emph{Proof of Proposition \ref{esiste2}.}
Under the assumption \eqref{normalfater} one checks easily that  $\MM^{\gots}\subset \ell^{q,\io}(\RRR)$ for any $\gots,q>0$. 
Then the assertion follows immediately from Lemma \ref{micnondorme}
 {
by fixing arbitrarily any $s_2\in(0,s)$ and
taking $\e_2$ so as to ensure that
$2\|u\|_{s,\star,\ell^\io(\RRR)}/(s-s_2)$ to be less than 1 for all $\e\in(0,\e_2)$.
}
\EP
%%%%%%%%%%%%%%%%%%%%%%%%%%%%%%%%%%%%%%%%%%%%%%%%%%%%%%%%%%%%%%%%%%%%%%%%%% 

To prove that Diophantine vectors -- in the sense of \eqref{diofa} -- are also $\star$-Bryuno we need a stronger condition on the sequence
$\{h_j\}_{j\in\ZZZ}$.

%%%%%%%%%%%%%%%%%%%%%%%%%%%%%%%%%%%%%%%%%%%%%%%%%%%%%%%%%%%%%%%%%%%%%%%%%% 
\begin{lemma}\label{taglio}
Assume $|\cdot|_\star$ to be as in \eqref{normalfabis} for some $\s>2$.
For any $\mu_1,\mu_2\in\RRR_+$ there are positive constants $K_1=K_1(\mu_1,\m_2,\s)$ and $K_2=K_2(\s)$ such that for all $N\ge 1$ one has
	\begin{equation}
		\label{tagliobis}
		\sup_{\substack{\nu\in \ZZZ^\ZZZ_f  \vspace{-.1cm}\\  |\nu |_{\star}\le N}} \prod_{j\in\ZZZ}(1+\jap{j}^{\mu_1} |\nu_j|^{\mu_2}) \le 
		K_1 e^{K_2N/(\log N)^{\s-1}}.
	\end{equation}
\end{lemma}
%%%%%%%%%%%%%%%%%%%%%%%%%%%%%%%%%%%%%%%%%%%%%%%%%%%%%%%%%%%%%%%%%%%%%%%%%% 

%%%%%%%%%%%%%%%%%%%%%%%%%%%%%%%%%%%%%%%%%%%%%%%%%%%%%%%%%%%%%%%%%%%%%%%%%% 
\prova
We start by writing
\[
\prod_{j\in\ZZZ}(1+\jap{j}^{\mu_1} |\nu_j|^{\mu_2}) \le \max\{\mu_1,\mu_2\}\exp\sum_{j\in\ZZZ} \log(1+ \jap{j}|\nu_j|)
\]
and splitting
\[ 
\sum_{j\in\ZZZ} \log(1+ \jap{j}|\nu_j|)=\sum_{\substack{j\in\ZZZ \\ |j|\le j_0}} \log(1+ \jap{j}|\nu_j|) + \sum_{\substack{j\in\ZZZ \\ |j|> j_0}} \log(1+ \jap{j}|\nu_j|) ,
\]
with $j_0= A_\s \, N/(\log N)^{\s}$, with $A_\s=\s^{\s-1}/\log \s$.
Using that $|\nu|_\star\le N$ and hence $|\nu_j|\le N$ for any $j\in\ZZZ$, the first sum is bounded as
\[
\sum_{\substack{j\in\ZZZ \\ |j|\le j_0}} \log(1+ \jap{j}|\nu_j|) \le (2j_0+1)\log(2j_0 N) \le 3 A_\s \, \frac{N}{(\log N)^{\s}} \log
\bigg( \frac{2A_\s N^2}{(\log N)^{\s}} \bigg) ,
\]
whereas the second sum is bounded as
\[
\sum_{\substack{j\in\ZZZ \\ |j|> j_0}} \log(1+ \jap{j}|\nu_j|) \le \sum_{\substack{j\in\ZZZ \\ |j|> j_0}}\frac{ \log(1+ \jap{j})^{\s}}{(\log j_0)^{\s-1}}|\nu_j|
\le \frac{N}{\log(A_\s\,N/(\log N)^{\s})^{\s-1}}.
\]
Then the assertion follows with $K_2=O(\s^{\s})$.
\EP
%%%%%%%%%%%%%%%%%%%%%%%%%%%%%%%%%%%%%%%%%%%%%%%%%%%%%%%%%%%%%%%%%%%%%%%%%% 

Proposition \ref{brubru} follows easily from Lemma \ref{taglio}.

\medskip 

%%%%%%%%%%%%%%%%%%%%%%%%%%%%%%%%%%%%%%%%%%%%%%%%%%%%%%%%%%%%%%%%%%%%%%%%%% 
\noindent\emph{Proof of Proposition \ref{brubru}.}
If we define 
\begin{equation} \label{betam}
\be_*(m)  := \frac{\g}{K_1} e^{-K_2 2^m (\log 2^m)^{-(\s-1)}} , \qquad m \in \NNN ,
\end{equation}
by Lemma \ref{taglio} {and by definition of $\be$, we have
\begin{equation} \label{boundbetam}
\frac{1}{\be^\star_\om(m)} \le \frac{1}{\be_*(m)} \le \frac{1}{\be(m)}
\end{equation}
%ù
and hence $\BB_\star(\om) \le \BB_* \le \BB$. Moreover we can bound
\[
\log \left( \frac{1}{\be^\star_\om(m)} \right) \le \log \left( \frac{1}{\be_*(m)} \right) = 
\log \frac{K_1}{\g} + C_\s 2^m m^{-(\s-1)} ,
\]
with $C_\s=K_2(\log 2)^{-(\s-1)}$. Thus, for $\s>2$, we find
\begin{equation} \label{bryuno-gamma}
\BB^\star(\om) := \sum_{m\ge1}\frac{1}{2^m}\log \left(\frac{1}{\be^\star_\om(m)}\right) 
\le \log \frac{K_1}{\g} + C_\s \sum_{m\ge1} \frac{1}{m^{\s-1}} < K_0 + \log\frac{K_1}{\g} ,
\end{equation}
for a suitable constant $K_0$, depending only on $\s$.
Since both $K_0$ and $K_1$ are independent of $\g$, 
the $\star$-Bryuno function $\BB^*(\om)$ is bounded uniformly as in \eqref{bryuno-gamma} for all $\om\in\DD_{\g,\mu_1,\mu_2}$.
\EP
%%%%%%%%%%%%%%%%%%%%%%%%%%%%%%%%%%%%%%%%%%%%%%%%%%%%%%%%%%%%%%%%%%%%%%%%%% 

%%%%%%%%%%%%%%%%%%%%%%%%%%%%%%%%%%%%%%%%%%%%%%%%%%%%%%%%%%%%%%%%%%%%%%%%%% 
\begin{rmk}
\emph{
The weaker the norm $\star$-norm, the larger is the bound in \eqref{tagliobis}.
For instance, if we take $h_j= \log(1+\jap{j})$ for all $j\in\ZZZ$ and define
\vspace{-.2cm}
\[
j_0 = \max\biggl\{ j \in \NNN : \sum_{j=-j_0}^{j_0} \log (1+\jap{j}) \le N \biggr\} , \qquad
\ZZZ_{j_0} := \{ j \in \ZZZ: |j| \le j_0\} ,
\vspace{-.2cm}
\]
then 
\vspace{-.2cm}
\[
\sup \sum_{j\in\ZZZ}\log (1+\jap{j}^{\mu_1} |\nu_j|^{\mu_2}) \quad \mbox{restricted to}\quad   \sum_{j\in\ZZZ}\log (1+\jap{j})|\nu_j| \le N
\vspace{-.2cm}
\] 
is larger that
\vspace{-.2cm}
\[
\min\{\mu_1,\mu_2\}\sup \sum_{j\in\ZZZ_{j_0}}\log (1+\jap{j} ) \quad \mbox{restricted to}\quad   \sum_{j\in\ZZZ_{j_0}}\log (1+\jap{j}) \le N ,
\vspace{-.2cm}
\]
which is trivially of order $N$, so that the Diophantine vectors in $\DD_{\g,\mu_1,\mu_2}$ are not $\star$-Bryuno if the $\star$-norm
is defined as in \eqref{normalfa} with $\s=0$.
}
\end{rmk}
%%%%%%%%%%%%%%%%%%%%%%%%%%%%%%%%%%%%%%%%%%%%%%%%%%%%%%%%%%%%%%%%%%%%%%%%%% 

%%%%%%%%%%%%%%%%%%%%%%%%%%%%%%%%%%%%%%%%%%%%%%%%%%%%%%%%%%%%%%%%%%%%%%%%%% 
%%%%%%%%%%%%%%%%%%%%%%%%%%%%%%%%%%%%%%%%%%%%%%%%%%%%%%%%%%%%%%%%%%%%%%%%%% 
\zerarcounters 
\section{The recursive equations} 
\label{brutale} 
%%%%%%%%%%%%%%%%%%%%%%%%%%%%%%%%%%%%%%%%%%%%%%%%%%%%%%%%%%%%%%%%%%%%%%%%%% 
%%%%%%%%%%%%%%%%%%%%%%%%%%%%%%%%%%%%%%%%%%%%%%%%%%%%%%%%%%%%%%%%%%%%%%%%%% 

To complete the proof of Theorem \ref{main}, we are left with the task of proving Propositions \ref{esiste1} and \ref{bieco}.

%%%%%%%%%%%%%%%%%%%%%%%%%%%%%%%%%%%%%%%%%%%%%%%%%%%%%%%%%%%%%%%%%%%%%%%%%% 
\begin{rmk} \label{content}
\emph{
The rest of the paper, from here up to Section \ref{R&C}, is mainly devoted to
the proof of Proposition \ref{esiste1} , which, from a technical point of view, represents the core of the analysis.
Proposition \ref{bieco} is easily proved in Section \ref{10} by relying on
Propositions \ref{esiste1} and \ref{esiste2} and on an explicit uniform estimate of the value of $\e_2$.
}
\end{rmk}
%%%%%%%%%%%%%%%%%%%%%%%%%%%%%%%%%%%%%%%%%%%%%%%%%%%%%%%%%%%%%%%%%%%%%%%%%% 

%%%%%%%%%%%%%%%%%%%%%%%%%%%%%%%%%%%%%%%%%%%%%%%%%%%%%%%%%%%%%%%%%%%%%%%%%% 
\begin{rmk}\label{invarianza}
\emph{
If $u\in \matH_\star^{s}(\TTT^\ZZZ,\ell^{\io}(\RRR))$ solves \eqref{lagrapde},
 then for all $\f_0\in\TTT^\ZZZ$  the function $u(\cdot+\f_0)$ solves 
\begin{equation} \label{lagrapde2}
(\om\cdot\del_\f)^2 u = -\e\del_\theta f(\f+\f_0+u) .
\end{equation}
Such a property is known as \emph{translation covariance}, and plays a key r\^ole in KAM theory
(see for instance \cite{llave} for the finite-dimensional case).
}
\end{rmk}
%%%%%%%%%%%%%%%%%%%%%%%%%%%%%%%%%%%%%%%%%%%%%%%%%%%%%%%%%%%%%%%%%%%%%%%%%% 

Decomposing
\begin{equation} \nonumber 
\matH_\star^{s}(\TTT^\ZZZ,\ell^{\io}(\RRR))=\matH_0^s(\TTT^\ZZZ,\ell^{\io}(\RRR))\oplus \matH^s_\perp(\TTT^\ZZZ,\ell^{\io}(\RRR)),
\end{equation}
with
\begin{subequations} \label{range+kernel}
\begin{align}
\matH^s_0(\TTT^\ZZZ,\ell^{\io}(\RRR)) & := {\rm Ker}((\om\cdot\del_\f)^2)= \{g\in \matH_\star^{s}(\TTT^\ZZZ,\ell^{\io}(\RRR))\, :\, g(\f) = \hbox{const} \} , \\
\matH^s_\perp(\TTT^\ZZZ,\ell^{\io}(\RRR)) & :={\rm Rank}((\om\cdot\del_\f)^2)= \biggl\{g  \in \matH_\star^{s}(\TTT^\ZZZ,\ell^{\io}(\RRR)) : % \avg{ g}:=
\fint_{\TTT^\ZZZ} g(\f) \, \der \f=0 \biggr\} ,
\end{align}
\end{subequations}
we introduce a Lyapunov-Schmidt decomposition and split \eqref{lagra} into the range and kernel equations,
given respectively by
\begin{subequations}\label{ranker}
\begin{align}
&(\om\cdot \del_\f)^2u = -\e \del_\theta f(\f+u) + \e \avg{ \del_\theta f(\cdot+u)} , \phantom{\big(}
\label{range1} \\
&  \e \avg{ \del_\theta f(\cdot+u)} =0. \phantom{\Big(}
\label{ker1}
\end{align}
\end{subequations}
%

%%%%%%%%%%%%%%%%%%%%%%%%%%%%%%%%%%%%%%%%%%%%%%%%%%%%%%%%%%%%%%%%%%%%%%%%%% 
\begin{rmk}\label{rinvarianza}
\emph{
As a consequence of Remark \ref{invarianza}, if $u\in \matH^s_\perp(\TTT^\ZZZ,\ell^{\io}(\RRR))$ solves \eqref{range1} and
$\om$ has rationally independent components, then for all $\f_0\in\TTT^\ZZZ$ the function $u(\cdot+\f_0)$ solves 
\begin{equation}\label{range2}
(\om\cdot \del_\f)^2u = -\e \del_\theta f(\f+\f_0+u) + \e \avg{ \del_\theta f(\cdot+\f_0+u )} .
\end{equation}
In the present paper we use the translation covariance property to prove Lemma \ref{metalemma} below,
which in turn not only ensures that the \emph{kernel equation} \eqref{ker1}
is automatically satisfied by any solution to the \emph{range equation} \eqref{range1},
but also implies the cancellation in Corollary \ref{cancellazione},
which plays a crucial r\^ole in the proof of the existence of the solution to \eqref{lagrapdetot}.
}
\end{rmk}
%%%%%%%%%%%%%%%%%%%%%%%%%%%%%%%%%%%%%%%%%%%%%%%%%%%%%%%%%%%%%%%%%%%%%%%%%% 

%%%%%%%%%%%%%%%%%%%%%%%%%%%%%%%%%%%%%%%%%%%%%%%%%%%%%%%%%%%%%%%%%%%%%%%%%% 
\begin{lemma}\label{metalemma}
Assume that there exists a solution $u\in \matH^s_\perp(\TTT^\ZZZ,\ell^{\io}(\RRR))$ to \eqref{range1}. Then $u$ solves \eqref{ker1} as well.
\end{lemma}
%%%%%%%%%%%%%%%%%%%%%%%%%%%%%%%%%%%%%%%%%%%%%%%%%%%%%%%%%%%%%%%%%%%%%%%%%% 

%%%%%%%%%%%%%%%%%%%%%%%%%%%%%%%%%%%%%%%%%%%%%%%%%%%%%%%%%%%%%%%%%%%%%%%%%% 
\prova
For all $\f_0\in\TTT^\ZZZ$, consider the Lagrangian action
\begin{equation}\label{lag}
\matL(u;\f_0):= \fint_{\TTT^\ZZZ} \left(
\frac{(\om\cdot\del_\f u(\f+\f_0))^2}{2}-\e f(u(\f+\f_0)+\f+\f_0)
 \right) \der\f.
\end{equation}
Under the regularity assumption on $u$, we can differentiate \eqref{lag} with respect to $\f_0$, so as to obtain
\begin{equation}\label{difff}
\begin{aligned}
\del_{\f_0}\matL(u;\f_0)&= \fint_{\TTT^\ZZZ} (\om\cdot\del_\f u (\f+\f_0)) \left( \om\cdot\del_\f\del_{\f_0}u (\f+\f_0) \right) \der \f \\
& \qquad\qquad -\fint_{\TTT^\ZZZ}\e \del_\theta f(u(\f+\f_0)+\f+\f_0)(1+\del_{\f_0}u)\der\f \\
& = -\fint_{\TTT^\ZZZ} \!\!\! \left( \left( (\om\cdot\del_\f)^2 u(\f \!+\! \f_0) +\e \del_\theta f(u(\f\!+\!\f_0)+\f\!+\!\f_0) \right) \del_{\f_0}u(\f\!+\!\f_0)\right)\der\f \\
& \qquad\qquad - \e \avg{\del_\theta f(u(\cdot + \f_0) + \cdot + \f_0)} ,
\end{aligned}
\end{equation}
where the last equality is obtained integrating by parts. Thanks to Remark \ref{rinvarianza},
inserting \eqref{range2} into \eqref{difff} gives
\begin{equation} \nonumber 
\begin{aligned}
\del_{\f_0}\matL(u;\f_0)
& = \fint_{\TTT^\ZZZ} \e\avg{\del_\theta f (u(\cdot+\f_0)+\cdot+\f_0)}\del_{\f_0}u (\f+\f_0)\, \der\f - \e \avg{\del_\theta f(u(\cdot + \f_0) + \cdot + \f_0)}   \\
&= \e \avg{\del_\theta f(u(\cdot + \f_0) + \cdot + \f_0)} 
\fint_{\TTT^\ZZZ} \del_{\f_0}u (\f+\f_0) \,\der \f - \e \avg{\del_\theta f(u(\cdot + \f_0) + \cdot + \f_0)} \\
& = - \e \avg{\del_\theta f(u(\cdot + \f_0) + \cdot + \f_0)} , \phantom{\fint_{\TTT^\ZZZ}}
\end{aligned}
\end{equation}
where, in the last equality we have used that, since $u\in \matH^s_\perp(\TTT^\ZZZ,\ell^{\io}(\RRR))$,
$$
\fint_{\TTT^\ZZZ} \del_{\f_0}u(\f+\f_0)\,d\f = \del_{\f_0}\fint_{\TTT^\ZZZ} u(\f+\f_0)\,\der\f= 0.
$$
Thus, we arrive at
\[
\begin{aligned}
\e\avg{\del_\theta f(\cdot+u(\cdot)) } &= \fint_{\TTT^\ZZZ} \e\del_\theta f(\f+u(\f)) \, \der\f 
=\fint_{\TTT^\ZZZ} \e\del_\theta f(\f+\f_0+u(\f+\f_0)) \, \der\f \\
& =-\del_{\f_0}\matL(u;\f_0) =-\del_{\f_0}\fint_{\TTT^\ZZZ} \left(
\frac{(\om\cdot\del_\f u(\f+\f_0))^2}{2}-\e f(u(\f+\f_0)+\f+\f_0)
 \right) \! \der\f\\
 &=-\del_{\f_0}\fint_{\TTT^\ZZZ} \left(
\frac{(\om\cdot\del_\f u(\f))^2}{2}-\e f(u(\f)+\f)
 \right) \! \der\f=0,
\end{aligned}
\]
so the assertion follows.
\EP
%%%%%%%%%%%%%%%%%%%%%%%%%%%%%%%%%%%%%%%%%%%%%%%%%%%%%%%%%%%%%%%%%%%%%%%%%% 

In view of Remark \ref{rinvarianza} and Lemma \ref{metalemma}, we look for a solution to \eqref{range2} of the form $u(\cdot+\f_0)$, such that
$u \in \matH^s_\perp(\TTT^\ZZZ,\ell^{\io}(\RRR))$ where $\f_0\in\TTT^\ZZZ$ is an arbitrary parameter.

First of all, we write $u(\f+\f_0)$ as a formal Lindstedt expansion, i.e.~we write
\begin{equation}\label{solou}
u(\f) = \sum_{k\ge 1}\e^k u^{(k)}(\f) =
\sum_{k\ge 1}\e^k \sum_{\nu\in\ZZZ^\ZZZ_f\setminus\{0\}} e^{\ii \nu\cdot\f} u^{(k)}_\nu ,
\end{equation}
and hence
\begin{equation}\label{uu}
u(\f + \f_0) =\sum_{k\ge 1}\e^k \sum_{\nu\in\ZZZ^\ZZZ_f\setminus\{0\}} e^{\ii \nu\cdot\f} u^{(k)}_\nu(\f_0),
\qquad u^{(k)}_\nu(\f_0):=e^{\ii\nu\cdot\f_0}u^{(k)}_\nu, 
\end{equation}
so that, plugging \eqref{uu} into \eqref{range2}, we obtain, for any $k\ge 1$, 
\begin{subequations}\label{ricorreancora}
\begin{align}
&(\om\cdot \nu)^2 u^{(k)}_\nu(\f_0) = [\del f(\f+\f_0+u(\f + \f_0)]^{(k-1)}_\nu,\qquad \nu\ne 0, \phantom{\big(}
\label{range} \\
& [\del f(\f+\f_0+u(\f+\f_0))]^{(k-1)}_0=0 , \phantom{\Big(}
\label{ker}
\end{align}
\end{subequations}
where, for all $k\ge 1$ and all $\nu\in\ZZZ^\ZZZ_f$, we have set
\begin{equation} \label{espando}
\begin{aligned}
& [\del f(\f+\f_0+u(\f+\f_0))]^{(k-1)}_\nu  \\
& \qquad\qquad := \sum_{p\ge1} \sum_{\substack{\nu_0,\ldots, \nu_p\in \ZZZ^\ZZZ_f \setminus\{0\} 
\\ \nu_0+\ldots +\nu_p=\nu}}
\sum_{\substack{k_1,\ldots,k_p \ge 1 \\ k_1+\ldots + k_p = k-1}} \frac{1}{p!} f_{\nu_0}(\f_0)(\ii \nu_0)^{p+1} 
[u^{(k_1)}_{\nu_1}(\f_0),\ldots, u^{(k_p)}_{\nu_p}(\f_0)]\,,
\end{aligned}
\end{equation}
with the $p$-linear tensor $(\ii\nu_0)^{p+1}$ defined as in \eqref{tensore0} and
\begin{equation}\label{fase}
f_{\nu}(\f_0):= f_\nu e^{\ii\nu\cdot\f_0} .
\end{equation}

If the sum in \eqref{uu} is absolutely convergent, then, if $u(\f+\f_0)$ solves \eqref{range},
it satisfies \eqref{ker} as well, as a consequence of Lemma \ref{metalemma}.
Therefore we can confine ourselves to look for a solution of the recursive equations \eqref{range}.
One can easily check, also by recursion, that
\begin{equation}\label{reale}
u^{(k)}_{-\nu} = \ol{u^{(k)}_\nu}
\end{equation}
which ensure $u^{(k)}(\f+\f_0)$ to be real-valued.

In order to prove the absolute convergence of the series \eqref{uu}, we introduce 
a convenient graphical representation for the coefficients $u^{(k)}_\nu$; this allows us to obtain
a bound 
\[
 {
\sum_{\nu\in\ZZZ^\ZZZ_f}\| u^{(k)}_\nu\| e^{s|\nu|_\star} \le C^k ,
}
\] 
which turns out to be enough to conclude the proof of Proposition \ref{esiste1}.

%%%%%%%%%%%%%%%%%%%%%%%%%%%%%%%%%%%%%%%%%%%%%%%%%%%%%%%%%%%%%%%%%%%%%%%%%% 
%%%%%%%%%%%%%%%%%%%%%%%%%%%%%%%%%%%%%%%%%%%%%%%%%%%%%%%%%%%%%%%%%%%%%%%%%% 
\zerarcounters 
\section{Tree expansion} 
\label{alberi} 
%%%%%%%%%%%%%%%%%%%%%%%%%%%%%%%%%%%%%%%%%%%%%%%%%%%%%%%%%%%%%%%%%%%%%%%%%% 
%%%%%%%%%%%%%%%%%%%%%%%%%%%%%%%%%%%%%%%%%%%%%%%%%%%%%%%%%%%%%%%%%%%%%%%%%% 

A graph $G$ is a pair $G=(V,L)$, where $V=V(G)$ is a non-empty set of elements, called \emph{vertices},
and $L=L(G)$ a family of couples of unordered elements of $V(G)$, called \emph{lines}.
A \emph{planar} graph is a graph $G$ which can be drawn on a plane without lines crossing.
Two vertices $v,w \in V(G)$ are said to be connected if either $(v,w)\in L(G)$ or there exist
$v_{0}=v,v_{1},\ldots,v_{n-1},v_{n}=w\in V(G)$ such that, defining the \emph{path}
$\calP(v,w):=\{(v_{0},v_{1}),\ldots,(v_{n-1},v_{n})\}$, then $\calP(v,w) \subset L(G)$.
A graph $G$ is \emph{connected} if for any $v,w\in V(G)$ there exists a path of lines connecting them,
and has a \emph{loop} if there exists $v\in V(G)$ such that either $(v,v)\in L(G)$ or there exists a path $\calP$ connecting $v$ to itself.
A \emph{subgraph} $S$ of a graph $G$ is a graph $S=(V,L)$ such that $V=V(S)\subseteq V(G)$ and $L=L(S)\subseteq L(G)$.

If the graph $G$ is \emph{oriented}, i.e.~if each couple in $L(G)$ is ordered,
we say that the line $\ell=(v,w)\in L(G)$ is oriented as well, and that $\ell$ \emph{exits} $v$ and \emph{enters} $w$.
An oriented line $\ell=(v,w)$ is drawn with an arrow from $v$ to $w$ superimposed.

Any subgraph $S$ of an oriented graph $G$ inherits the orientation of $G$. We say that $\ell=(v,w)\in L(G)$
\emph{exits} a subgraph $S$ of $G$ if $v\in V(S)$ and $w\in V(G)\setminus V(S)$, and that $\ell'=(v',w')$ \emph{enters} $S$
if $v'\in V(G)\setminus V(S)$, while $w'\in V(S)$. The lines entering and exiting $S$ are called the \emph{external lines} of $S$,
while the lines $\ell\in L(S)$ are the \emph{internal lines} of $S$. If $S$ has only one exiting line $\ell$
and only one entering line $\ell'$, we set $\ell_S=\ell$ and $\ell_S'=\ell'$.

A \emph{rooted tree graph} $\vartheta$ is a planar connected oriented graph
with no loops and with a special vertex $r$, called the \emph{root},
such that only one line enters $r$ and for any $v\in V(\vartheta)\setminus\{r\}$ either $(v,r)\in L(\vartheta)$ or there are $v_1,\ldots,v_n\in V(\vartheta)$
such that $\calP(v,r)=\{(v,v_1),(v_1,v_2),\ldots ,(v_n,r)\}\subseteq L(\vartheta)$.
We call \emph{root line} the line $\ell_{\vartheta}$ entering the root $r$ and $v_{\vartheta}$ the node that the root line exits.
We set also $N(\vartheta):=V(\vartheta)\setminus\{r\}$ and call \emph{nodes} the elements of $N(\vartheta)$.
Given any subgraph $S \subseteq \vartheta$, let $k(S)$ denote the number of nodes of $S$;
we say that $k(S)$ is the \emph{order} of $S$.

The orientation on a rooted tree graph $\vartheta$ provides a partial ordering relation on the vertices:
given $v,w\in V(\vartheta)$, we write $v\preceq w$ if $\calP(v,w)\subseteq L(\vartheta)$.
Since any line $\ell \in L(\vartheta)$ is uniquely determined by the vertex $v$ which it exits,
we may write $\ell=\ell_{v}$ and say that $\ell_{v} \preceq w$ if $v \preceq w$.
An example of rooted tree graph is represented in Figure \ref{albero}: by convention, the tree is drawn so that the root is the leftmost vertex,
and the arrows of all the lines point from right to left toward the root; moreover, to distinguish further the root from the other vertices,
the root is drawn as a gray bullet, while the nodes are drawn as black bullets.

%%%%%%%%%%%%%%%%%%%%%%%%%%%%%%%%%%%%%%%%%%%%%%%%%%%%%%%%%%%%%%%%%%%%%%%% 
\begin{figure}[ht] 
\centering 
\ins{098pt}{-084.5pt}{$\vartheta=$}
\includegraphics[width=3in]{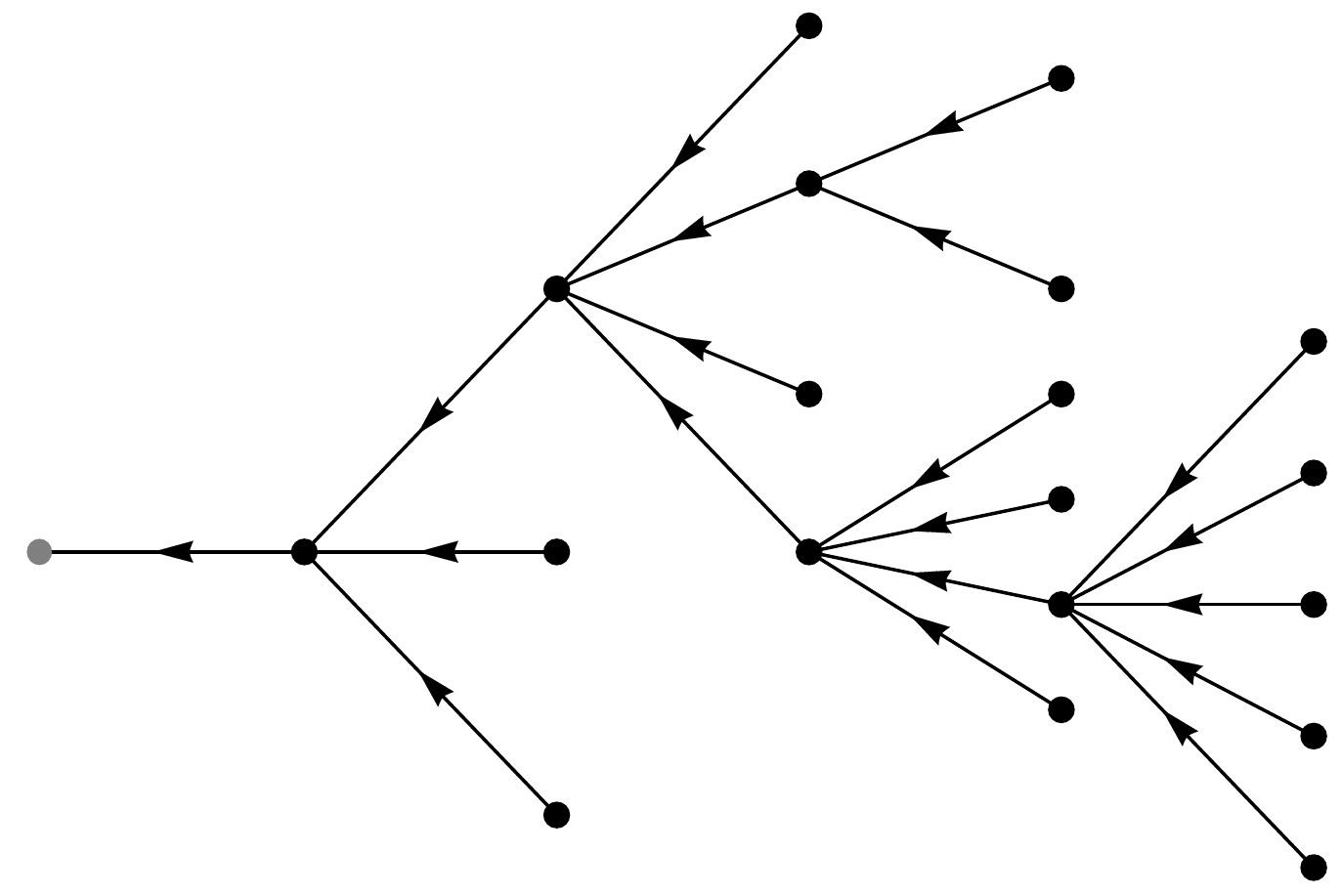}
\caption{A rooted  tree graph $\vartheta$ with $19$ nodes.}
\label{albero} 
\end{figure} 
%%%%%%%%%%%%%%%%%%%%%%%%%%%%%%%%%%%%%%%%%%%%%%%%%%%%%%%%%%%%%%%%%%%%%%%% 

A \emph{labelled rooted tree graph} -- simply \emph{tree} in the following -- is a rooted tree graph together with a label function
defined on $V(\vartheta)$ and $L(\vartheta)$.  A subgraph of a tree $\vartheta$ is called a \emph{subtree} if it is a tree.
The labels are assigned according to the following rules.
\begin{itemize}[topsep=1ex]
\itemsep0em
\item
With each node $v$ we associate a \emph{mode} label $\nu_v\in\ZZZ^\ZZZ_f$ and a \emph{branching} label $p_v \ge 0$
which denotes the number of lines entering $v$. We then associate with $v$ a \emph{node factor}
\begin{equation}\label{nodo}
\calF_v:= \frac{1}{p_v!} f_{\nu_v} (\ii \nu_v)^{p_v+1} .
\end{equation}
Note that \eqref{nodo} is a $p_v$-linear tensor in $\ell^\io(\CCC)$ and reduces to an element of $\ell^\io(\CCC)$ for $p_v=0$.
\item
With each line $\ell$ we associate a \emph{scale} label $n_\ell\in\NNN$ and
a \emph{momentum} label $\nu_\ell\in\ZZZ^\ZZZ_f$, with the constraints that
$\nu_\ell=0$ is allowed only for the root line and each momentum $\nu_\ell$ satisfies the \emph{conservation law}
\begin{equation}\label{conserva}
\nu_\ell= \sum_{\substack{v \in N(\vartheta) \\ v\preceq \ell}} \nu_v .
\end{equation}
Then, we associate with $\ell$ a \textit{propagator}\footnote{The name is borrowed from Quantum Field Theory \cite{Galla,GGM}.}
$\calG_\ell:=\calG_{n_\ell}(\om\cdot \nu_\ell)$,
where
\begin{equation}\label{prop}
\calG_{n}(x) :=\left\{
\begin{aligned}
&\frac{\Psi_{n}(x)}{x^2} , \quad &x\ne0 , \\
&1 ,  &x=0 ,
\end{aligned}
\right.
\end{equation}
with the \emph{scale function} $\Psi_n(x)$ defined as follows. Introduce the sequence $\{m_n\}_{n\in\ZZZ_+}$ such that
$m_0=0$ and $m_{n+1}=m_n+i_n+1$, where
\[
i_n:=\max\{ i\in\NNN\;:\; \be^\star_{\om}({m_n})<2 \be^\star_{\om}({m_n+i})\}.
\]
In other words, $\be^\star_{\om}({m_{n+1}})$ is the largest element in $\{\be^\star_{\om}({m})\}_{m\ge0}$
smaller than  $\frac{1}{2}\be^\star_{\om}({m_n})$. Then we set 
\begin{equation}\label{taglia}
\Psi_n(x):= 
\left\{
\begin{aligned}
&1,\qquad \frac{1}{4}\be^\star_{\om}({m_{n_\ell}}) \le |x| < \frac{1}{4}\be^\star_{\om}({m_{n_\ell-1}}), \\
&0,\qquad \mbox{ otherwise,}
\end{aligned}
\right.
\end{equation}
where $\be^\star_{\om}(m_{-1})$ has to be interpreted as $+\io$.
\end{itemize}

Finally, for any $v\in N(\vartheta)\setminus\{v_\vartheta\}$, 
if $\pi(v)$ is the node which the line exiting from $v$ enters,
we define the \emph{value} of a tree $\vartheta$ as
\begin{equation}\label{val}
\begin{aligned}
\Val(\vartheta) 
& := 
\Bigg(\prod_{v\in N(\vartheta)} \calF_v\Bigg) \Bigg(\prod_{\ell\in L(\vartheta)} \calG_\ell\Bigg) \\
& =  \frac{1}{p_{v_\vartheta}!}
\ii \nu_{v_\vartheta} f_{\nu_{v_\vartheta}} \Biggl( \prod_{v\in N(\vartheta)\setminus \{v_\vartheta\}} 
\frac{1}{p_{v}!}(\ii \nu_{\pi(v)}\cdot \ii \nu_v ) f_{\nu_v} \Biggl)
\Bigg(\prod_{\ell\in L(\vartheta)} \!\!\calG_\ell\Bigg) ,
\end{aligned}
\end{equation}
where the symbol
\[
\prod_{v\in N(\vartheta)} \calF_v
 \]
has to be interpreted as an ordered application of the $p_v$-linear tensors $\calF_v$, 
with the ordering determined by the partial ordering of the tree. This means that, if $k(\vartheta)=1$, we have
$\calF_v = \ii \nu_v f_{\nu_v}$, with $\{v\}=N(\vartheta)$, while, for $k(\vartheta)>1$, we define, recursively,
\begin{equation} \label{defrecF}
\prod_{v\in N(\vartheta)} \calF_v = \calF_{v_\vartheta}
\Biggl[ \prod_{v\in N(\vartheta_1)} \calF_v , \; \ldots, \prod_{v\in N(\vartheta_{p_{v_\vartheta}})} \calF_v \Biggr]
\end{equation}
where $\vartheta_1,\ldots,\vartheta_{p_{v_\vartheta}}$ are the subtrees whose root lines enter $v_\vartheta$. 

%%%%%%%%%%%%%%%%%%%%%%%%%%%%%%%%%%%%%%%%%%%%%%%%%%%%%%%%%%%%%%%%%%%%%%%% 
\begin{rmk} \label{lospostoqui}
\emph{
For any $v\in N(\vartheta)$ with $p_v \ge 1$, we have\
\[
\Bigg( \prod_{\substack{ w \in N(\vartheta) \\ w \preceq v }} \calF_w \Bigg)
\Bigg( \prod_{\substack{ \ell\in L(\vartheta) \\ \ell \preceq v}} \calG_\ell\Bigg) = 
\frac{1}{p_v !}  f_{\nu_v} (\ii\nu_v)^{p_v+1} [\Val(\vartheta_{v,1}),\ldots,\Val(\vartheta_{v,p_v})] ,
\]
where $\vartheta_{v,1},\ldots,\vartheta_{v,p_v}$ denote the subtrees of $\vartheta$ whose root lines enter $v$.
}
\end{rmk}
%%%%%%%%%%%%%%%%%%%%%%%%%%%%%%%%%%%%%%%%%%%%%%%%%%%%%%%%%%%%%%%%%%%%%%%% 

%%%%%%%%%%%%%%%%%%%%%%%%%%%%%%%%%%%%%%%%%%%%%%%%%%%%%%%%%%%%%%%%%%%%%%%% 
\begin{ex}\label{ordine1}
\emph{
As an example consider the trees with $4$ nodes represented in Figure \ref{4alberellini}:
for the tree on the top left one has
\[
\prod_{v\in N(\vartheta_1)} \calF_v = 
f_{\nu_1} (\ii\nu_1)^2 [ f_{\nu_2} (\ii\nu_2)^2[ f_{\nu_3} (\ii\nu_3)^2 [f_{\nu_4} \ii\nu_4]]]
=  \ii^7f_{\nu_1}f_{\nu_2}f_{\nu_3} f_{\nu_4} \nu_1 (\nu_1\cdot \nu_2)(\nu_2\cdot\nu_3) (\nu_3\cdot\nu_4) ,
\]
for the tree on the top right one has
\[
\prod_{v\in N(\vartheta_2)} \calF_v = 
\frac{1}{3!} f_{\nu_1} (\ii\nu_1)^4 [f_{\nu_2} \ii\nu_2 , f_{\nu_3} \ii\nu_3 ,  f_{\nu_4} \ii\nu_4]
= \frac{1}{3!} 
\ii^7f_{\nu_1}f_{\nu_2}f_{\nu_3} f_{\nu_4} \nu_1 (\nu_1\cdot\nu_2)(\nu_1\cdot \nu_3)(\nu_1\cdot \nu_4) .
\]
for the tree in the bottom left one has
\[
\prod_{v\in N(\vartheta_3)} \calF_v = 
f_{\nu_1} (\ii\nu_1)^2 \left[\frac{1}{2} f_{\nu_1} (\ii\nu_2)^3[  f_{\nu_3} i\nu_3 ,  f_{\nu_1} \ii\nu_4 ]\right]
= \frac{1}{2} \ii^7f_{\nu_1}f_{\nu_2}f_{\nu_3} f_{\nu_4} \nu_1 (\nu_1\cdot\nu_2)(\nu_2\cdot \nu_3)(\nu_2\cdot \nu_4) ,
\]
and, finally, for the tree in the bottom right one has
\[
\prod_{v\in N(\vartheta_4)} \calF_v = 
\frac{1}{2} f_{\nu_1} (\ii\nu_1)^3 [f_{\nu_3} \ii\nu_3 ,  f_{\nu_2} (\ii\nu_2)^2[f_{\nu_4}\ii\nu_4]]
= \frac{1}{2} \ii^7f_{\nu_1}f_{\nu_2}f_{\nu_3} f_{\nu_4}
\nu_1 (\nu_1\cdot\nu_2)(\nu_1\cdot\nu_3) (\nu_2\cdot\nu_4) ,
\]
}
\end{ex} 
%%%%%%%%%%%%%%%%%%%%%%%%%%%%%%%%%%%%%%%%%%%%%%%%%%%%%%%%%%%%%%%%%%%%%%%% 

%%%%%%%%%%%%%%%%%%%%%%%%%%%%%%%%%%%%%%%%%%%%%%%%%%%%%%%%%%%%%%%%%%%%%%%% 
\begin{rmk}\label{taglione}
\emph{
According to \eqref{val}, the definition \eqref{taglia} combined with \eqref{prop} implies that if $\Val(\vartheta)\ne0$
then, for any. line $\ell\in L(\vartheta)$, the scale label $n_\ell$ identifies uniquely the size of the \emph{small divisor} $(\om\cdot\nu_\ell)^2$
which appears in the propagator $\calG_\ell$.
}
\end{rmk}
%%%%%%%%%%%%%%%%%%%%%%%%%%%%%%%%%%%%%%%%%%%%%%%%%%%%%%%%%%%%%%%%%%%%%%%% 

%%%%%%%%%%%%%%%%%%%%%%%%%%%%%%%%%%%%%%%%%%%%%%%%%%%%%%%%%%%%%%%%%%%%%%%% 
\begin{figure}[ht] 
\vspace{-.8cm} 
\centering 
\ins{041pt}{-044.5pt}{$\vartheta_1=$}
\ins{100pt}{-054pt}{$\nu_1$}
\ins{133pt}{-054pt}{$\nu_2$}
\ins{167pt}{-054pt}{$\nu_3$}
\ins{201pt}{-054pt}{$\nu_4$}
\ins{239pt}{-44.5pt}{$\vartheta_2=$}
\ins{396pt}{-054pt}{$\nu_2$}
\ins{330pt}{-054pt}{$\nu_1$}
\ins{396pt}{-022pt}{$\nu_3$}
\ins{396pt}{-086pt}{$\nu_4$}
\null
\hspace{0.8cm}
\includegraphics[width=2.0in]{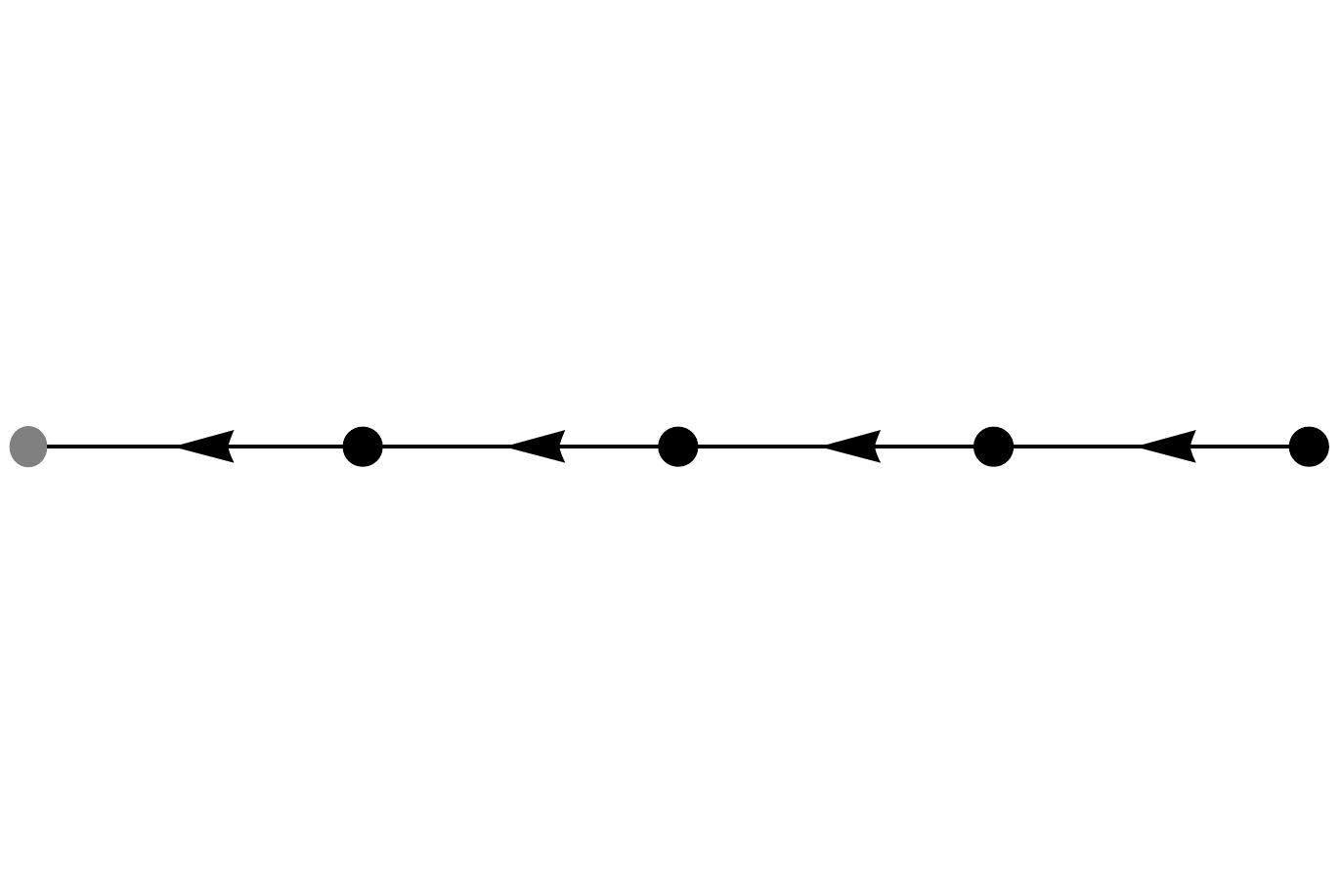}
\hspace{1.6cm} 
\includegraphics[width=2.0in]{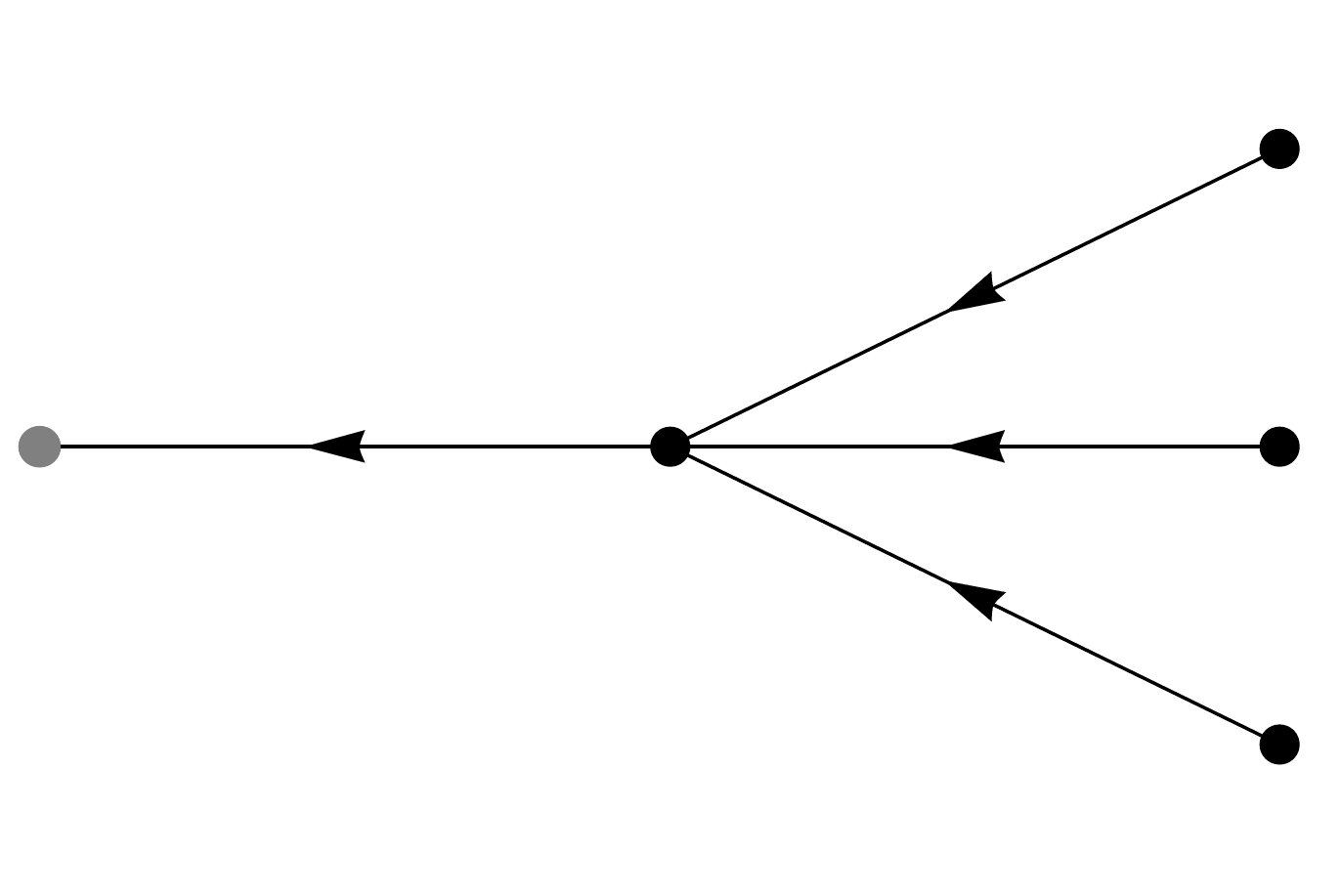}
\end{figure} 
%%%%%%%%%%%%%%%%%%%%%%%%%%%%%%%%%%%%%%%%%%%%%%%%%%%%%%%%%%%%%%%%%%%%%%%% 
%
%%%%%%%%%%%%%%%%%%%%%%%%%%%%%%%%%%%%%%%%%%%%%%%%%%%%%%%%%%%%%%%%%%%%%%%% 
\begin{figure}[H] 
\vspace{-1.2cm} 
\ins{041pt}{-044.5pt}{$\vartheta_3=$}
\ins{113pt}{-053pt}{$\nu_1$}
\ins{157pt}{-053pt}{$\nu_2$}
\ins{203pt}{-021pt}{$\nu_3$}
\ins{203pt}{-086pt}{$\nu_4$}
\ins{239pt}{-044.5pt}{$\vartheta_4=$}
\ins{308pt}{-053pt}{$\nu_1$}
\ins{355pt}{-032pt}{$\nu_3$}
\ins{355pt}{-075pt}{$\nu_2$}
\ins{400pt}{-096pt}{$\nu_4$}
\hspace{2.2cm}
\includegraphics[width=2.0in]{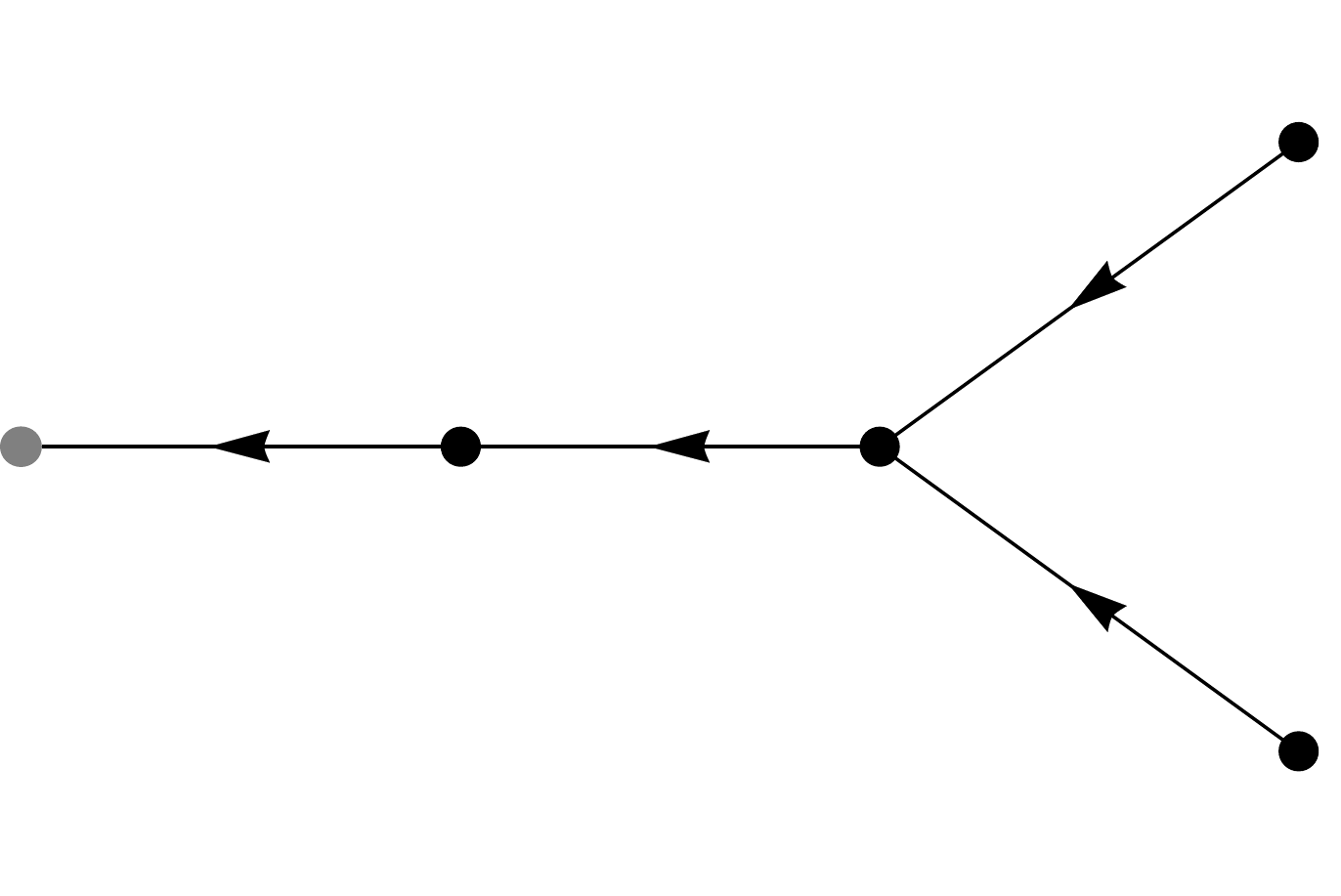}
\hspace{1.6cm} 
\includegraphics[width=2.0in]{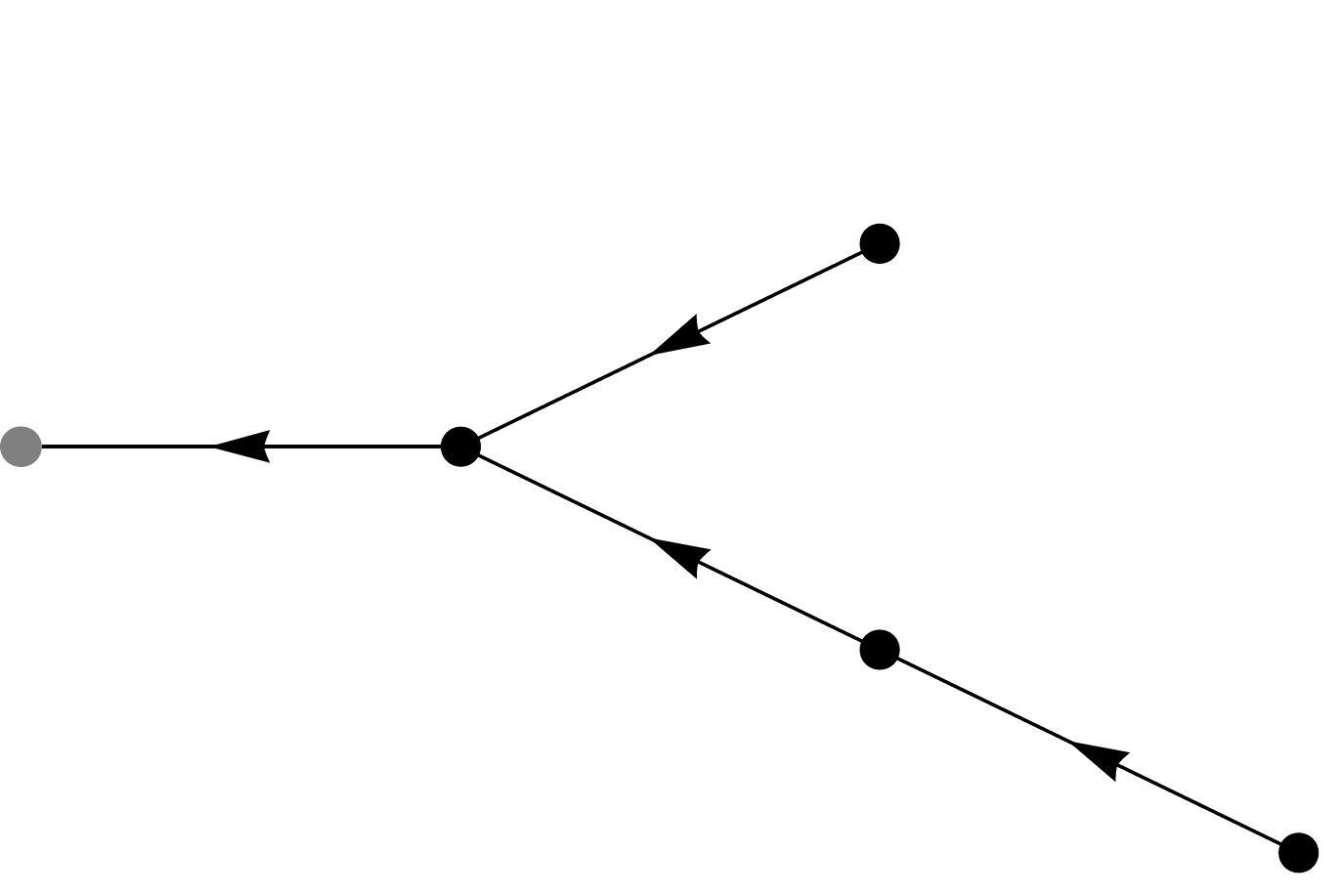}
\vskip.2truecm 
\caption{Trees with $4$ nodes (only the mode labels are showed).}
\label{4alberellini} 
\end{figure} 
%%%%%%%%%%%%%%%%%%%%%%%%%%%%%%%%%%%%%%%%%%%%%%%%%%%%%%%%%%%%%%%%%%%%%%%% 

%%%%%%%%%%%%%%%%%%%%%%%%%%%%%%%%%%%%%%%%%%%%%%%%%%%%%%%%%%%%%%%%%%%%%%%% 
\begin{rmk}\label{supporto}
\emph{
For technical reasons, it is actually more convenient to consider the scale functions to be $C^\io$ functions with compact support 
rather than step functions as in \eqref{taglia}. Indeed, it turns out to be useful to consider derivatives $\del_x\calG_n(x)$
(see the forthcoming Remark \ref{vanishing} for further comments).
}
\end{rmk}
%%%%%%%%%%%%%%%%%%%%%%%%%%%%%%%%%%%%%%%%%%%%%%%%%%%%%%%%%%%%%%%%%%%%%%%% 

We may perform the following operation on trees: we detach the line exiting a node and reattach it to another node (see Figure \ref{attach},
where the line shifted, i.e.~detached and reattached, is the root line. As a consequence of the operation, we reverse the arrows of the lines, if needed,
in such a way that all the arrows point toward the node which the line has been reattached to. For instance, in the example of
Figure \ref{attach}, after detaching the root line from the node with mode label $\nu_1$ of the tree $\vartheta$ (tree on the left) 
and reattaching to the node with mode label $\nu_2$, the tree $\vartheta'$ (tree in the middle) is obtained such that
the arrow of the line connecting the two nodes is reverted; finally, on the right the tree $\vartheta'$ 
s redrawn so as to make the root the leftmost vertex and all the arrows pointing from right to left.

%%%%%%%%%%%%%%%%%%%%%%%%%%%%%%%%%%%%%%%%%%%%%%%%%%%%%%%%%%%%%%%%%%%%%%%% 
\begin{figure}[ht] 
\vskip.3truecm 
\centering 
\ins{027pt}{-028pt}{$\vartheta=$}
\ins{188pt}{-027.5pt}{$\vartheta'=$}
\ins{260pt}{-031.5pt}{$=$}
\ins{082pt}{-037pt}{$\nu_1$}
\ins{126pt}{-011pt}{$\nu_2$}
\ins{126pt}{-065pt}{$\nu_3$}
\ins{208pt}{-037pt}{$\nu_1$}
\ins{252pt}{-011pt}{$\nu_2$}
\ins{252pt}{-065pt}{$\nu_3$}
\ins{353pt}{-037pt}{$\nu_1$}
\ins{312pt}{-037pt}{$\nu_2$}
\ins{394pt}{-037pt}{$\nu_3$}
\includegraphics[width=1.8in]{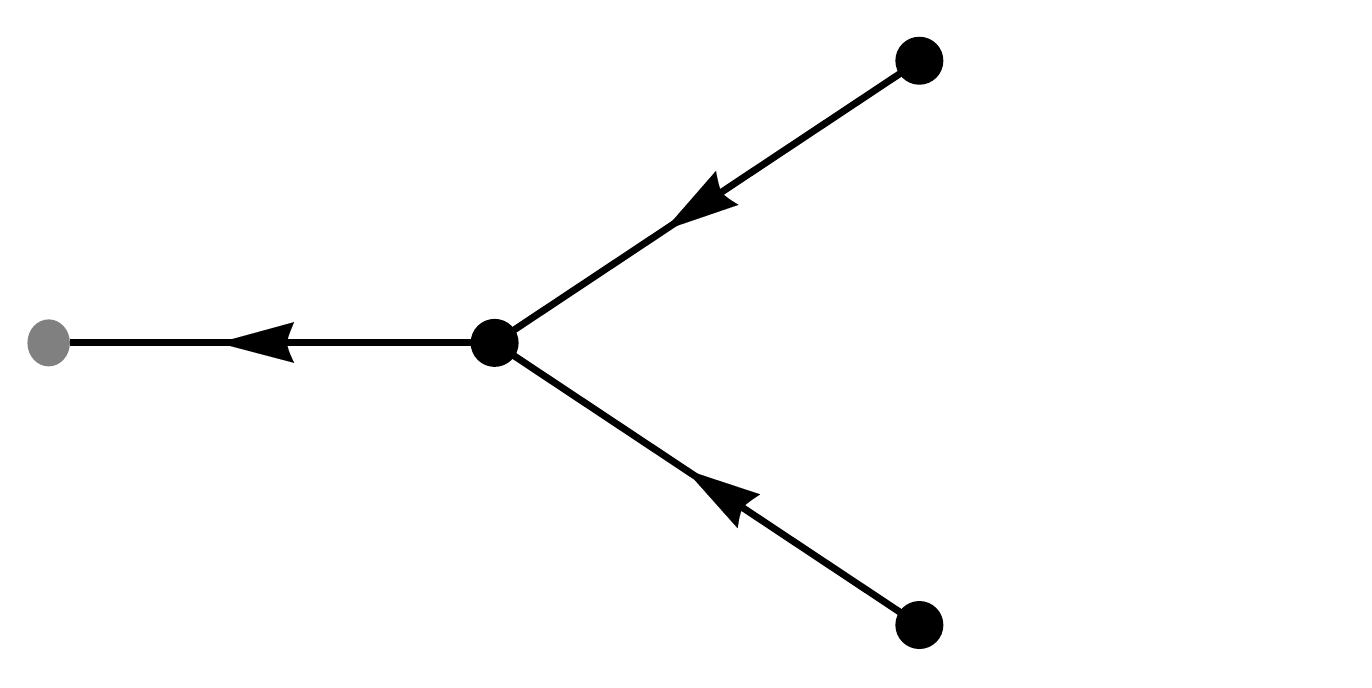}
\hspace{-.4cm}
\includegraphics[width=1.8in]{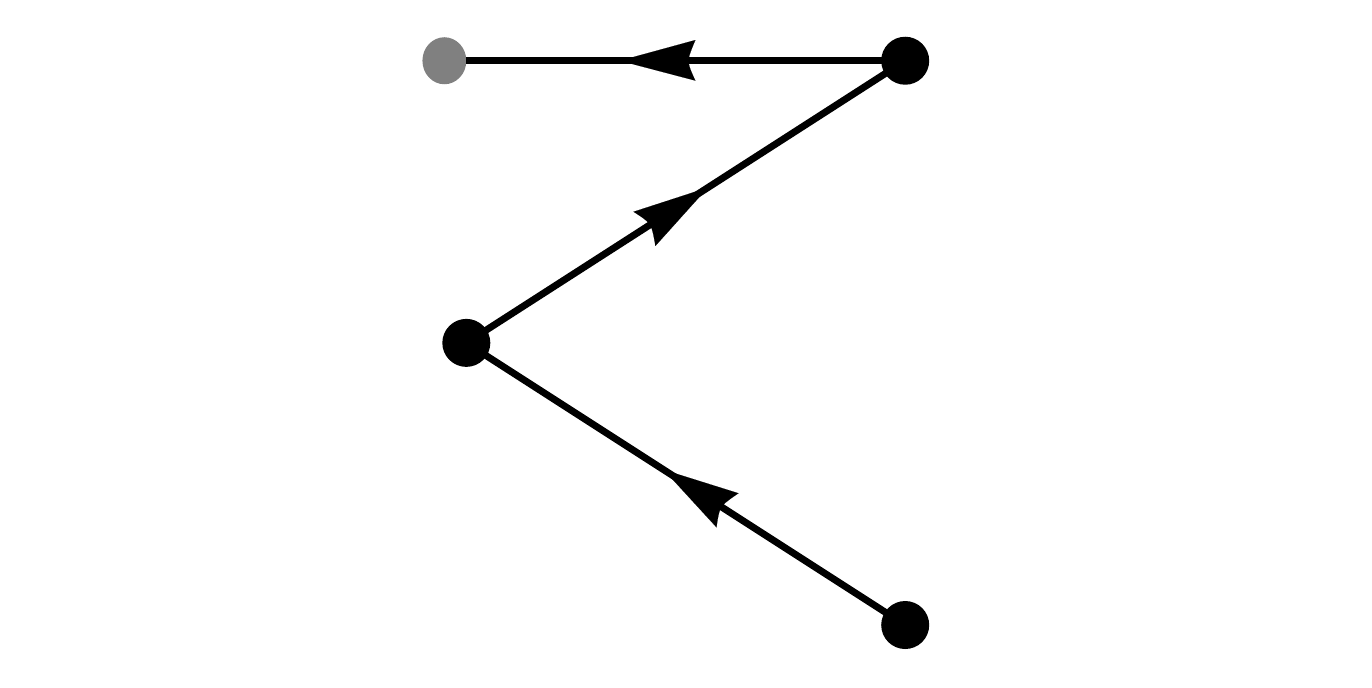}
\hspace{-1.2cm}
\includegraphics[width=1.8in]{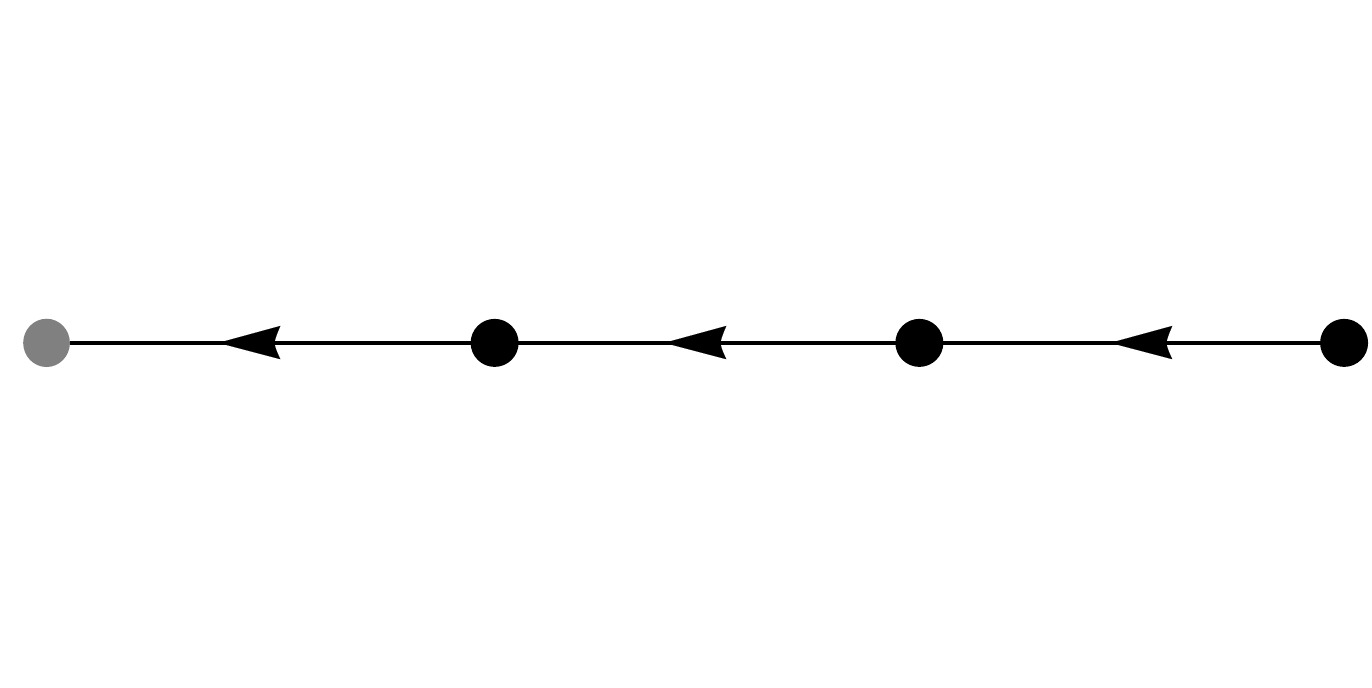}
\vspace{.3cm}
\caption{The tree $\vartheta'$ obtained from $\vartheta$ by shifting the root line (only the mode labels are shown).}
\label{attach} 
\end{figure} 
%%%%%%%%%%%%%%%%%%%%%%%%%%%%%%%%%%%%%%%%%%%%%%%%%%%%%%%%%%%%%%%%%%%%%%%% 

Let $\Theta^{(k)}_\nu$ be the set of trees $\vartheta$ with $k(\vartheta)=k$, i.e.~with $k$ nodes, such that the root line has momentum $\nu$.
Then it is not difficult to check by induction that
\begin{equation}\label{rappre}
u^{(k)}_\nu = \sum_{\vartheta\in\Theta^{(k)}_\nu}\Val(\vartheta),\qquad \nu\ne 0 ,
\end{equation}
with $\Val(\vartheta)$ defined in \eqref{val}.
Our aim is to prove a bound like
\begin{equation}\label{hovinto}
\sum_{\nu\in\ZZZ^\ZZZ_f\setminus\{0\}}    \Big\| \sum_{\vartheta\in\Theta^{(k)}_\nu}\Val(\vartheta)\Big\| e^{s|\nu|_\star} \leq C^k
\end{equation}
for some constant $C>0$. Indeed, if \eqref{hovinto} holds true, Proposition \ref{esiste1} follows with $\e_1= C^{-1}$.

We need to bound the sum over all possible tree values in
order to obtain the convergence of the Lindstedt series \eqref{uu}. The number of trees with $k$ nodes
is easily bounded by $4^k$ \cite{GM2}.  Moreover one can prove by induction that
\begin{equation}\label{nome}
\sum_{v\in N(\vartheta)}p_v =k(\vartheta) -1 .
\end{equation}

Regarding the product of the node factors, we have
\begin{equation}\label{sononumeri}
\biggl\| \prod_{v\in N(\vartheta)} \calF_v \biggr\| \le  \prod_{v\in N(\vartheta)}\frac{1}{p_v!} \|\nu_v^{p_v+1}\|_{\rm op} |f_{\nu_v}| .
\end{equation}
Note that the  product in the r.h.s.~of \eqref{sononumeri}
is a finite product of positive real numbers, so, while the product in the l.h.s.~is suitably ordered,
in the r.h.s.~the order is irrelevant. In the light of the bound \eqref{sononumeri},
the only possible obstruction to the convergence of the series \eqref{solou}
is due to the presence of the propagators, i.e.~to the small divisors.

%%%%%%%%%%%%%%%%%%%%%%%%%%%%%%%%%%%%%%%%%%%%%%%%%%%%%%%%%%%%%%%%%%%%%%%%%% 
%%%%%%%%%%%%%%%%%%%%%%%%%%%%%%%%%%%%%%%%%%%%%%%%%%%%%%%%%%%%%%%%%%%%%%%%%% 
\zerarcounters 
\section{Resonant clusters} 
\label{sec} 
%%%%%%%%%%%%%%%%%%%%%%%%%%%%%%%%%%%%%%%%%%%%%%%%%%%%%%%%%%%%%%%%%%%%%%%%%% 
%%%%%%%%%%%%%%%%%%%%%%%%%%%%%%%%%%%%%%%%%%%%%%%%%%%%%%%%%%%%%%%%%%%%%%%%%% 

As argued at the end of the previous section, the factors that give rise to difficulties in proving the convergence of the
Lindstedt series are the propagators. This happens even in the finite-dimensional case and with Diophantine frequencies.
Let us report two significant examples of trees with $k$ nodes for which the product of the propagators 
is of order $k^k$, in the simpler case of $\om\in\RRR^d$ satisfying \eqref{diop}
(note that, being a finite-dimensional case, in the two example discussed below the choice of the norm is irrelevant).

%%%%%%%%%%%%%%%%%%%%%%%%%%%%%%%%%%%%%%%%%%%%%%%%%%%%%%%%%%%%%%%%%%%%%%%%%%
\begin{ex}\label{nonprob} 
\emph{
Let $\vartheta$ be the linear tree with $k$ nodes $v_1,\ldots,v_k$,  depicted in Figure \ref{noprobbro},
such that $\nu_{v_1}=\ldots=\nu_{v_{k-1}}=0$ and $\nu_{v_k}=\nu$, with $\nu\in\ZZZ^d$ such that
\begin{equation}\label{sim}
|\om\cdot\nu| \sim \frac{\g}{\|\nu\|_1^\tau} . \end{equation}
The conservation law \eqref{conserva} implies that all the lines (including the root-line) have momentum equal to $\nu$.
If $n$ is such that $\Psi_{n}(\om\cdot\nu)=1$  and $n_{\ell}=n$ for all $\ell\in L(\vartheta)$,
the product of the propagators times $f_\nu$ is therefore equal to
\begin{equation} \nonumber f_\nu \prod_{\ell\in L(\vartheta)} \calG_\ell = \left(\frac{1}{(\om\cdot\nu)^2}\right)^k
\sim \h_\nu e^{-2s \|\nu\|_1} \left(\frac{\|\nu\|_1^{2\tau}}{\gamma^2}\right)^k \sim \h_\nu e^{-s \|\nu\|_1} \frac{\lceil 2\tau k \rceil !}{\gamma^{2k} s^{2\tau k}} ,
\end{equation}
where we have written $f_\nu = \h_\nu e^{-2s \|\nu\|_1}$, with $\h_\nu$ bounded in $\nu$ (in fact, summable over $\nu$) because of \eqref{fouriero}.
However such a tree is not really a source of problems because $\calF_v=0$ for all $v$ except for $v_k$, and thus $\Val(\vartheta)=0$.
}
\end{ex}
%%%%%%%%%%%%%%%%%%%%%%%%%%%%%%%%%%%%%%%%%%%%%%%%%%%%%%%%%%%%%%%%%%%%%%%%%%

%%%%%%%%%%%%%%%%%%%%%%%%%%%%%%%%%%%%%%%%%%%%%%%%%%%%%%%%%%%%%%%%%%%%%%%%%%
% Figure 5
%%%%%%%%%%%%%%%%%%%%%%%%%%%%%%%%%%%%%%%%%%%%%%%%%%%%%%%%%%%%%%%%%%%%%%%%%%
\vspace{-.4cm}
\begin{figure}[ht] 
\centering 
\ins{022pt}{-30.5pt}{$\vartheta =$}
\ins{100pt}{-40pt}{$0$}
\ins{158pt}{-40pt}{$0$}
\ins{216pt}{-40pt}{$0$}
\ins{273pt}{-40pt}{$0$}
\ins{331pt}{-40pt}{$0$}
\ins{392pt}{-41pt}{$\nu$}
\includegraphics[width=5in]{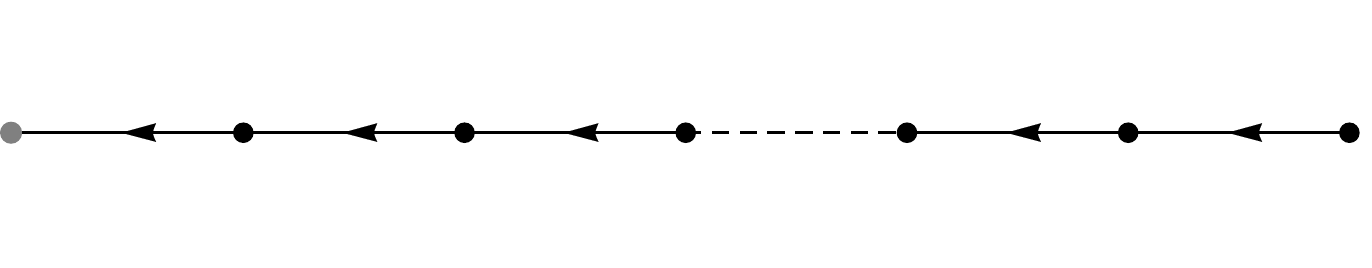} 
\caption{Graphical representation of the tree $\vartheta$ discussed in Example \ref{nonprob}.} 
\label{noprobbro} 
\end{figure} 
%%%%%%%%%%%%%%%%%%%%%%%%%%%%%%%%%%%%%%%%%%%%%%%%%%%%%%%%%%%%%%%%%%%%%%%%%%

The next example is somehow a generalisation of the previous one, and it contains all the crucial difficulties.
It shows that there are trees $\vartheta$ of any order $k$ such that $\Val(\vartheta)$
has a non-summable dependence on $k$.
 
%%%%%%%%%%%%%%%%%%%%%%%%%%%%%%%%%%%%%%%%%%%%%%%%%%%%%%%%%%%%%%%%%%%%%%%%
\begin{ex}\label{aiuto}
\emph{
Let $\vartheta$ be the tree with $k$ nodes in Figure \ref{probbro}, with $k$ odd,
and the mode labels arranged in such a way that $\vartheta$ has $(k+1)/{2}$ lines with momentum $\nu$,
with $\nu$ satisfying \eqref{sim}, and $(k-1)/{2}$ lines with momentum $\nu_0$, with $\|\nu_0\|_1 \sim 1$.
Assume also that the scale labels are such that $\Psi_{n_\ell}(\om\cdot\nu)=1$ for all $\ell\in L(\vartheta)$ -- otherwise $\Val(\vartheta)=0$.
Then, we have
\begin{equation} \label{above0}
\Val(\vartheta) = 
\ii^k \nu \left( - \frac{1}{2} |f_{\nu_0}|^2 \right)^{\frac{k-1}{2}} \!\!\!\! f_\nu
\left( -\nu_0\cdot\nu_0 \right)^{\frac{k-1}{2}} \left(\nu_0\cdot\nu_0\right)^{\frac{k-3}{2}} 
\left(\frac{1}{(\om\cdot\nu_0)^2}\right)^{\frac{k-1}{2}}
    	\left(\frac{1}{(\om\cdot\nu)^2}\right)^{\frac{k+1}{2}} 
\end{equation}
so that, if $f$ satisfies \eqref{fouriero} with $|\cdot|_\star = \|\cdot\|_1$, we may bound,
writing again $f_\nu = \h_\nu e^{-2s \|\nu\|_1}$ as in Example \ref{nonprob},
\begin{equation}\label{above}
\begin{aligned}
\|\Val(\vartheta)\|
& \sim 
\frac{\|\nu_0\|_1^{2(k-1)} |f_{\nu_0}|^{k-1}}{2^{\frac{k-1}{2}}|\om\cdot\nu_0|^{k-1}}
\frac{\|\nu\|_1|f_\nu|}{|\om\cdot\nu|^{k+1}} 
\sim C^k (\h_{\nu_0})^{k-1} \h_\nu
 \g^{-(k+1)} \|\nu\|_1^{\tau(k+1)+1} e^{-2s\|\nu\|_1} \\
& \sim (\h_{\nu_0})^{k-1} \h_\nu e^{-s\|\nu\|_1} {C}^k \g^{-(k+1)} s^{-(\tau(k+1)+1)} \lceil \tau(k+1)+1 \rceil ! ,
\phantom{\int}
\end{aligned}
\end{equation}
for a suitable constant ${C}$ independent of $\g$ and $\tau$.
Differently from Example \ref{nonprob}, there is no reason for this term to vanish, 
and hence the presence of such a contribution is something to worry about if one wants to prove the convergence
of the Lindstedt series.
}
\end{ex}
%%%%%%%%%%%%%%%%%%%%%%%%%%%%%%%%%%%%%%%%%%%%%%%%%%%%%%%%%%%%%%%%%%%%%%%%

%%%%%%%%%%%%%%%%%%%%%%%%%%%%%%%%%%%%%%%%%%%%%%%%%%%%%%%%%%%%%%%%%%%%%%%% 
% Figure 6
%%%%%%%%%%%%%%%%%%%%%%%%%%%%%%%%%%%%%%%%%%%%%%%%%%%%%%%%%%%%%%%%%%%%%%%% 
\begin{figure}[ht] 
\centering 
\ins{023pt}{-01.5pt}{$\vartheta =$}
\ins{098pt}{-10pt}{$\nu_0$}
\ins{122pt}{-68pt}{$-\nu_0$}
\ins{155pt}{-10pt}{$\nu_0$}
\ins{181pt}{-68pt}{$-\nu_0$}
\ins{213pt}{-10pt}{$\nu_0$}
\ins{270pt}{-10pt}{$\nu_0$}
\ins{300pt}{-68pt}{$-\nu_0$}
\ins{236pt}{-68pt}{$-\nu_0$}
\ins{328pt}{-10pt}{$\nu_0$}
\ins{353pt}{-68pt}{$-\nu_0$}
\ins{391pt}{-12pt}{$\nu$}
\includegraphics[width=5in]{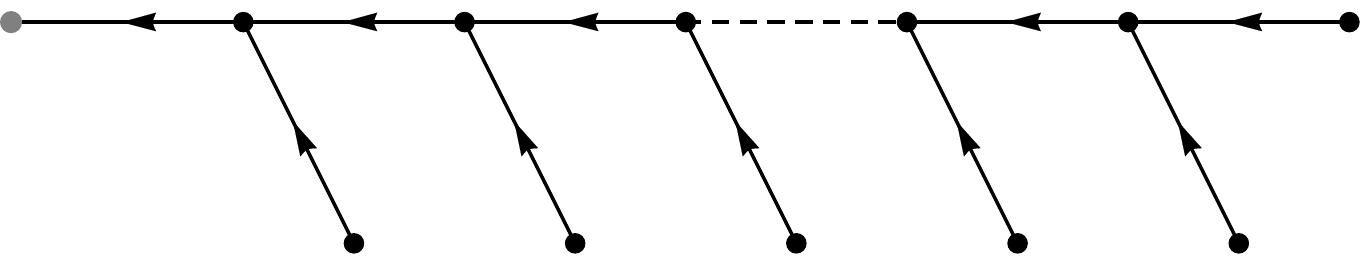} 
\vskip.2truecm 
\caption{Graphical representation of of the tree $\vartheta$ discussed in Example \ref{aiuto}.} 
\label{probbro} 
\end{figure} 
%%%%%%%%%%%%%%%%%%%%%%%%%%%%%%%%%%%%%%%%%%%%%%%%%%%%%%%%%%%%%%%%%%%%%%%%

%%%%%%%%%%%%%%%%%%%%%%%%%%%%%%%%%%%%%%%%%%%%%%%%%%%%%%%%%%%%%%%%%%%%%%%%
\begin{rmk}\label{modulofuori}
\emph{
Example \ref{aiuto} above suggests that a bound like 
\begin{equation}\label{magara}
\sum_{\nu\in\ZZZ^\ZZZ_f\setminus\{0\}} \sum_{\vartheta\in\Theta^{(k)}_\nu}\|\Val(\vartheta)\|e^{s|\nu|_\star} \leq C^k
\end{equation}
is impossible to achieve even in the simpler finite-dimensional case with Diophantine frequency vector.
However, we need the much weaker bound \eqref{hovinto} in order to prove the convergence of the Lindstedt series.
This means that, if by any chance there were to exist other trees $\vartheta_1,\vartheta_2,\ldots$ in $\Theta^{(k)}_\nu$ such that
\begin{equation}\label{speranza}
\Val(\vartheta)+\Val(\vartheta_1) + \Val(\vartheta_2) + \ldots = 0,
\end{equation} 
and the values of the remaining trees were summable over $k$, 
the problematic contribution would be deleted by other (equally problematic) contributions. 
In fact, the sum in \eqref{speranza} does not need to vanish: it suffices that it is small enough to be estimated by a summable quantity.
Of course, one first needs to understand what  a problematic contribution ultimately is.
}
\end{rmk}
%%%%%%%%%%%%%%%%%%%%%%%%%%%%%%%%%%%%%%%%%%%%%%%%%%%%%%%%%%%%%%%%%%%%%%%%

To better understand the problem related to trees like that in Example \ref{aiuto}, let us redraw Figure \ref{probbro} as in
Figure \ref{catenabis}: the momentum $\nu$ of each line entering an encircled subgraph is such that the small divisor $(\om\cdot\nu)^2$
can be very small, while both the small divisor $(\om\cdot\nu_0)^2$ and the factor $|f_{\nu_0}|^2$ associated to any subgraph are of order 1.

%%%%%%%%%%%%%%%%%%%%%%%%%%%%%%%%%%%%%%%%%%%%%%%%%%%%%%%%%%%%%%%%%%%%%%%% 
% Figure RC catena
%%%%%%%%%%%%%%%%%%%%%%%%%%%%%%%%%%%%%%%%%%%%%%%%%%%%%%%%%%%%%%%%%%%%%%%% 
\begin{figure}[H] 
\centering 
\ins{028pt}{-17pt}{$\vartheta =$}
\ins{108pt}{-11pt}{$\nu_0$}
\ins{119pt}{-73pt}{$-\nu_0$}
\ins{165pt}{-11pt}{$\nu_0$}
\ins{176pt}{-73pt}{$-\nu_0$}
\ins{221pt}{-11pt}{$\nu_0$}
\ins{280pt}{-11pt}{$\nu_0$}
\ins{290pt}{-73pt}{$-\nu_0$}
\ins{232pt}{-73pt}{$-\nu_0$}
\ins{335pt}{-11pt}{$\nu_0$}
\ins{347pt}{-73pt}{$-\nu_0$}
\ins{394pt}{-26pt}{$\nu$}
\includegraphics[width=5in]{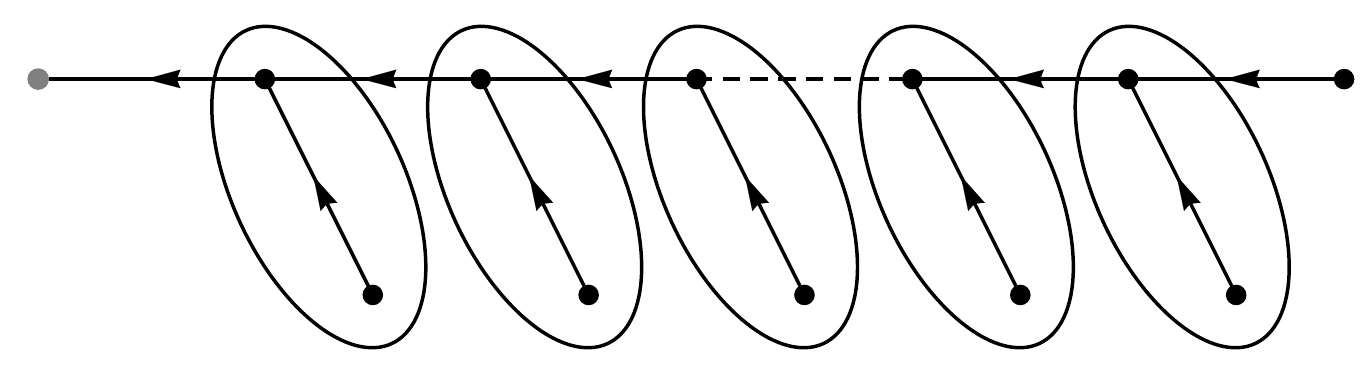} 
\vskip.2truecm 
\caption{Graphical representation of of the tree $\vartheta$ discussed in Example \ref{aiuto}.} 
\label{catenabis} 
\end{figure} 
%%%%%%%%%%%%%%%%%%%%%%%%%%%%%%%%%%%%%%%%%%%%%%%%%%%%%%%%%%%%%%%%%%%%%%%%

From the examples above, we can argue that problems arise (at least) whenever a tree $\vartheta$ contains 
subgraphs $T$, each of which has the following structure:
\begin{enumerate}[topsep=1ex]
\itemsep0em
\item $T$ has one exiting line $\ell_T$ and only one entering line $\ell'_T$;
\item the mode labels of the nodes $v\in N(T)$ sum up to zero, so that
$\nu_{\ell_T}=\nu_{\ell'_T}$;
\item  the norms $|\nu_v|_\star$ corresponding to the nodes $v\in N(T)$
are so small that the product of the Fourier coefficients
\[
\prod_{v \in N(T)} f_{\nu_v} \sim e^{-2s \sum_{v\in N(T)} |\nu_v|_\star} \prod_{v \in N(T)} \h_{\nu_v}
\]
of the node factors of the nodes in $T$ is not enough to balance the small divisor $|\om\cdot\nu_{\ell_T}|^2$ of the exiting line;
\item for all $\ell\in L(T)$ one has $|\om\cdot\nu_{\ell_T}| < |\om\cdot\nu_{\ell}|$.
\end{enumerate}

%%%%%%%%%%%%%%%%%%%%%%%%%%%%%%%%%%%%%%%%%%%%%%%%%%%%%%%%%%%%%%%%%%%%%%%%
\begin{rmk} \label{sottolineiamolo}
\emph{
Example \ref{aiuto} shows that it is possible to create a chain arbitrarily long in which such subgraphs are repeated,
as highlighted in Figure \ref{catenabis}, and this produces an accumulation of small divisors. In fact,
Lemma \ref{brjuno} below, together with the consequent estimates of Lemma \ref{convergerebbe},
yields that the subgraphs described above are the only source of problems to deal with.
}
\end{rmk}
%%%%%%%%%%%%%%%%%%%%%%%%%%%%%%%%%%%%%%%%%%%%%%%%%%%%%%%%%%%%%%%%%%%%%%%% 

Let us now make formal the argument above by introducing the notion of resonant cluster.\footnote{Also called
self-energy cluster or self-energy graph in the literature, again in analogy with Quantum Field Theory.}
First of all, for any tree $\vartheta$ and any subgraph $S$ of $\vartheta$, we define 
\[
K(S):= \sum_{v\in N(S)} |\nu_v|_\star .
\]

%%%%%%%%%%%%%%%%%%%%%%%%%%%%%%%%%%%%%%%%%%%%%%%%%%%%%%%%%%%%%%%%%%%%%%%% 
\begin{defi}\label{RC}
Given a tree $\vartheta$, a \emph{resonant cluster (RC)} is a connected subgraph $T$ of $\vartheta$ such that
\begin{enumerate}[topsep=1ex]
\itemsep0em
\item $T$ has only one entering line $\ell'_T$ and one exiting line $\ell_T$;
\item\label{item2} $\nu_{\ell_T}=\nu_{\ell'_T}$;
\item\label{item3} $K(T) < 2^{m_{\und{n}_T}-1}$;
\item\label{item4} $\ol{n}_T<\und{n}_T$, with $\ol{n}_T:=\max\{n_\ell :  \ell\in L(T)\}$ and $\und{n}_T :=\min\{n_{\ell_T},n_{\ell'_T}\}$.
\end{enumerate}
We say that $\ol{n}_T$ is the scale of $T$. Furthermore, we call \emph{resonant line} any line $\ell\in L(\vartheta)$ exiting a RC in $\vartheta$
and \emph{non-resonant line} any line which is not a resonant line. 
\end{defi}
%%%%%%%%%%%%%%%%%%%%%%%%%%%%%%%%%%%%%%%%%%%%%%%%%%%%%%%%%%%%%%%%%%%%%%%% 

%%%%%%%%%%%%%%%%%%%%%%%%%%%%%%%%%%%%%%%%%%%%%%%%%%%%%%%%%%%%%%%%%%%%%%%% 
\begin{rmk}\label{ordine2}
\emph{
With the scale functions \eqref{taglia}, item \ref{item2} in Definition \ref{RC} implies that,
if $\Val(\vartheta) \neq 0$ and $T$ is a RC in $\vartheta$, then $\und{n}_T =n_{\ell_T}=n_{\ell'_T}$.
This is not the case if one uses $C^\io$ functions with compact support (see Remark \ref{supporto}).
}
\end{rmk}
%%%%%%%%%%%%%%%%%%%%%%%%%%%%%%%%%%%%%%%%%%%%%%%%%%%%%%%%%%%%%%%%%%%%%%%% 

%%%%%%%%%%%%%%%%%%%%%%%%%%%%%%%%%%%%%%%%%%%%%%%%%%%%%%%%%%%%%%%%%%%%%%%% 
\begin{ex} \label{RC1}
\emph{
Examples of RCs $T$ with $K(T)=2$ are given in Figure \ref{sec2}, where only the subgraphs are depicted,
ignoring the rest of the tree; in all three cases one has the constraints that
$\nu_1+\nu_2=0$ and that the external lines have scale larger than the internal line.
Then, if we add to the tree $\vartheta$ of Example \ref{aiuto} the scale labels ands assume that $n_0$ and $n$ for which
$\Psi_{n_0}(\om\cdot\nu_0)\neq 0$ and $\Psi_{n}(\om\cdot\nu)\neq 0$ are such that $n_0<n$,
we obtain a chain of $(k-1)/2$ RCs, as illustrated in Figure \ref{catenabis}.
Note that each RC in Figura \ref{catenabis} is as the RC in the middle of Figure \ref{sec2}, with $\nu_1=-\nu_2=\nu_0$.
}
\end{ex}
%%%%%%%%%%%%%%%%%%%%%%%%%%%%%%%%%%%%%%%%%%%%%%%%%%%%%%%%%%%%%%%%%%%%%%%% 

%%%%%%%%%%%%%%%%%%%%%%%%%%%%%%%%%%%%%%%%%%%%%%%%%%%%%%%%%%%%%%%%%%%%%%%% 
% Figure RC
%%%%%%%%%%%%%%%%%%%%%%%%%%%%%%%%%%%%%%%%%%%%%%%%%%%%%%%%%%%%%%%%%%%%%%%% 
\begin{figure}[ht] 
\centering 
\ins{062pt}{-45pt}{$\nu_1$}
\ins{107pt}{-45pt}{$\nu_2$}
\ins{229pt}{-43pt}{$\nu_1$}
\ins{275pt}{-68pt}{$\nu_2$}
\ins{362pt}{-42pt}{$\nu_1$}
\ins{411pt}{-17pt}{$\nu_2$}
\includegraphics[width=2in]{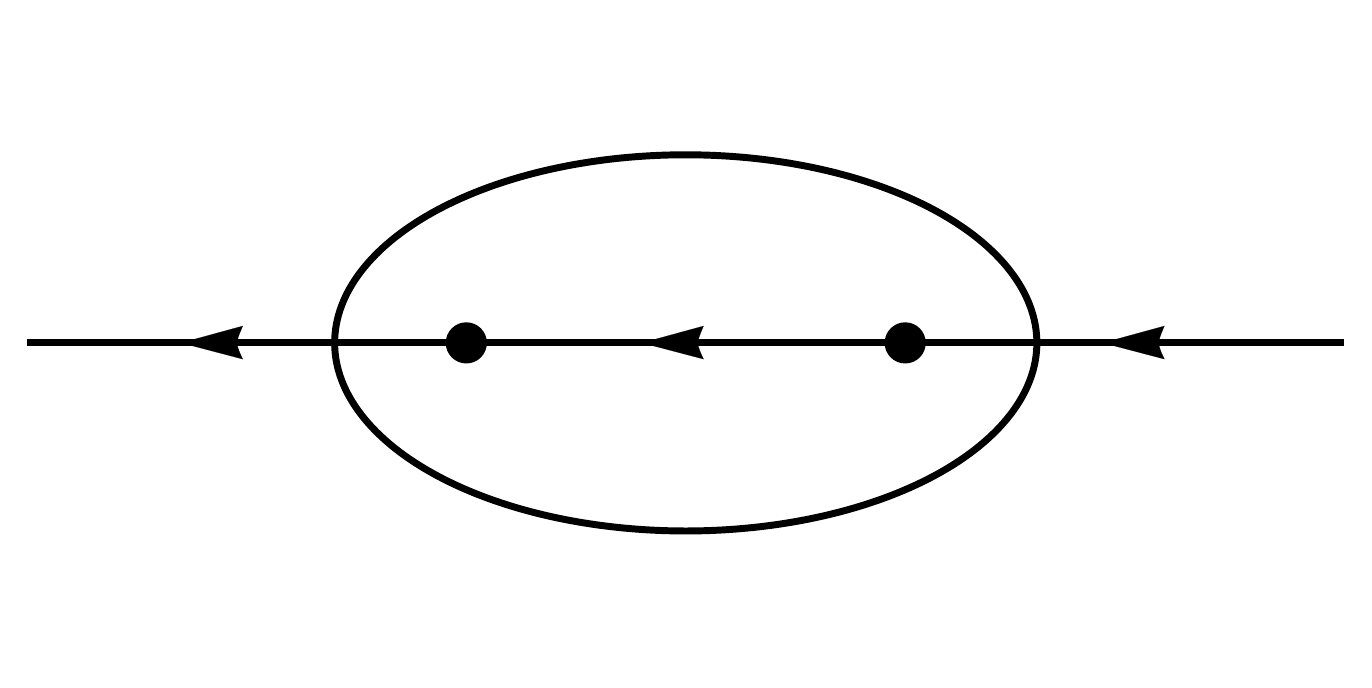}\vspace{-.2cm}
\hspace{.4cm}
\includegraphics[width=1.6in]{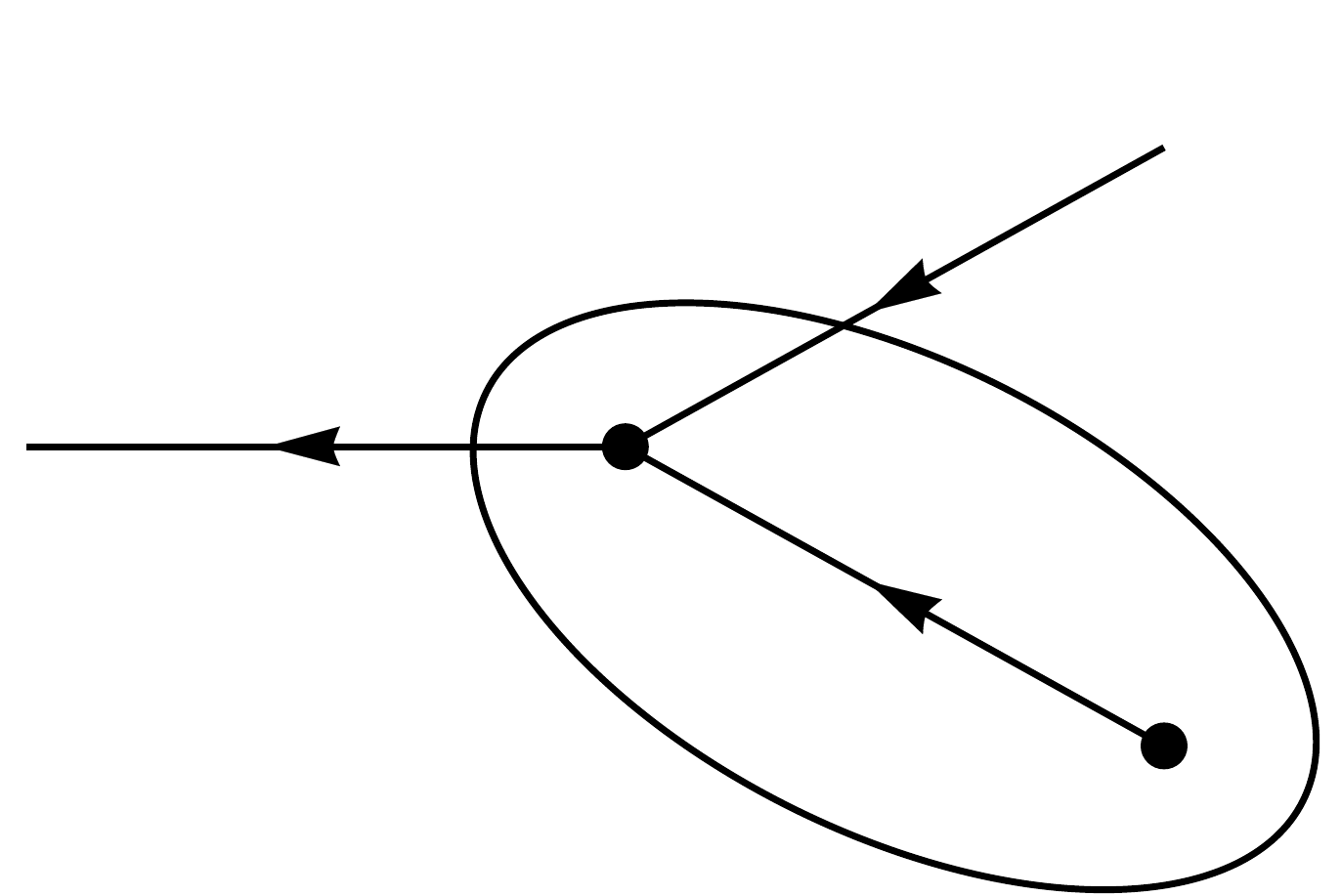} 
\hspace{.4cm}
\includegraphics[width=1.6in]{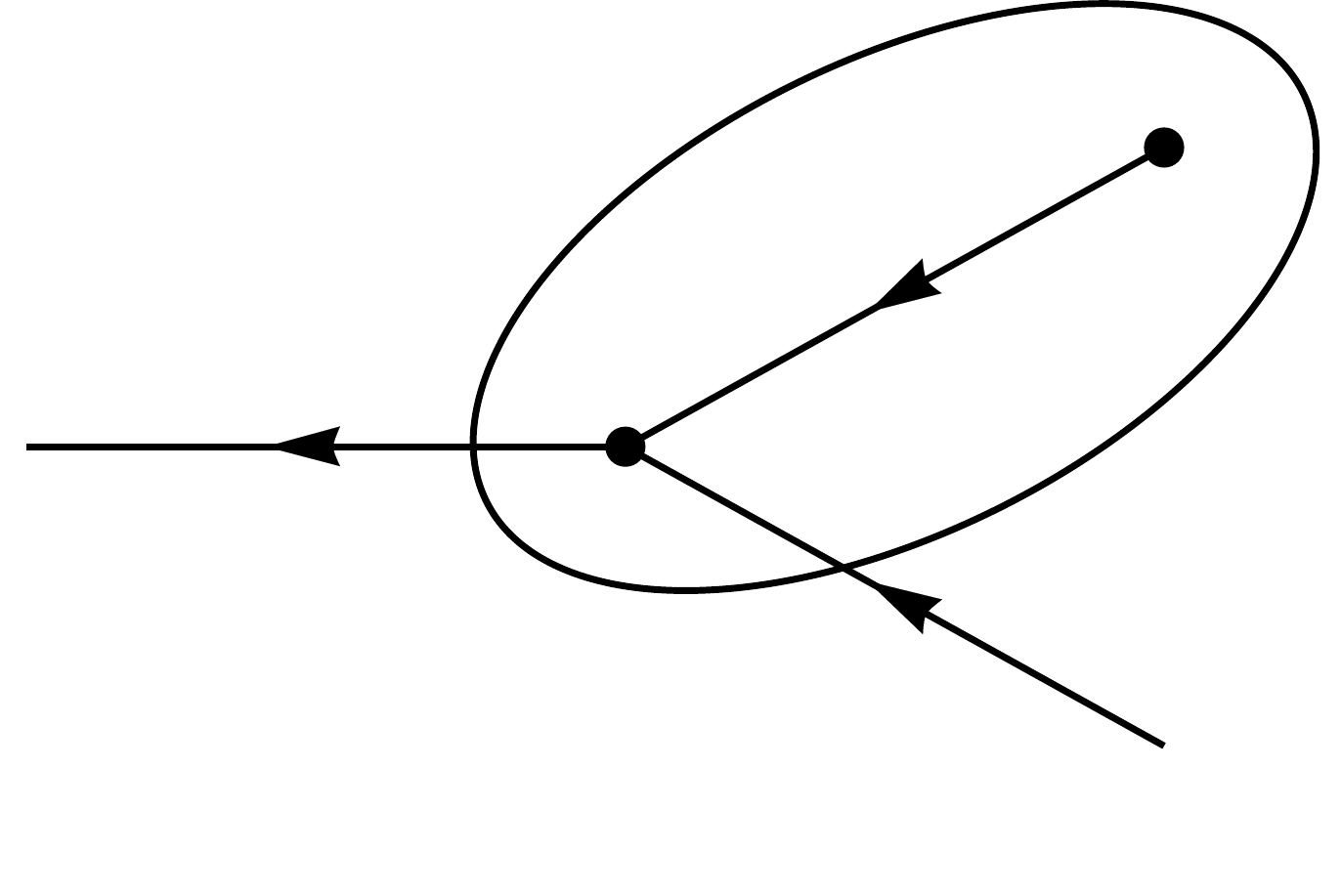} 
\vskip.2truecm 
\caption{Self-energy clusters of order 2 (see the text for details).}
\label{sec2} 
\end{figure} 
%%%%%%%%%%%%%%%%%%%%%%%%%%%%%%%%%%%%%%%%%%%%%%%%%%%%%%%%%%%%%%%%%%%%%%%%

Call $L_{\!N\!R}(\vartheta)$ the set of non-resonant lines in $\vartheta$ and, for $n\ge 0$, define
$\gotN_n(\vartheta)$ as the number of lines $\ell\in L_{\!N\!R}(\vartheta)$ such that $n_\ell \ge n$.
Now, we are in the position to prove the following version of the so-called Bryuno lemma, which gives
a bound on the non-resonant lines.

%%%%%%%%%%%%%%%%%%%%%%%%%%%%%%%%%%%%%%%%%%%%%%%%%%%%%%%%%%%%%%%%%%%%%%%%%%
\begin{lemma}\label{brjuno}
Let $\vartheta$ be such that $\Val(\vartheta)\ne0$.
For all $n\ge1$ one has
 \begin{equation}\label{bru}
 \gotN_n(\vartheta)\le \max\Bigl\{ 4 K(\vartheta) 2^{-m_n} -1,0\Bigr\} .
 \end{equation}
\end{lemma}
%%%%%%%%%%%%%%%%%%%%%%%%%%%%%%%%%%%%%%%%%%%%%%%%%%%%%%%%%%%%%%%%%%%%%%%%%%

%%%%%%%%%%%%%%%%%%%%%%%%%%%%%%%%%%%%%%%%%%%%%%%%%%%%%%%%%%%%%%%%%%%%%%%%%%
\prova
The proof is by induction on the number of nodes $k$ of the tree. If $k=1$ then there is only one line $\ell\in L(\vartheta)$.
If $n_\ell<n$ the bound trivially holds, otherwise, since $\Val(\vartheta)\ne0$ one has 
\[
|\om\cdot\nu_\ell| < \frac{1}{4} \be^\star_\om(m_{n_\ell-1})\le \frac{1}{2}\be^\star_\om(m_n-1)
\]
and hence $K(\vartheta)=|\nu_\ell|_\star>2^{m_n-1}$, so that the bound \eqref{bru} holds.

Assume now \eqref{bru} to hold for all $k'<k$ and consider a tree $\vartheta\in\Theta^{(k)}_\nu$.
If $n_{\ell_\vartheta}<n$, then the assertion follows
by the inductive hypothesis, otherwise let $\ell_1,\ldots,\ell_p$ be the lines on scale $\ge n$
closest to $\ell_\vartheta$, and let $\vartheta_1,\ldots,\vartheta_p$
be the subtrees of $\vartheta$ having $\ell_1,\ldots,\ell_p$ as root lines, respectively. If there is no such line,
the bound follows once more. If $p\ge2$ one has
\[
\gotN_n(\vartheta)\le\sum_{j=1}^p \gotN_n(\vartheta_j) \le 4 \sum_{j=1}^p K(\vartheta_j) 2^{-m_n} -p  \le 4 K(\vartheta) 2^{-m_n} -1.
\]
If $p=1$ let $T$ be the subgraph having $\ell_\vartheta$ as exiting line and $\ell_1$ as entering line. If $T$ is not a RC then,
since by construction
$\ol{n}_T <n\le n_{\ell_1}$, either item \ref{item2} or item \ref{item3} or both of Definition \ref{RC} fail.
In any case, one has $K(T)\ge 2^{m_n-1}$. Indeed, if $\nu_{\ell_\vartheta}\ne \nu_{\ell_1}$ one has
\begin{equation}\label{questo}
|\om\cdot(\nu_{\ell_\vartheta}-\nu_{\ell_1})| \le |\om\cdot\nu_{\ell_\vartheta}| + |\om\cdot\nu_{\ell_1}|
\le  \frac{1}{4} \be^\star_\om(m_{n_{\ell_\vartheta}-1}) +
\frac{1}{4} \be^\star_\om(m_{n_{\ell_1}-1}) \le \be^\star_\om(m_n-1) ,
\end{equation}
so that
\[
K(T) \ge |\nu_{\ell_\vartheta} - \nu_{\ell_1}|_\star > 2^{m_n -1}
\]
and hence
\[
\gotN_n(\vartheta)=1+\gotN_n(\vartheta_1) \le 4 K(\vartheta_1) 2^{-m_n} = 4 (K(\vartheta) - K(T)) 2^{-m_n} \le 4 K(\vartheta) 2^{-m_n} -1 .
\]
Otherwise, if $T$ is a RC then $\ell_\vartheta$ is resonant, so that
\[
\gotN_n(\vartheta)=\gotN_n(\vartheta_1) \le 4 K(\vartheta_1) 2^{-m_n} -1 \le 4 K(\vartheta) 2^{-m_n} -1.
\]
In all cases the assertion follows.
\EP

Thanks to Lemma \ref{brjuno} above we deduce, as aforementioned, that
RCs are in fact the only source of problems.

%%%%%%%%%%%%%%%%%%%%%%%%%%%%%%%%%%%%%%%%%%%%%%%%%%%%%%%%%%%%%%%%%%%%%%%%%%
\begin{lemma}\label{convergerebbe}
For any tree $\vartheta$, define 
\begin{equation} \nonumber
\Val_{\!N\!R}(\vartheta):= \Biggl( \prod_{v\in N(\vartheta)}\calF_v \Biggr) \Biggl( \prod_{\ell\in L_{\!N\!R}(\vartheta)}\calG_\ell \Biggr)
\end{equation} 
Then there is a positive constant $C_0$ such that
\[
\sum_{\nu\in\ZZZ^\ZZZ_f\setminus\{0\}} e^{s|\nu|_\star} \sum_{\vartheta\in\Theta^{(k)}_{\nu}}\|\Val_{\!N\!R}(\vartheta)\|\leq C_0^k.
\]
\end{lemma}
%%%%%%%%%%%%%%%%%%%%%%%%%%%%%%%%%%%%%%%%%%%%%%%%%%%%%%%%%%%%%%%%%%%%%%%%%%

%%%%%%%%%%%%%%%%%%%%%%%%%%%%%%%%%%%%%%%%%%%%%%%%%%%%%%%%%%%%%%%%%%%%%%%%%%
\prova
For any $n_0\ge0$ one has 
\begin{equation} \label{ultimaformulanumerata}
\begin{aligned}
\prod_{\ell\in L_{\!N\!R}(\vartheta)}|\calG_\ell|&\leq\prod_{n\ge0}\left(\frac{4}{(\be^\star_\om(m_n))^2}\right)^{\gotN_n(\vartheta)} 
\le \left(\frac{2}{\be^\star_\om(m_{n_0})}\right)^{2k}\prod_{n>n_0}\left(\frac{1}{\be^\star_\om(m_n)}\right)^{2\gotN_n(\vartheta)} \\
& \le (A(n_0))^{k} \prod_{n>n_0}\left(\frac{1}{\be^\star_\om(m_n)}\right)^{2^{-(m_n-3)}K(\vartheta) } \\
&= (A(n_0))^{k} \exp\Biggl( 8 K(\vartheta) \sum_{n>n_0 } \frac{1}{2^{m_n}}\log \left(\frac{1}{\be^\star_\om(m_n)}\right)
\Biggr)
\end{aligned}
\end{equation}
with $A(n_0) := (2/\be_\om^\star(m_{n_0}))^2$.
Then, since $\om$ is a $\star$-Bryuno vector, we can choose $n_0$ such that
\begin{equation} \label{bryunoo}
8 \sum_{n>n_0 } \frac{1}{2^{m_n}}\log \left(\frac{1}{\be^\star_\om(m_n)} \right)< \frac{s}{2} .
\end{equation}
Therefore, using that the number of unlabelled trees of order $k$ is bounded by $4^k$, we obtain 
\[
\begin{aligned}
& \sum_{\nu\in\ZZZ^\ZZZ_f\setminus\{0\}} e^{s|\nu|_\star}   \sum_{\vartheta\in\Theta^{(k)}_{\nu}} \|\Val_{\!N\!R}(\vartheta) \| \\     
& \qquad \qquad 
\stackrel{\eqref{sononumeri}}{\le} 
(A(n_0))^k \sum_{\nu\in\ZZZ^\ZZZ_f\setminus\{0\}} e^{s|\nu|_\star} 
\sum_{\vartheta\in\Theta^{(k)}_{\nu}}e^{s K(\vartheta)/2} \!\! \prod_{v\in N(\vartheta)}\frac{1}{p_v!}\|\nu_v^{p_v+1}\|_{\rm op}|f_{\nu_v}| \\
& \qquad \qquad 
\stackrel{\phantom{\eqref{sononumeri}}}{\le} 
(A(n_0))^k    \sum_{\nu\in\ZZZ^\ZZZ_f\setminus\{0\}}  
\sum_{\vartheta\in\Theta^{(k)}_{\nu}}e^{2s K(\vartheta)} \!\! 
\prod_{v\in N(\vartheta)}\frac{1}{p_v!}\|\nu_v\|^{p_v+1}_{1}|f_{\nu_v}|e^{-s|\nu_v|_\star/2}  \\
& \qquad \qquad 
\stackrel{\eqref{jedoernome}}{\le} 
(A(n_0))^k   \sum_{\nu\in\ZZZ^\ZZZ_f\setminus\{0\}}  
\sum_{\vartheta\in\Theta^{(k)}_{\nu}}e^{2s K(\vartheta)} \!\! 
\prod_{v\in N(\vartheta)}|f_{\nu_v}|
(p_v+1) \Bigg(\frac{2c_1}{s}\Bigg)^{(p_v+1)} \\
& \qquad \qquad 
\stackrel{\eqref{nome}}{\le}  \Bigg(\frac{16c_1^2}{s^2}A(n_0)\Bigg)^k   \sum_{\nu\in\ZZZ^\ZZZ_f\setminus\{0\}}  
\sum_{\vartheta\in\Theta^{(k)}_{\nu}}
\prod_{v\in N(\vartheta)}|f_{\nu_v}|e^{2s |\nu_v|_\star} \\
& \qquad \qquad 
\stackrel{\phantom{\eqref{sononumeri}}}{\le} 
\Bigg(\frac{64c_1^2}{s^2}A(n_0)\Bigg)^k \sum_{\nu_1,\ldots,\nu_k\in\ZZZ^\ZZZ_f} \prod_{i=1}^k |f_{\nu_i}|e^{2s |\nu_i|_\star} 
\le \Bigg(\frac{64c_1^2}{s^2}A(n_0) \|f\|_{2s,\star,\RRR}\Bigg)^k,
\end{aligned}
\]
with $c_1$ as in \eqref{jedoernome}, so that, setting
\begin{equation} \label{C0}
C_0 = \left( \frac{16 c_1}{s \be_\om^\star(m_{n_0})} \right)^2 \|f\|_{2s,\star,\RRR} .
\end{equation}
the assertion follows.
\EP
%%%%%%%%%%%%%%%%%%%%%%%%%%%%%%%%%%%%%%%%%%%%%%%%%%%%%%%%%%%%%%%%%%%%%%%%%%

%%%%%%%%%%%%%%%%%%%%%%%%%%%%%%%%%%%%%%%%%%%%%%%%%%%%%%%%%%%%%%%%%%%%%%%%%%
\begin{rmk} \label{badC}
\emph{
By reasoning as in the proof of Lemma \ref{convergerebbe} one can prove that
\[
\sum_{\nu\in\ZZZ^\ZZZ_f\setminus\{0\}} e^{s|\nu|_\star} \sum_{\vartheta\in\Theta^{(k)}_{\nu}} \|\Val(\vartheta) \|  \leq (C(k))^k ,
\]
with $C(k)\to+\io$ as $k\to+\io$. Such a bound implies that, for all $k$ and all $\nu\ne0$, the coefficient
$u^{(k)}_\nu$ is well defined in $\ell^\io(\CCC)$.
However, this is not enough to obtain the convergence of the Lindstedt series.
}
\end{rmk}
%%%%%%%%%%%%%%%%%%%%%%%%%%%%%%%%%%%%%%%%%%%%%%%%%%%%%%%%%%%%%%%%%%%%%%%%%%

%%%%%%%%%%%%%%%%%%%%%%%%%%%%%%%%%%%%%%%%%%%%%%%%%%%%%%%%%%%%%%%%%%%%%%%%%%
%%%%%%%%%%%%%%%%%%%%%%%%%%%%%%%%%%%%%%%%%%%%%%%%%%%%%%%%%%%%%%%%%%%%%%%%%%
\zerarcounters 
\section{Symmetries and cancellations} 
\label{canc} 
%%%%%%%%%%%%%%%%%%%%%%%%%%%%%%%%%%%%%%%%%%%%%%%%%%%%%%%%%%%%%%%%%%%%%%%%%%
%%%%%%%%%%%%%%%%%%%%%%%%%%%%%%%%%%%%%%%%%%%%%%%%%%%%%%%%%%%%%%%%%%%%%%%%%%

In order to deal with the RCs, let us go back to Examples \ref{nonprob} and \ref{aiuto}.
We can think of the tree in Example \ref{nonprob} as a tree with a chain of \emph{trivial RCs},
whose values ``incidentally'' vanish. In fact, we can interpret the vanishing of the value of the trivial RC
as due to the fact that, according to \eqref{fase}, it can be regarded as the derivative w.r.t.~$\f_0$ of 
a function, $f_0(\f_0)$, which, being an average, is $\f_0$-independent.
In general, if, instead of one RC, we consider the sum of the  values of all RCs with the same order
and the same scale, the leading order of such quantity also turns out to be the derivative of an average 
and hence it vanishes as well. Let us clarify the last statement.

First of all, we introduce some notation. Given a tree $\vartheta$ and a subgraph $S$ in $\vartheta$
with only one exiting line $\ell_S$ and only one exiting line $\ell_{S}'$ (see Section \ref{alberi}), let
$\vartheta_{\ell_S}$ and $\vartheta_{\ell'_S}$ denote the subtrees of $\vartheta$ having $\ell_S$ and $\ell'_S$ as root line, respectively.
Then we can write
\begin{equation}\label{comelovedo}
\Val(\vartheta_{\ell_S}) =\calG_{\ell_S}\matW(S) [\Val(\vartheta_{\ell'_S})],
\end{equation}
which implicitly defines the linear operator
\[
\matW(S):\ell^\io(\CCC) \longrightarrow \ell^\io(\CCC) .
\]
By comparing \eqref{val} with \eqref{comelovedo}, if we write
$v_S$ for the node which $\ell_S$ exits, and $v'_S$ for the node which $\ell_S'$ enters, we see that
\begin{equation}\label{valS}
\matW(S) [z]
=  \frac{1}{p_{v_S}!}\ii \nu_{v_S} f_{\nu_{v_S}} \Biggl( \prod_{v\in N(S)\setminus \{v_S\}} 
\frac{1}{p_{v}!} (\ii \nu_{\pi(v)}\cdot \ii \nu_v ) f_{\nu_v} \Biggl)
\Bigg(\prod_{\ell\in L(S)} \calG_\ell\Bigg) \ii \nu_{v'_S} \cdot z ,
\end{equation}
We call $\matW(S)$ the \emph{value} of the subgraph $S$.

%%%%%%%%%%%%%%%%%%%%%%%%%%%%%%%%%%%%%%%%%%%%%%%%%%%%%%%%%%%%%%%%%%%%%%%%%%
\begin{rmk} \label{piS}
\emph{
Since for all $v\in N(\vartheta)$ there exists a path $\calP(r,v)$, with $r$ being the root of $\vartheta$,
the subgraph $S$ uniquely determines a subgraph $\pi(\vartheta,S)$,
with only one entering line and one exiting line, given respectively by $\ell_S$ and $\ell_\vartheta$, such that
\begin{equation}\label{comelovedoB}
\Val(\vartheta) = \calG_{\ell_\vartheta}\matW(\pi(\theta,S)) [\Val(\vartheta_{\ell_S})] =
\calG_{\ell_\vartheta}\matW(\pi(\vartheta,S)) [\calG_{\ell_S}\matW(S) [\Val(\vartheta_{\ell'_S})]].
\end{equation}
}
\end{rmk}
%%%%%%%%%%%%%%%%%%%%%%%%%%%%%%%%%%%%%%%%%%%%%%%%%%%%%%%%%%%%%%%%%%%%%%%%%%

%%%%%%%%%%%%%%%%%%%%%%%%%%%%%%%%%%%%%%%%%%%%%%%%%%%%%%%%%%%%%%%%%%%%%%%%%%
\begin{rmk} \label{piSbis}
\emph{
We can write \eqref{valS} as
\begin{equation}\label{valSbis}
\matW(S) [z] = \Bigg( \prod_{v\in N(S)} \calF_v\Bigg) [z]
\Bigg(\prod_{\ell\in L(S)} \calG_\ell\Bigg) ,
\end{equation}
where the product of the node factors is defined recursively on the order $k(S)$ as follows:
\begin{itemize}[topsep=0.5ex]
\itemsep0em
\item for any $v\in N(S)$ such that $p_v \ge 1$ and $\ell_v \notin \calP_S$,
if $\vartheta_{v,1},\ldots,\vartheta_{v,p_v}$ are the trees entering $v$, the product
\[
\Bigg( \prod_{\substack{ w \in N(S) \\ w \preceq v }} \calF_w \Bigg)
\Bigg( \prod_{\substack{ \ell\in L(S) \\ \ell \preceq v}} \calG_\ell\Bigg) =
\frac{1}{p_v !}  f_{\nu_v} (\ii\nu_v)^{p_v+1} [\Val(\vartheta_{v,1}),\ldots,\Val(\vartheta_{v,p_v})] ,
\]
is meant as in Remark \ref{lospostoqui}, since $v_S' \notin N(\vartheta_{\ell_v})$ in such a case;
\item or any node $v\neq v_S'$ such that $\ell_v \in \calP_S$ one has
\[
\begin{aligned}
& \Bigg( \prod_{\substack{ w \in N(S) \\ w \preceq v }} \calF_w \Bigg) [z]
\Bigg( \prod_{\substack{ \ell\in L(S) \\ \ell \preceq v}} \calG_\ell\Bigg) 
\\ & \qquad \qquad 
=
\frac{1}{p_v !}  f_{\nu_v} (\ii\nu_v)^{p_v+1}
[\Val(\vartheta_{v,1}),\ldots,\Val(\vartheta_{v,r-1}), \calG_{\ell_{v,r}} \matW(S_v)[z] ,\Val(\vartheta_{v,r+1}),\ldots,\Val(\vartheta_{v,p_v})] ,
\end{aligned}
\]
if $\vartheta_{v,1},\ldots,\vartheta_{v,p_v}$ are the trees entering $v$,
with $\vartheta_{v,r}$ being the subtree with root line $\ell_{\vartheta_{v,r}} \in \calP_S$,
and $S_v$ is the subgraph with one exiting line $\ell_{S_v}=\ell_{\vartheta_{v,r}}$ and only one entering line $\ell_{S_v}'=\ell_{S}'$;
\item for $v=v_S'$ one has
\[
\begin{aligned}
& \Bigg( \prod_{\substack{ w \in N(S) \\ w \preceq v_S'}} \calF_w \Bigg) [z]
\Bigg( \prod_{\substack{ \ell\in L(S) \\ \ell \preceq v_S'}} \calG_\ell\Bigg)
\\& \qquad \qquad  
=
\frac{1}{p_{v_S'}!} f_{\nu_{v_S'}} ( \ii \nu_{v_S'} )^{p_v+1}
[\Val(\vartheta_{v_{S}',1}),\ldots,\Val(\vartheta_{v_{S}',r-1}),z,\Val(\vartheta_{v_{S}',r+1}),\Val(\vartheta_{v_{S}',p_v-1})],
\end{aligned}
\]
if $\vartheta_{v_S',1},\ldots,\vartheta_{v_S',p_v}$ are the trees entering $v_{S}'$ and $\vartheta_{v_S',r}=\vartheta_{\ell_S'}$.
\end{itemize}
In particular, while the product of the node factors in \eqref{val} is an element of $\ell^\io(\CCC)$, the product in \eqref{valSbis}
is an element of $\LL(\ell^\io(\CCC),\ell^\io(\CCC))$.
}
\end{rmk}
%%%%%%%%%%%%%%%%%%%%%%%%%%%%%%%%%%%%%%%%%%%%%%%%%%%%%%%%%%%%%%%%%%%%%%%%%%

Next, given a RC $T$, for all $\ell\in L(T)$ we set
\begin{equation}\label{nu0}
\nu_\ell^0 :=\left\{
\begin{aligned}
&\nu_\ell,\qquad & \ell\notin\calP_T , \\
&\nu_\ell-\nu_{\ell'_T} ,\qquad & \ell \in \calP_T ,
\end{aligned}
\right.
\end{equation}
and introduce the \emph{symbol function} $\x_\ell:\RRR\to\RRR$ as
\begin{equation}\label{xi}
\x_\ell (x) := \left\{
\begin{aligned}
&\om\cdot\nu_\ell^0,\qquad & \ell\notin\calP_T,  \\
&\om\cdot\nu_\ell^0 + x ,\qquad & \ell\in\calP_T.
\end{aligned}
\right.
\end{equation}
Then, we define\footnote{Given two Banach spaces $X$ and $Y$, $\LL(X,Y)$
denotes the space of continuous linear operators $L:X\to Y$.} $\matW_T:\RRR\to \LL(\ell^\io(\CCC),\ell^{\io}(\CCC))$ as
\begin{equation}\label{valfun}
\matW_T(x)[z]:=
\frac{1}{p_{v_T}!}\ii \nu_{v_T} f_{\nu_{v_T}} \Biggl( \prod_{v\in N(\vartheta)\setminus \{v_T\}} 
\frac{1}{p_{v}!} (\ii \nu_{\pi(v)}\cdot \ii \nu_v ) f_{\nu_v} \Biggl)
\Bigg(\prod_{\ell\in L(T)} \calG_{n_\ell}(\x_\ell(x))\Bigg)
\ii \nu_{v'_T} \cdot z ,
\end{equation}
so that $\matW(T) = \matW_T(\om\cdot\nu_{\ell'_T})$. We write
\begin{equation}\label{quello}
\matW(T)= \matW_T(\om\cdot\nu_{\ell'_T}) =  \matW_T(0) + \del_x  \matW_T(0)\om\cdot\nu_{\ell'_T} + (\om\cdot\nu_{\ell'_T})^2
R_T (\om\cdot\nu_{\ell'_T}) ,
\end{equation}
with $R_T(\cdot) \in\LL(\ell^\io(\CCC),\ell^\io(\CCC))$.

%%%%%%%%%%%%%%%%%%%%%%%%%%%%%%%%%%%%%%%%%%%%%%%%%%%%%%%%%%%%%%%%%%%%%%%%%%
\begin{rmk} \label{vanishing}
\emph{
In order to bound $R_T (\om\cdot\nu_{\ell'_T})$ it comes in handy to write it in terms of the second derivative of $\matW_T(x)$,
by using the integral form for the remainder of the Taylor formula (see also Remark \ref{whysmooth} below). Therefore,
when computing both $\del_x  \matW_T(0)$ and $R_T (\om\cdot\nu_{\ell'_T})$,
it may happen that for some line $\ell\in L(T)$ the quantity $\xi_\ell(x)$ crosses the boundary of the support of the scale function $\Psi_{n_\ell}(\xi_\ell(x))$,
and, as consequence, differentiating the scale functions of the form \eqref{taglia} could be a problem.
In principle one could envisage to construct a suitable partition of unity, adapted to the frequency vector $\om$,
such that the small divisor $\xi_\ell(x)$ never go out the support of $\Psi_{n_\ell}$
(as explicitly achieved in the finite-dimensional case \cite{GG}). However, actually,
it is more convenient to work with smooth scale functions, as outlined in Remark \ref{supporto} -- and as first used
by Gentile and Mastropietro \cite{GM1} -- both because all bounds at fixed $\om$ are derived straightforwardly and because it is more suited
to deal with when studying the dependence on $\om$ of the torus. Here and in the following, 
to simplify the discussion, we use step scale functions and neglect all contributions in which
the derivatives act on the scale functions, by referring to refs. \cite{G1} and \cite{CGP},
for the finite- and infinite-dimensional cases, respectively.
}
\end{rmk}
%%%%%%%%%%%%%%%%%%%%%%%%%%%%%%%%%%%%%%%%%%%%%%%%%%%%%%%%%%%%%%%%%%%%%%%%%%

If we could prove that $ \matW_T(0) = \del_x  \matW_T(0) = 0$ and, moreover, that $\|R_T (\om\cdot\nu_{\ell'_T})\|_{\rm op}\le C^{|V(T)|}$,
then inserting \eqref{questo} into \eqref{comelovedo} would lead to
\begin{equation}\label{codesto}
\|\Val(\vartheta_{\ell_T})\| \le C^{k(T)}\|\Val(\vartheta_{\ell'_T})\|,
\end{equation}
for a suitable constant $C$, that is we would be able to compensate the propagator of the resonant line $\ell_T$,
thanks to the \emph{gain factor} $(\om\cdot\nu_{\ell'_T})^2$
and we would be able to obtain a bound like \eqref{magara};
unfortunately, as Example \ref{aiuto} shows, in general a bound like \eqref{codesto} does not hold since
there are RCs $T$ such that both $\matW_T(0)$ and $\partial_x \matW_T(0)$ do not vanish.
However, as already pointed out in Remark \ref{modulofuori}, in order to prove Proposition \ref{esiste1},
the bound \eqref{magara} is much stronger than what we really need, since there can be cancellations among
the values of the trees. In the same spirit, we can imagine to collect together all trees in which any RC is replaced
by some other RC with the same external lines and the same number of nodes. This can be achieved in the following way.

Let $\Theta^{(k)}_\nu(n)$ denote the set of trees $\vartheta\in\Theta^{(k)}_\nu$ such that $n_\ell <n$ for all $\ell\in L(\vartheta)$.
Any RC $T$, with $k(T)=k$, $\und{n}_T=n$ and $\nu_{\ell_T}=\nu_{\ell'_T}=\nu$, can be imagined to be constructed
from a tree $\vartheta\in\Theta^{(k)}_0(n)$ by fixing a node $w\in N(\vartheta)$ and attaching to $w$ a line $\ell'$ carrying momentum $\nu$.
As a consequence of such an operation,
\begin{enumerate}[topsep=.5ex]
\itemsep0em
\item $\ell_\vartheta$ and $\ell'$ becomes the exiting line $\ell_T$ and the entering line $\ell_T'$ of $T$;
\item $N(T)=N(\vartheta)$, $L(T)=L(\vartheta)\setminus \ell_\vartheta$;
\item the branching label of $w$ increases by 1 and the corresponding node factor \eqref{nodo} becomes
\[
\calF_w = \frac{1}{(p_w+1)!} f_{\nu_w} (\ii \nu_w)^{p_w+2} ,
\]
while the branching labels and the node factors of all the other nodes $v\in N(\vartheta)\setminus\{w\}$ remain unchanged;
\item if $\nu_\ell^0$ denoted the momentum of any line $\ell\in L(\vartheta)$,
then the momentum of each line $\ell \in \calP(r,w)$ is modified from $\nu_\ell^0$ into $\nu_\ell^0+\nu$,
while the momentum of each line $\ell \not\in \calP(r,w)$ is still $\nu_\ell^0$;
\item the line $\ell_T$ has the same momentum $\nu$ as the entering line $\ell'$.
\end{enumerate}
We define the set $\SSSS^{(k)}(n)$ of the RCs $T'$ with $k(T')=k$ and $\ol{n}_{T'}<n$
as the set of RCs which can be obtained by considering all trees $\vartheta\in \Theta^{(k)}_0(n)$ and attaching 
a line to any node $v\in N(\vartheta)$. 

%%%%%%%%%%%%%%%%%%%%%%%%%%%%%%%%%%%%%%%%%%%%%%%%%%%%%%%%%%%%%%%%%%%%%%%%%%
\begin{rmk} \label{pp+1a}
\emph{
Let $T$ be a RC obtained from the tree $\vartheta$ by attaching the line $\ell'$ to the node $w\in N(\vartheta)$. 
If we compare $\Val(\vartheta)$, as given by \eqref{val}, with $\matW(T)$, as given by \eqref{valS} with $S=T$,
the branching label $p_w$ has a different value in the two expressions: if $p_w=p$ denotes the branching label of $w\in N(\vartheta)$,
then, when $w=v_T'$ is seen as a node in $N(T)$, we have $p_w=p+1$.
}
\end{rmk}
%%%%%%%%%%%%%%%%%%%%%%%%%%%%%%%%%%%%%%%%%%%%%%%%%%%%%%%%%%%%%%%%%%%%%%%%%%

%%%%%%%%%%%%%%%%%%%%%%%%%%%%%%%%%%%%%%%%%%%%%%%%%%%%%%%%%%%%%%%%%%%%%%%%%%
\begin{rmk} \label{pp+1b}
\emph{
For any fixed $w\in N(\vartheta)$ there are $p_{w}+1$ RCs that can be created by following the procedure above,
because the line $\ell'$ can be inserted in $p_{w}+1$ ways w.r.t.~the existing lines entering $w$.
This means that, if we sum together the values of the RCs which differ only because of the position of the attached line, 
since such values are equal to each other, we obtain a factor $p_w+1$ which, multiplied by the combinatorial factor $1/(p_w+1)!$ 
returns the factor $1/p_w!$ appearing in $\Val(\vartheta)$.
}
\end{rmk}
%%%%%%%%%%%%%%%%%%%%%%%%%%%%%%%%%%%%%%%%%%%%%%%%%%%%%%%%%%%%%%%%%%%%%%%%%%

%%%%%%%%%%%%%%%%%%%%%%%%%%%%%%%%%%%%%%%%%%%%%%%%%%%%%%%%%%%%%%%%%%%%%%%%%%
\begin{ex} \label{inserting}
\emph{
As an example illustrating what observed in Remark \ref{pp+1b} see Figure \ref{insert1}.
By attaching a line to a node $v\in N(\vartheta_1)$ we obtain the RCs in Figure \ref{sec2}, while
by attaching a line to a node $v\in N(\vartheta_2)$ we obtain the RCs in Figure \ref{insert2}).
In particular, in the latter case, there are three RCs which are obtained by attaching the line $\ell'$ to the node with mode label $\nu_1$,
and all of them have the same value.
}
\end{ex}
%%%%%%%%%%%%%%%%%%%%%%%%%%%%%%%%%%%%%%%%%%%%%%%%%%%%%%%%%%%%%%%%%%%%%%%% 

%%%%%%%%%%%%%%%%%%%%%%%%%%%%%%%%%%%%%%%%%%%%%%%%%%%%%%%%%%%%%%%%%%%%%%%% 
% Figure RC
%%%%%%%%%%%%%%%%%%%%%%%%%%%%%%%%%%%%%%%%%%%%%%%%%%%%%%%%%%%%%%%%%%%%%%%% 
\begin{figure}[ht] 
\centering 
\ins{052pt}{-36.5pt}{$\vartheta_1=$}
\ins{254pt}{-36.5pt}{$\vartheta_2=$}
\ins{131pt}{-46pt}{$\nu_1$}
\ins{182pt}{-46pt}{$\nu_2$}
\ins{326pt}{-46pt}{$\nu_1$}
\ins{379pt}{-10pt}{$\nu_2$}
\ins{379pt}{-84pt}{$\nu_3$}
\null\hspace{.5cm}
\includegraphics[width=1.65in]{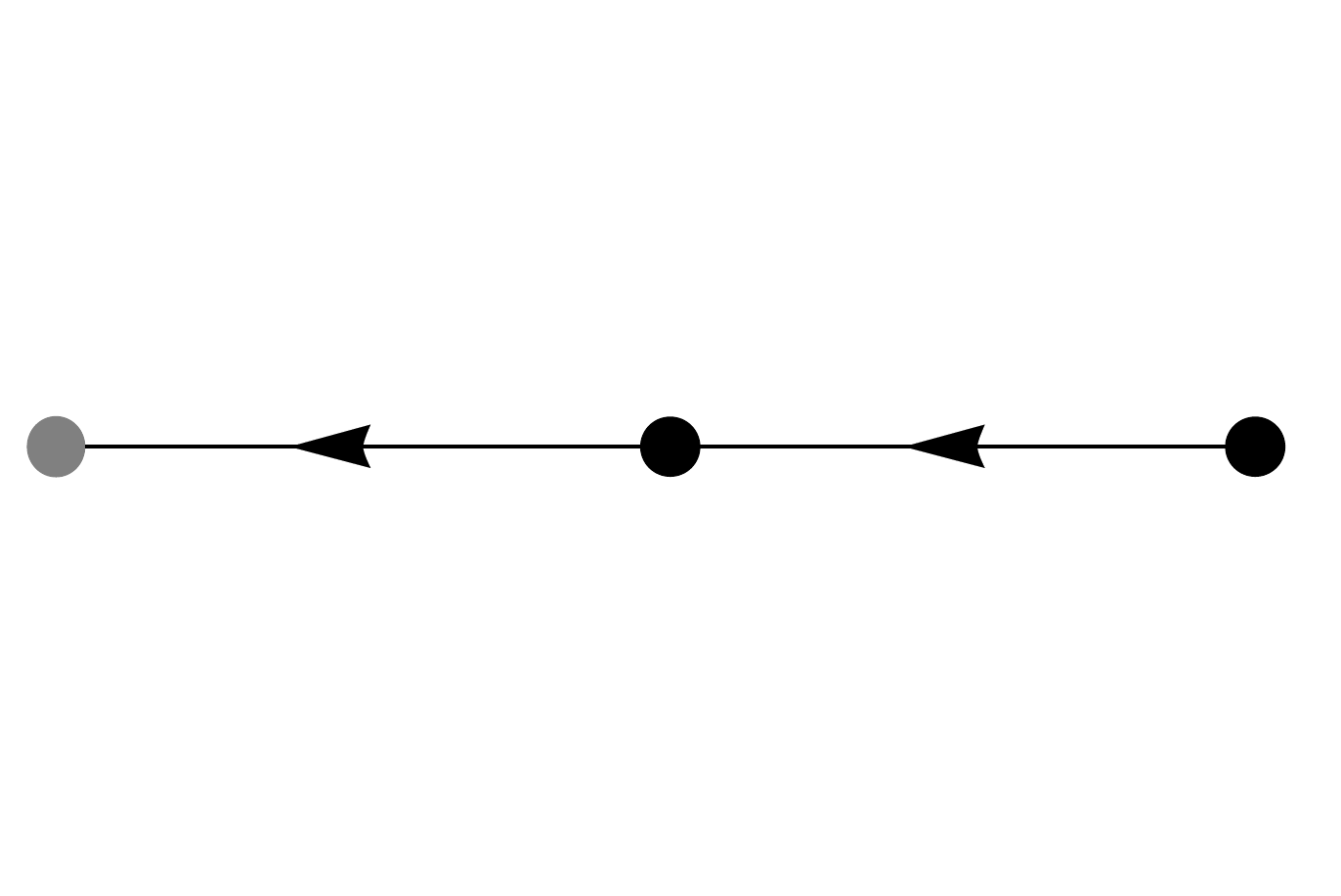}
\hspace{1.cm}
\includegraphics[width=2.2in]{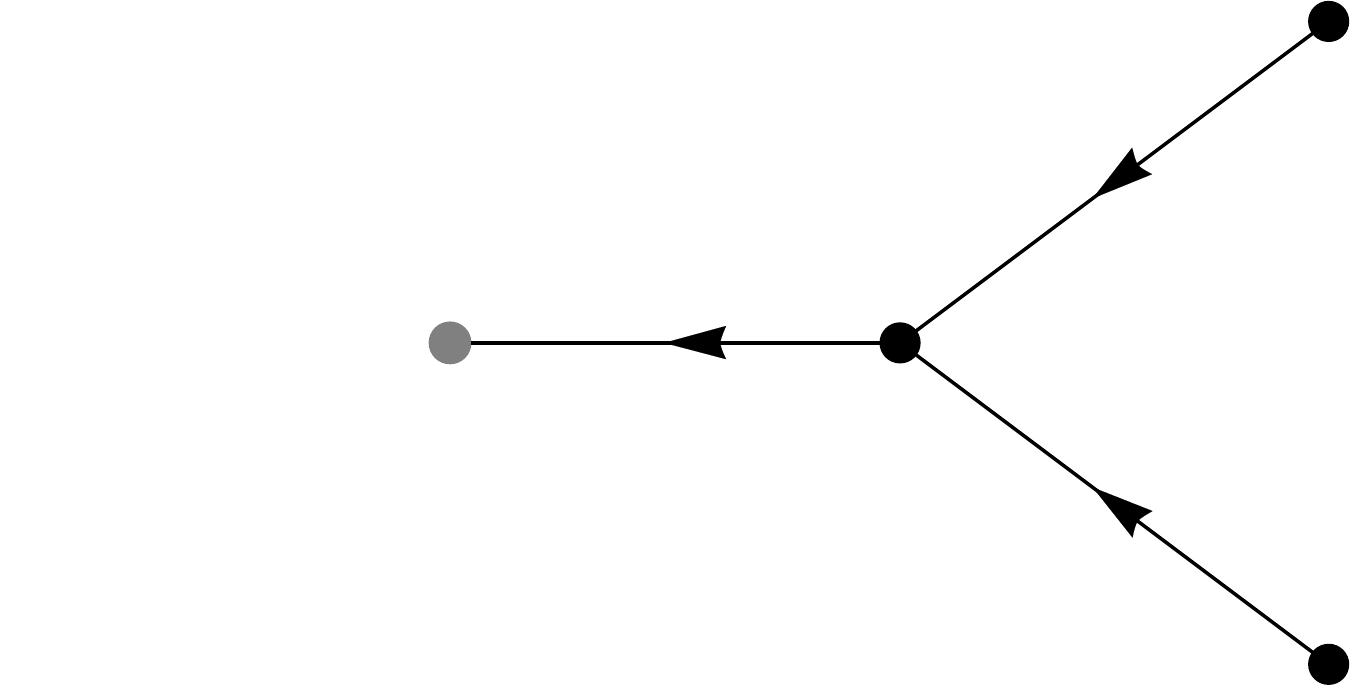} 
\caption{Trees $\vartheta_1$ and $\vartheta_2$ of order 2 and 3, respectively.}
\label{insert1} 
\end{figure} 
%%%%%%%%%%%%%%%%%%%%%%%%%%%%%%%%%%%%%%%%%%%%%%%%%%%%%%%%%%%%%%%%%%%%%%%%

%%%%%%%%%%%%%%%%%%%%%%%%%%%%%%%%%%%%%%%%%%%%%%%%%%%%%%%%%%%%%%%%%%%%%%%% 
% Figure RC
%%%%%%%%%%%%%%%%%%%%%%%%%%%%%%%%%%%%%%%%%%%%%%%%%%%%%%%%%%%%%%%%%%%%%%%% 
\begin{figure}[ht] 
\centering 
\ins{031pt}{-50pt}{$\nu_1$}
\ins{114pt}{-50pt}{$\nu_1$}
\ins{209pt}{-50pt}{$\nu_1$}
\ins{293pt}{-50pt}{$\nu_1$}
\ins{377pt}{-50pt}{$\nu_1$}
\ins{061pt}{-35pt}{$\nu_2$}
\ins{144pt}{-35pt}{$\nu_2$}
\ins{239pt}{-35pt}{$\nu_2$}
\ins{323pt}{-35pt}{$\nu_2$}
\ins{407pt}{-35pt}{$\nu_2$}
\ins{053pt}{-63pt}{$\nu_3$}
\ins{136pt}{-63pt}{$\nu_3$}
\ins{233pt}{-63pt}{$\nu_3$}
\ins{315pt}{-63pt}{$\nu_3$}
\ins{399pt}{-63pt}{$\nu_3$}
\includegraphics[width=1.3in]{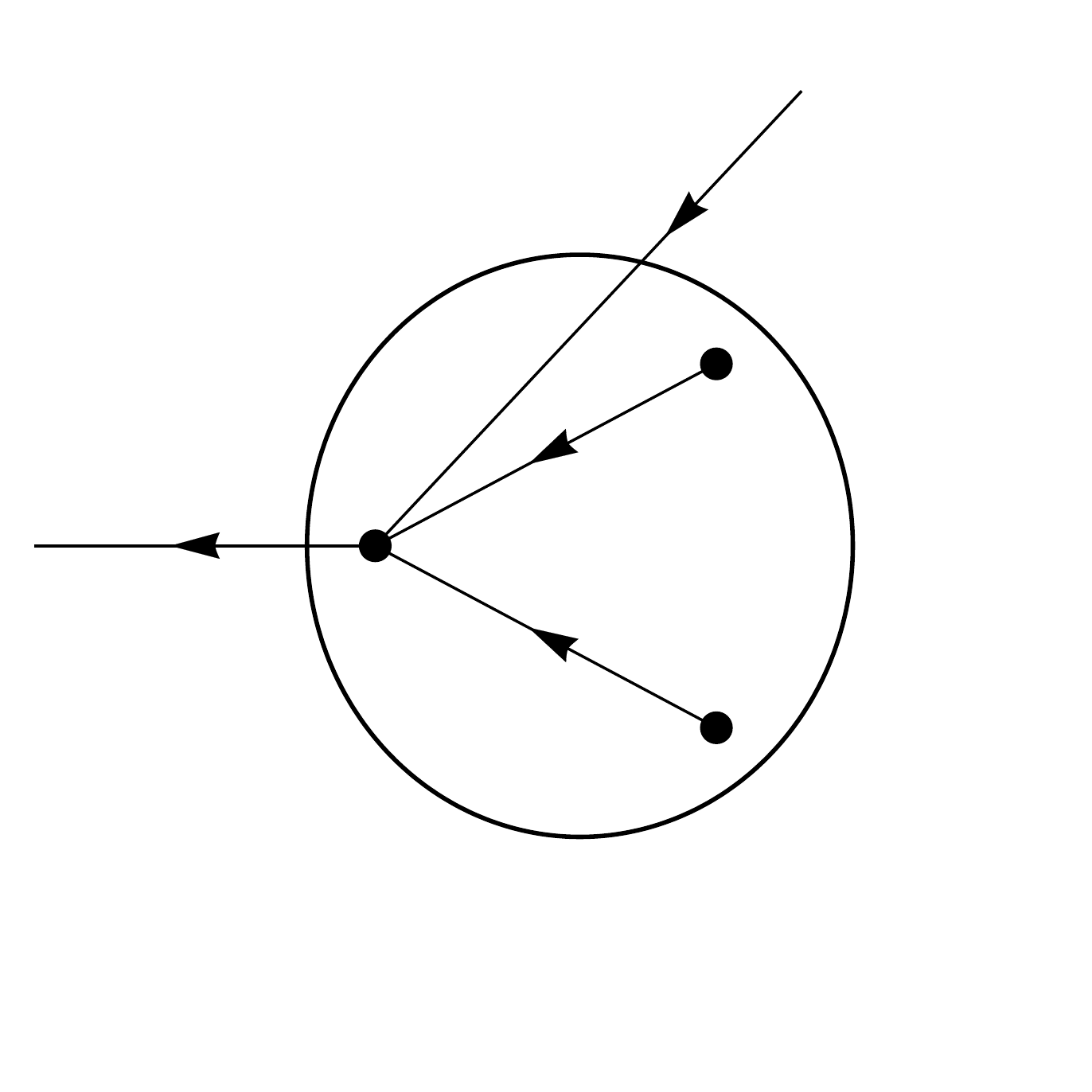}
\hspace{-.6cm}
\includegraphics[width=1.3in]{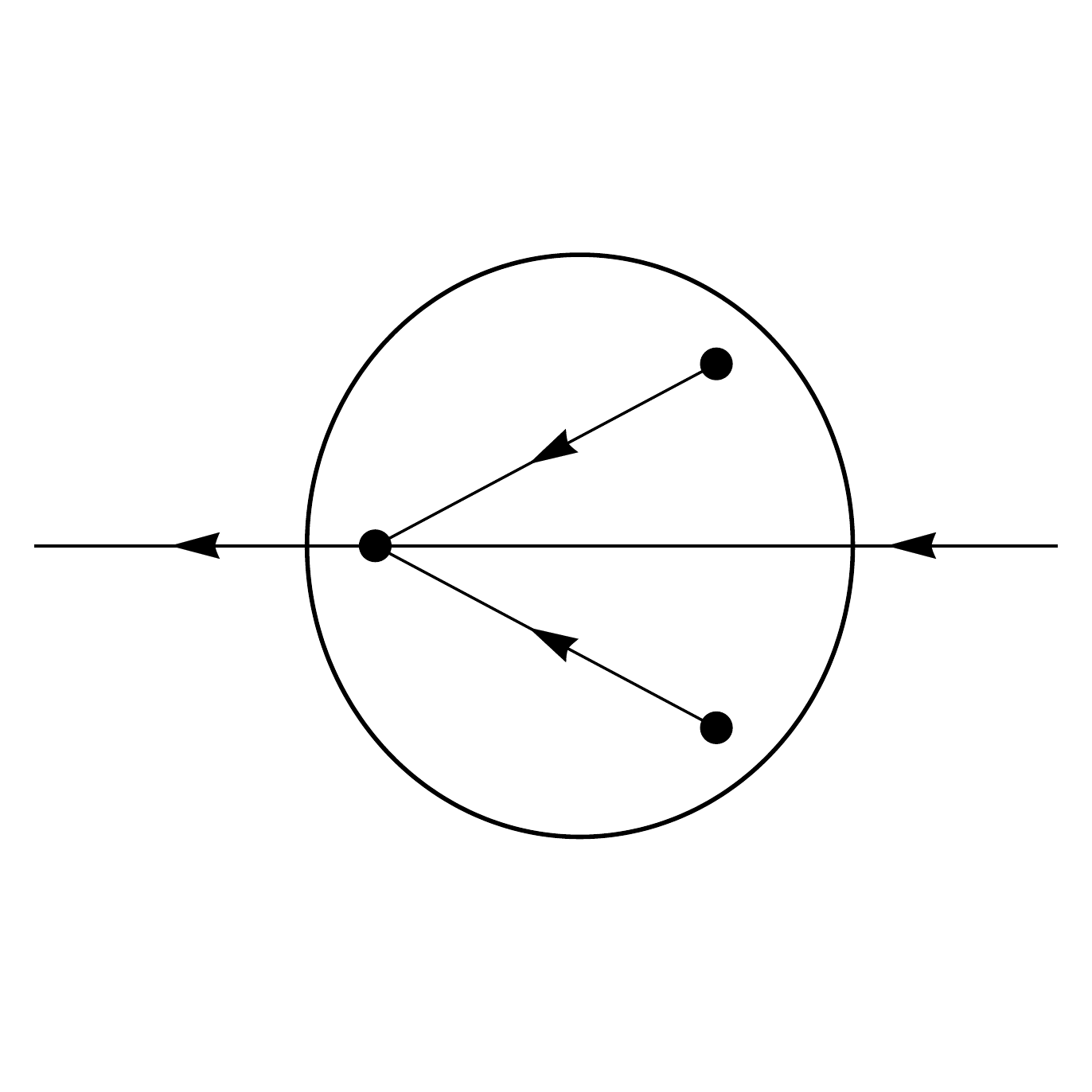}
\hspace{-.2cm}
\includegraphics[width=1.3in]{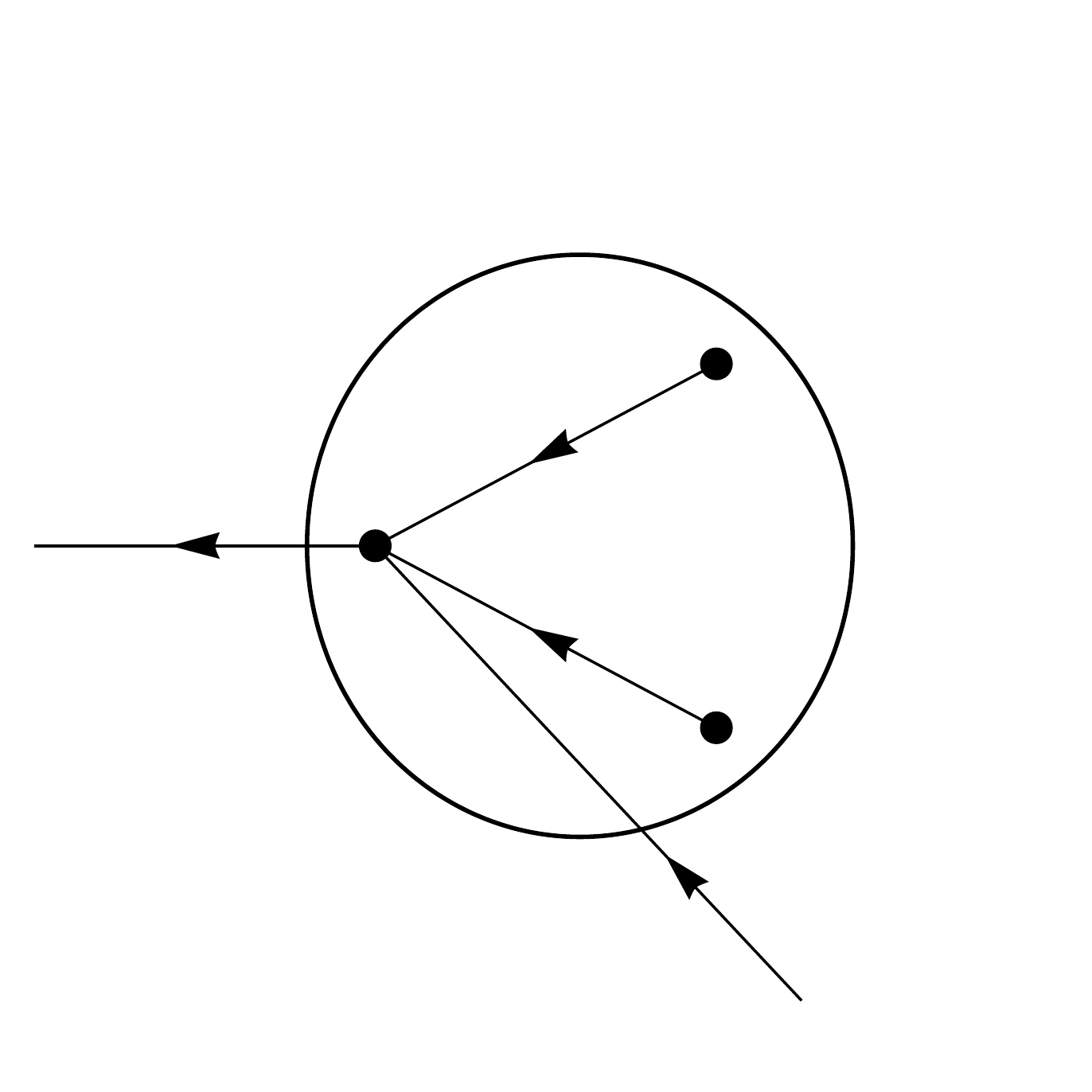}
\hspace{-.6cm}
\includegraphics[width=1.3in]{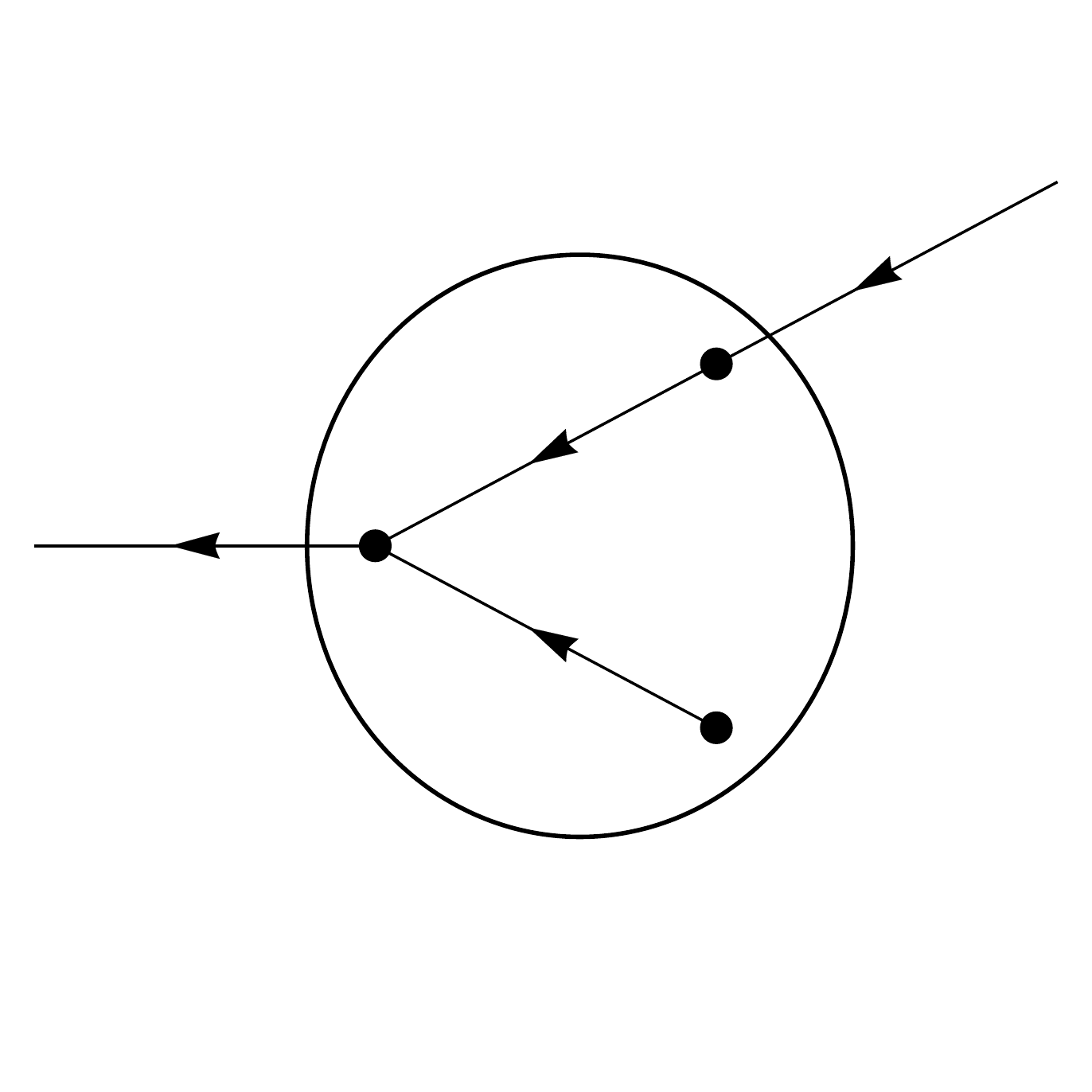}
\hspace{-.6cm}
\includegraphics[width=1.3in]{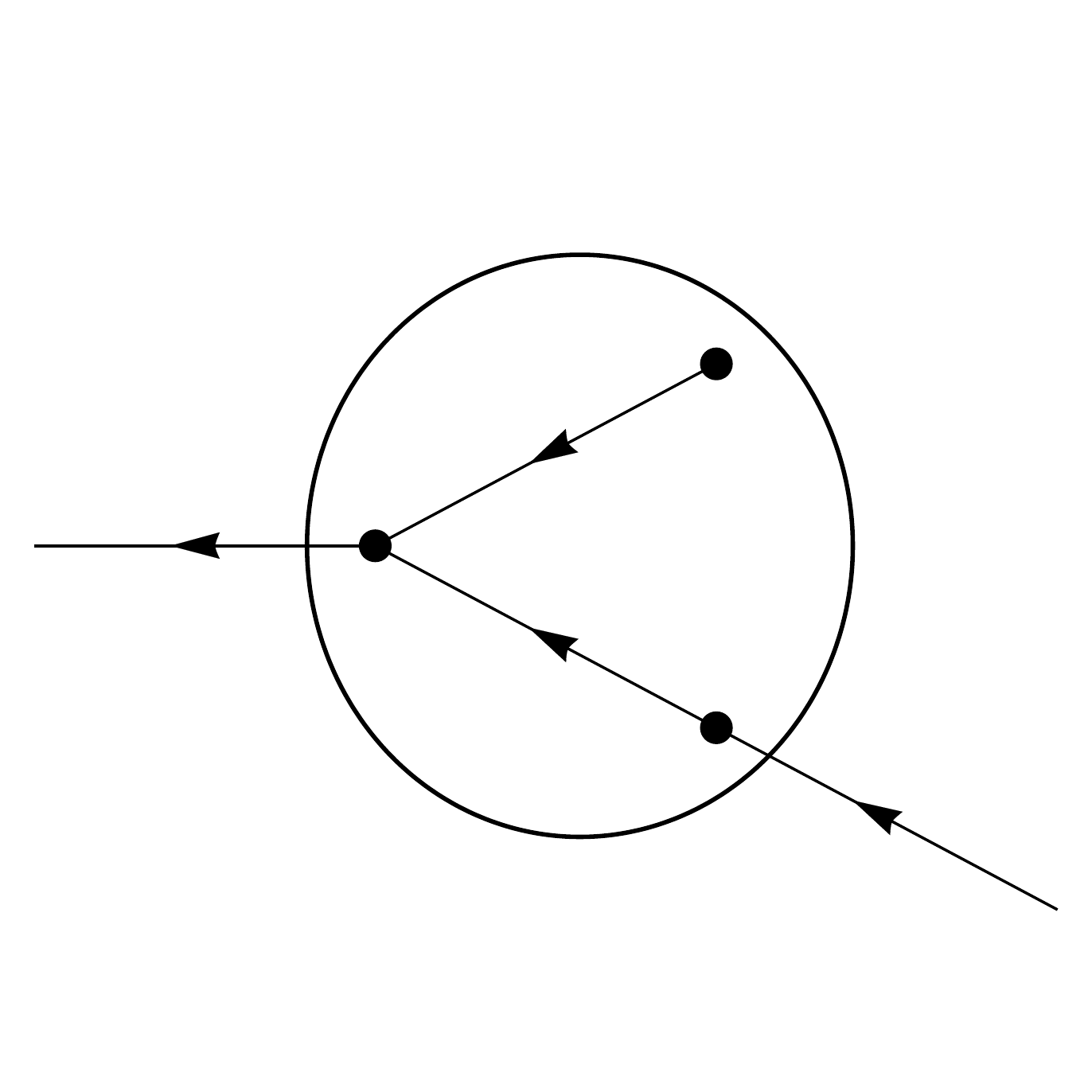}
\caption{Resonant clusters obtained by attaching a line to a node of the tree $\vartheta_2$ in Figure \ref{sec2}.}
\label{insert2} 
\end{figure} 
%%%%%%%%%%%%%%%%%%%%%%%%%%%%%%%%%%%%%%%%%%%%%%%%%%%%%%%%%%%%%%%%%%%%%%%%

With this in mind, given a tree $\vartheta$ containing a RC $T$,
if we sum together all trees which can be obtained from $\vartheta$ by replacing $T$ with any other RC of order $k(T)$,
by replacing $T$ with any RC $T'$ such that  $k(T')=k(T)=k$ and $\ol{n}_{T'}<\und{n}_T=n$,
then, using the notation in Remark \ref{piS}, we find
\begin{equation}\label{sommo}
\sum_{\vartheta'\in\gotF(\vartheta,T)}\Val(\vartheta) = \calG_{\ell_\vartheta} 
\matW(\pi(\vartheta,T)) [\calG_{\ell_T} \matM_{n}^{(k)}(\om\cdot\nu_{\ell'_T})[\Val(\vartheta_{\ell'_T})] ] ,
\end{equation}
where
\begin{equation}\label{emme}
\matM^{(k)}_n(x) := \sum_{T'\in \SSSS^{(k)}(n)}\matW_{T'}(x) .
\end{equation}
%

%%%%%%%%%%%%%%%%%%%%%%%%%%%%%%%%%%%%%%%%%%%%%%%%%%%%%%%%%%%%%%%%%%%%%%%%%%
\begin{rmk} \label{Skn-example1}
\emph{
For instance, if $T$ is the RC on the left of Figure \ref{sec2}, so that $k(T)=2$, then $\SSSS^{(2)}(n)$ contains all the RCs which,
neglecting the labels, appear as one of the three graphs in Figure \ref{sec2},  and, when the labels are taken into account,
are obtained by summing over all possible assignments of the mode labels $\nu_1$ and $\nu_2$
satisfying the constraints in Definition \ref{RC}, with $\und{n}_T=n$.
}
\end{rmk}
%%%%%%%%%%%%%%%%%%%%%%%%%%%%%%%%%%%%%%%%%%%%%%%%%%%%%%%%%%%%%%%%%%%%%%%%%%

%%%%%%%%%%%%%%%%%%%%%%%%%%%%%%%%%%%%%%%%%%%%%%%%%%%%%%%%%%%%%%%%%%%%%%%%%%
\begin{ex} \label{SKn-example2}
\emph{
If $T$ is a RC in $\vartheta$, the subgraph $\pi(\vartheta,T)$ can still be a RC or a chain of RCs,
so that the decomposition \eqref{sommo} can be iterated.
This happens, for instance, if we consider the tree $\vartheta$ in Example \ref{aiuto}
(refer to Figure \ref{catenabis}).
Set $p=(k+1)/2$ and call $\ell_1,\ldots,\ell_p$ the lines with momentum $\nu$ ordered from left to right,
$T_1,\ldots,T_{p-1}$ the RCs such that $\ell_{T_i}=\ell_i$ and $\ell_{T_i}'=\ell_{i+1}$ for $i=1,\ldots,p-1$,
and $\vartheta_1,\ldots,\vartheta_p$ the subtrees with root lines $\ell_1,\ldots,\ell_p$, respectively,
so that $\ell_1=\ell_\vartheta$ and $\vartheta_1=\vartheta$ and,
with the notation in Remark \ref{piS}, one has $\pi(\vartheta,T_{i+1})= \vartheta \setminus \vartheta_{i+1}$
and $\pi(\vartheta_i,T_{i+1})=T_i$ for all $i=1,\ldots,p-1$. We have already stressed (see Example \ref{RC1})
that each $T_i$, for $i=1,\ldots,p-1$, is as in the middle of Figure \ref{sec2}.
Then, analogously to \eqref{comelovedo}, we have
\[
\Val(\vartheta)=\calG_{\ell_1} \matW(T_1) [\calG_{\ell_2} \matW(T_2) [ \ldots [\calG_{\ell_p} [\ii \nu f_\nu]]]] ,
\]
where, for all $i=1,\ldots,p-1$,
\begin{equation} \label{WT1}
\matW(T_i)[z] = 
\frac{1}{2} \frac{(\ii \nu_0) |\nu_0|^2 |f_{\nu_0}|^2}{(\om\cdot\nu_0)^2}  \ii\nu_0 \cdot z 
\end{equation}
assumes the same value and equals $\matW_{T_i}(0)[z]$; in other words,
$\matW_{T_i}(x)$ is the same for all $i=1,\ldots,p-1$ and is independent of $x$.
If we replace a RC $T_i$, for some $i$, with the RC $T_i'$ as on the left of Figure \ref{sec2} (see Remark \ref{Skn-example1}),
then $\matW(T_i)[z]$ is substituted with
\begin{equation} \label{WT2}
\matW(T_i')[z] = 
\frac{(\ii \nu_0) |\nu_0|^2 |f_{\nu_0}|^2}{(\om\cdot\nu_0 + \om\cdot\nn)^2}  (-\ii\nu_0 \cdot z) , 
\end{equation}
which shows that $\matW_{T_i'}(x)$ depends explicitly on $x$ in such a case,
while, if the RC $T_i'$ replacing $T_i$ is as on the right of Figure \ref{sec2},
then $\matW(T_i)[z]$ is substituted with
\begin{equation} \label{WT3}
\matW(T_i')[z] = 
\frac{1}{2} \frac{(\ii \nu_0) |\nu_0|^2 |f_{\nu_0}|^2}{(\om\cdot\nu_0)^2}  \ii\nu_0 \cdot z = \matW(T_i)[z] .
\end{equation}
In conclusion, if we take into account Remark \ref{Skn-example1}, we obtain
\begin{equation}\label{emme2}
\begin{aligned}
\matM^{(2)}_n(x)[z] & :=
\sum_{\substack{ \nu_0 \in \ZZZ^\ZZZ_f\setminus\{0\} \\ n'<n}} 
\left( \frac{1}{2} (\ii \nu_0) |\nu_0|^2 |f_{\nu_0}|^2 \frac{\Psi_{n'}(\om\cdot\nu_0)}{(\om\cdot\nu_0)^2} \right. \\
& - 
\left. (\ii \nu_0) |\nu_0|^2 |f_{\nu_0}|^2 \frac{\Psi_{n'}(\om\cdot\nu_0+x)}{(\om\cdot\nu_0 + x)^2} +
\frac{1}{2} (\ii \nu_0) |\nu_0|^2 |f_{\nu_0}|^2
\frac{\Psi_{n'}(\om\cdot\nu_0)}{(\om\cdot\nu_0)^2}  \right) \ii\nu_0 \cdot z .
\end{aligned}
\end{equation}
For each fixed $\nu_0$ the sum in \eqref{emme2} vanish when $x=0$, so that $\matM^{(2)}_n(0)= 0$.
One also easily checks, provided we ignore the contributions of derivatives of the scale functions
(see Remark \ref{vanishing}), that $\partial_x \matM^{(2)}_n(0)= 0$.
}
\end{ex}
%%%%%%%%%%%%%%%%%%%%%%%%%%%%%%%%%%%%%%%%%%%%%%%%%%%%%%%%%%%%%%%%%%%%%%%%%%

In general, since we are interested in the case in which $\om\cdot\nu_{\ell'_T}$ is small, we may write
\begin{equation}\label{taylor}
\matM_{n}^{(k)} (\om\cdot\nu_{\ell'_T}) = 
\matM_{n}^{(k)} (0) + (\om\cdot\nu_{\ell'_T}) \,\del_x\matM_{n}^{(k)} (0) + (\om\cdot\nu_{\ell'_T})^2\matR^{(k)}_n (\om\cdot\nu_{\ell'_T}) ,
\end{equation}
with $\matR^{(k)}_n\!:\ell^\io(\CCC)\to\ell^\io(\CCC)$ a continuous linear operator and
\[
\del_x\matM_{n}^{(k)} (0) = \sum_{T'\in \SSSS^{(k)}(n)} \del_x \matW_{T'}(0) ,
\]
where, according to \eqref{valfun},
\begin{equation} \label{valfunder}
\del_x \matW_T(0) =
\frac{1}{p_{v_T}!}\ii \nu_{v_T} f_{\nu_{v_T}} \Biggl( \prod_{v\in N(\vartheta)\setminus \{v_T\}} 
\frac{1}{p_{v}!} (\ii \nu_{\pi(v)}\cdot \ii \nu_v ) f_{\nu_v} \Biggl) 
\del_x \!\!\! 
\prod_{\ell\in L(T)} \calG_{n_\ell}(\x_\ell(x))) 
\Big|_{x=0} \!\!\!\!\! \ii \nu_{v'_T} ,
\end{equation}
with 
\[
\del_x \!\!\! 
\prod_{\ell\in L(T)} \!\!\! \calG_{n_\ell}(\x_\ell(x))) 
\Big|_{x=0}  \!\!\!\! = \!\!\!\!
\sum_{v\in N(T)} \!\!\! \left. \del_x \calG_{n_\ell}(\x_\ell(x))) \right|_{x=0} \!\!\!\!\!\!\!\!
\prod_{\ell' \in L(T)\setminus\{\ell\}} \!\!\!\!\!\!\! \calG_{n_{\ell'}}(\x_{\ell'(0)})) , \qquad
\del_x \calG_{n_\ell}(\x_\ell(0))) = - \frac{2 \Psi_{n_\ell}(\om\cdot\nu_\ell^0)}{(\om\cdot\nu_\ell^0)^3} .
\]

Then the bound \eqref{hovinto} follows if we are able to prove that $ \matM_{n}^{(k)} (0) = \del_x\matM_{n}^{(k)} (0) = 0$
and that $\|\matR^{(k)}_n (\om\cdot\nu_{\ell'_T}) \|_{\rm op}<C^k$ for some positive constant $C$.

%%%%%%%%%%%%%%%%%%%%%%%%%%%%%%%%%%%%%%%%%%%%%%%%%%%%%%%%%%%%%%%%%%%%%%%%%%
\begin{rmk} \label{sottolineiamo anche questo}
\emph{
The cancellations $ \matM_{n}^{(k)} (0) = 0$ and $\del_x\matM_{n}^{(k)} (0) = 0$,
that we verified explicitly for the RCs considered in Example \ref{Skn-example1},
actually hold for all $k\ge 1$ and $n\ge 0$.
They follow from deep structural properties of the equations of motion \eqref{lagra}, as the next argument shows.
}
\end{rmk}
%%%%%%%%%%%%%%%%%%%%%%%%%%%%%%%%%%%%%%%%%%%%%%%%%%%%%%%%%%%%%%%%%%%%%%%%%%

Recalling the translation covariance discussed in Remark \ref{invarianza}, if we define
\begin{subequations} \label{senta}
\begin{align}
U^{(k)}_n(\f + \f_0;\e) & := \sum_{h=1}^k\e^h  \sum_{\nu\in\ZZZ^\ZZZ_f\setminus\{0\}} e^{\ii\nu\cdot\f}
\sum_{\vartheta\in\Theta^{(h)}_\nu(n)} \Val(\vartheta;\f_0), \phantom{\sum_N^N}
\label{sentaa} \\
G^{(k)}_n(\f_0) & := \sum_{\vartheta\in\Theta^{(k)}_0(n)} \Val(\vartheta;\f_0) , \phantom{\sum_N^N}
\label{sentab}
\end{align}
\end{subequations}
where we have set
\begin{equation}\label{Fv-f0}
\calF_v(\f_0):=\calF_v e^{\ii \nu_v \cdot \f_0} = \frac{1}{p_v!} f_{\nu_v} e^{\ii \nu_v \cdot \f_0} (\ii \nu_v)^{p_v+1} =
\frac{1}{p_v!} f_{\nu_v} (\f_0) \, (\ii \nu_v)^{p_v+1} 
\end{equation}
and
\begin{equation}\label{val-f0}
\Val(\vartheta;\f_0) := \Bigg(\prod_{v\in N(\vartheta)} \!\! \calF_v (\f_0) \Bigg) 
\Bigg(\prod_{\ell\in L(\vartheta)} \!\! \calG_\ell\Bigg) ,
\end{equation}
then for all $k\ge1$ and all $\nu\in\ZZZ^\ZZZ_f$, we have
$\Val(\vartheta;\f_0)=e^{\ii \nu \cdot \f_0} \Val(\vartheta)$ for all $\vartheta \in \Theta^{(k)}_\nu(n)$.
Moreover, by Lemma \ref{metalemma}, for all $k\ge 1$ and all $\f_0\in\TTT^\ZZZ$, we have
\begin{equation}\label{mediato}
G^{(k)}_n(\f_0)= [\e \del f(\f+\f_0+U_n^{(k)}(\f+\f_0;\e))]^{(k)}_0 = 0.
\end{equation}

The first cancellation is a consequence of the following result.

%%%%%%%%%%%%%%%%%%%%%%%%%%%%%%%%%%%%%%%%%%%%%%%%%%%%%%%%%%%%%%%%%%%%%%%%%%
\begin{lemma}\label{emmeladerivatadigi}
One has
%
%\[
$\matM_n^{(k)}(0) = \del_{\f_0} G^{(k)}_n(\f_0)$.
%\]
%
\end{lemma}
%%%%%%%%%%%%%%%%%%%%%%%%%%%%%%%%%%%%%%%%%%%%%%%%%%%%%%%%%%%%%%%%%%%%%%%%%%

%%%%%%%%%%%%%%%%%%%%%%%%%%%%%%%%%%%%%%%%%%%%%%%%%%%%%%%%%%%%%%%%%%%%%%%%%%
\prova
Since $\Val(\vartheta;\f_0)$ depends on $\f_0$ only through the node factors \eqref{Fv-f0}, 
for $z\in \ell^\io(\RRR)$ we have
\begin{equation}\label{valB}
\del_{\f_0} \! \Val(\vartheta;\f_0)\cdot z := \mathtt{d}_{\f_0} \!\Val(\vartheta;\f_0)[z]  := 
\sum_{v\in N(\vartheta)}(\ii\nu_v\cdot z )\Bigg(\prod_{v\in N(\vartheta)} \calF_v (\f_0) \Bigg)
\Bigg(\prod_{\ell\in L(\vartheta)} \calG_\ell\Bigg).
\end{equation} 
Now, we note that for $\nu\in \ZZZ^\ZZZ_f$ and all $z, z_1,\dots z_p\in\CCC^\ZZZ$
\[
\sum_{r=1}^{p+1}\nu^{p+2}[z_1,\dots,z_{r-1},z, z_{r}\dots,z_p] =(p+1)(\nu\cdot z )\nu^{p+1}[z_1,\dots,z_p]
\]
Given any $\vartheta\in \Theta^{(k)}_0(n)$ and any $v\in N(\vartheta)$,
we have (see Remark \ref{piSbis})
\[
\Bigg( \prod_{\substack{ w \in N(\vartheta) \\ w \preceq v }} \calF_w (\f_0) \Bigg)
\Bigg( \prod_{\substack{ \ell\in L(\vartheta) \\ \ell \preceq v}} \calG_\ell\Bigg) = 
\frac{1}{p_v !}  f_{\nu_v}(\f_0) (\ii\nu_v)^{p_v+1} [\Val(\vartheta_1),\ldots,\Val(\vartheta_{p_v})] ,
\]
where $\vartheta_1,\ldots,\vartheta_{p_v}$ denote the subtrees of $\vartheta$ whose root lines enter $v$.

Then we obtain
\begin{equation}\label{derivo}
\begin{aligned}
&  ( \mathtt{d}_{\f_0} \calF_v) [\Val(\vartheta_1),\ldots, \Val(\vartheta_{p_v})] [z] \phantom{\sum_{r=1}^{p_v}}  \\
& \qquad\qquad 
=  \frac{1}{p_v !}  f_{\nu_v}(\f_0)  (\ii \nu_v\cdot z)(\ii\nu_v)^{p_v+1} [\Val(\vartheta_1),\ldots,\Val(\vartheta_{p_v})]
\phantom{\sum_{r=1}} \\
& \qquad\qquad 
= (p_v+1)  \frac{1}{(p_v+1) !}  f_{\nu_v}(\f_0) (\ii\nu_v)^{p_v+2}  [\Val(\vartheta_1),\ldots,\Val(\vartheta_{p_v}),z] \phantom{\sum_{r=1}^{p_v}} \\
 & \qquad\qquad 
 = \frac{1}{(p_v+1) !}  f_{\nu_v}(\f_0) \sum_{r=1}^{p_v+1} (\ii\nu)^{p_v+2}[\Val(\vartheta_1),\dots,\Val(\vartheta_{r-1}),z, \Val(\vartheta_{r}),\dots,\Val(\vartheta_{p_v})] .
 \end{aligned}
\end{equation}
We can regard each summand in \eqref{derivo} as the node factor of a RC obtained from $\vartheta$
by attaching to the node $v$ an extra line $\ell$ with momentum $\nu_\ell=0$ (see Remark \ref{pp+1a}).
The value of each such RCs gives a contribution to $\matM_n^{(k)}(0)[z]$. 
If we consider all possible trees $\vartheta\in\Theta^{(k)}_0(n)$ and apply \eqref{derivo} to all nodes $v\in\vartheta$,
we produce all possible RC $T\in\SSSS^{(k)}(n)$. Therefore, \eqref{valB} equals $\matM_n^{(k)}(0)[z]$. 
Hence the assertion follows.
\EP
%%%%%%%%%%%%%%%%%%%%%%%%%%%%%%%%%%%%%%%%%%%%%%%%%%%%%%%%%%%%%%%%%%%%%%%%%%

Then we deduce the following result.

%%%%%%%%%%%%%%%%%%%%%%%%%%%%%%%%%%%%%%%%%%%%%%%%%%%%%%%%%%%%%%%%%%%%%%%%%%
\begin{coro}\label{cancellazione}
One has $\matM_n^{(k)}(0)=0$.
\end{coro}
%%%%%%%%%%%%%%%%%%%%%%%%%%%%%%%%%%%%%%%%%%%%%%%%%%%%%%%%%%%%%%%%%%%%%%%%%%

%%%%%%%%%%%%%%%%%%%%%%%%%%%%%%%%%%%%%%%%%%%%%%%%%%%%%%%%%%%%%%%%%%%%%%%%%%
\prova
Combine Lemma \ref{emmeladerivatadigi} with \eqref{mediato}.
\EP
%%%%%%%%%%%%%%%%%%%%%%%%%%%%%%%%%%%%%%%%%%%%%%%%%%%%%%%%%%%%%%%%%%%%%%%%%%

The second cancellation is the content of the following result.

%%%%%%%%%%%%%%%%%%%%%%%%%%%%%%%%%%%%%%%%%%%%%%%%%%%%%%%%%%%%%%%%%%%%%%%%%%
\begin{lemma}\label{reverse}
One has $\del_x \matM_n^{(k)}(0)= 0$.
\end{lemma}
%%%%%%%%%%%%%%%%%%%%%%%%%%%%%%%%%%%%%%%%%%%%%%%%%%%%%%%%%%%%%%%%%%%%%%%%%%

%%%%%%%%%%%%%%%%%%%%%%%%%%%%%%%%%%%%%%%%%%%%%%%%%%%%%%%%%%%%%%%%%%%%%%%%%%
\prova
Given any RC $T$, set $\calP_T:=\calP(v_{T},v_{T}')$, where $v_T$ and $v'_T$ are the node
which $\ell_T$ exits and $\ell_T'$ enters, respectively.
If $\calP_T=\emptyset$, then $\matW_T(x)=\matW_T(0)$ for all $x$, so that $\del_x\matW_T(0)=0$. 
If, instead, $\calP_T \neq\emptyset$, then $\del_x \matW_T(0)$ is given by \eqref{valfunder}.

For any fixed line $\ell\in\calP_T$, write $N(T)=N_1(T)\sqcup N_2(T)$, with $N_1(T):=\{ v\in N(T) : v\preceq \ell\}$
and $V_2(T):=N(T)\setminus N_1(T)$, and
define $T_1$ as the subgraph of $T$ such that $N(T_1)=N_1(T)$
and $L(T_1)=\{ \ell' \in L(T) : \ell' \prec \ell\}$, and $T_2$ as the subgraph of $T$ such that $N(T_2)=N_2(T)$ and
$L(T_2)=L(T)\setminus (L(T_1)\cup\{\ell\})$ (see Figure \ref{canc2});
by construction one has $L(T)=L(T_1) \sqcup \{\ell\} \sqcup L(T_2)$.
Since $T$ is a RC, by item \ref{item2} in Definition \ref{RC} one has
\begin{equation}\label{pari}
\sum_{v\in N(T_1)}\nu_v + \sum_{v\in N(T_2)} \nu_v=0.
\end{equation}

Now, consider the family $\matF_1(T)$ consisting of all RCs which can be obtained by detaching both $\ell_T$ and $\ell'_T$
and reattaching $\ell_T$ to a node in $V_2(T)$ and $\ell'_T$ to a node in $V_1(T)$.
Consider also the family $\matF_2(T)$ consisting of all RCs which can be obtained by detaching both
$\ell_T$ and $\ell'_T$ and reattaching $\ell_T$ to a node in $V_1(T)$ and $\ell'_T$ to a node in $V_2(T)$
(refer to Examples \ref{effe1} and \ref{effe2} below).
Any $T'\in\matF_1(T)\sqcup\matF_2(T)$ gives rise to a contribution $\matW_{T'}(x)$ to $\matM^{(k)}_n(x)$.

One checks easily, by taking into account Remark \ref{pp+1b}, that
\[
\sum_{T'\in\matF_1(T)} \prod_{v\in N(T')}\calF_v = \sum_{T'\in\matF_2(T)} \prod_{v\in N(T')}\calF_v,
\]
Moreover, since $\calG_n(x)$ is an even function,
because of \eqref{pari}, all the propagators except the differentiated one do not change,
while the differentiated propagator for $T'\in\calF_2(T)$ changes sign w.r.t.~the $T''\in\calF_1(T)$, i.e.
\[
\left. \del_x \calG_{n_\ell}(\om\cdot\nu_\ell^0+x)\right|_{x=0}
\prod_{\substack{ \ell'\in L(T') \\ \ell'\ne\ell}} \calG_{n_{\ell'}}(\om\cdot\nu_{\ell'}^0) 
 = -
\left. \del_x \calG_{n_\ell}(\om\cdot\nu_\ell^0+x)\right|_{x=0}
\prod_{\substack{ \ell'\in L(T'\!') \\ \ell'\ne\ell}} \calG_{n_{\ell'}}(\om\cdot\nu_{\ell'}^0) 
\]
hence the overall sum 
\[
\sum_{T'\in\matF_1(T)} \partial_x \matW_{T'}(x) + \sum_{T'\in\matF_2(T)} \partial_x \matW_{T'}(x) ,
\]
when computed at $x=0$, vanishes.
\EP
%%%%%%%%%%%%%%%%%%%%%%%%%%%%%%%%%%%%%%%%%%%%%%%%%%%%%%%%%%%%%%%%%%%%%%%%%%

%%%%%%%%%%%%%%%%%%%%%%%%%%%%%%%%%%%%%%%%%%%%%%%%%%%%%%%%%%%%%%%%%%%%%%%% 
% Figure T1 + T2
%%%%%%%%%%%%%%%%%%%%%%%%%%%%%%%%%%%%%%%%%%%%%%%%%%%%%%%%%%%%%%%%%%%%%%%% 
\begin{figure}[ht] 
\vspace{-.3cm}
\centering 
\ins{200pt}{-100pt}{$T$}
\ins{160pt}{-090pt}{$T_2$}
\ins{270pt}{-090pt}{$T_1$}
\ins{090pt}{-063pt}{$\ell_T$}
\ins{350pt}{-063pt}{$\ell_T'$}
\ins{220pt}{-065pt}{$\ell$}
\includegraphics[width=4.2in]{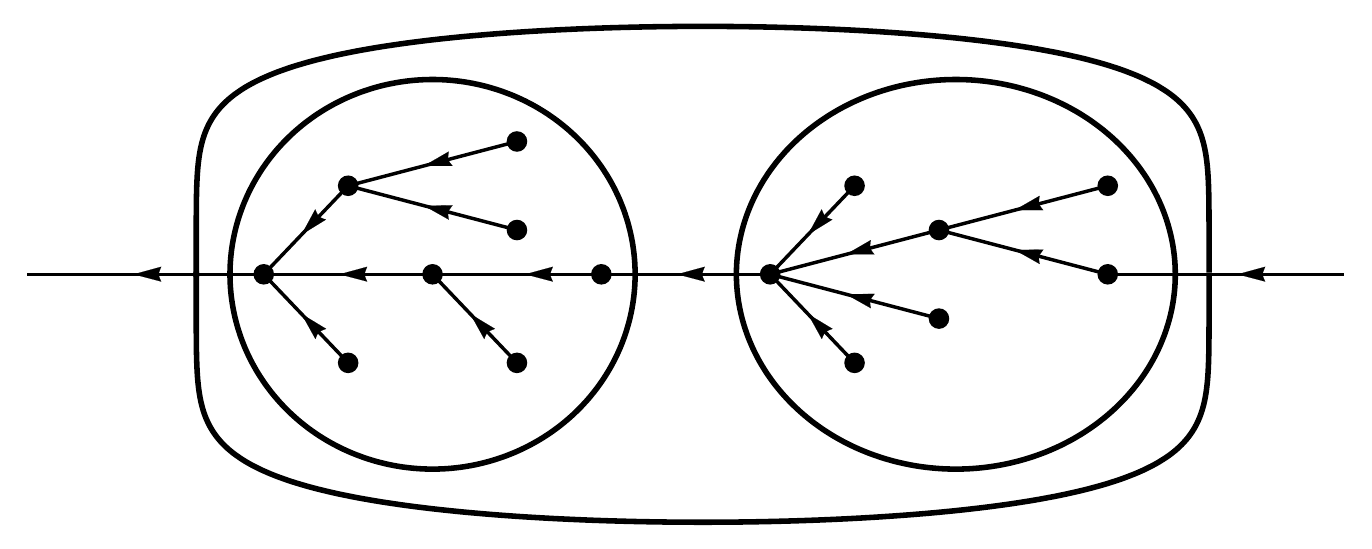} 
\vskip.2truecm 
\caption{The sets $T_1$ and $T_2$ inside a RC $T$.}
\label{canc2} 
\end{figure} 
%%%%%%%%%%%%%%%%%%%%%%%%%%%%%%%%%%%%%%%%%%%%%%%%%%%%%%%%%%%%%%%%%%%%%%%%

%%%%%%%%%%%%%%%%%%%%%%%%%%%%%%%%%%%%%%%%%%%%%%%%%%%%%%%%%%%%%%%%%%%%%%%%%%
\begin{ex} \label{effe1}
\emph{
An example of RC $T'\in \matF_1(T)$, with $T$ as in Figure \ref{canc2}, is given in Figure \ref{canc2bis};
the arrows of the lines along the path connecting $v_T$ and $v_{T'}$ are reverted.
}
\end{ex}
%%%%%%%%%%%%%%%%%%%%%%%%%%%%%%%%%%%%%%%%%%%%%%%%%%%%%%%%%%%%%%%%%%%%%%%%%%

%%%%%%%%%%%%%%%%%%%%%%%%%%%%%%%%%%%%%%%%%%%%%%%%%%%%%%%%%%%%%%%%%%%%%%%% 
% Figure T' 1
%%%%%%%%%%%%%%%%%%%%%%%%%%%%%%%%%%%%%%%%%%%%%%%%%%%%%%%%%%%%%%%%%%%%%%%% 
\begin{figure}[ht] 
\vspace{-.3cm}
\centering 
\ins{200pt}{-100pt}{$T'$}
\ins{160pt}{-090pt}{$T_2'$}
\ins{270pt}{-090pt}{$T_1'$}
\ins{090pt}{-045pt}{$\ell_{T'}$}
\ins{350pt}{-074pt}{$\ell_{T'}'$}
\ins{220pt}{-065pt}{$\ell$}
\includegraphics[width=4.2in]{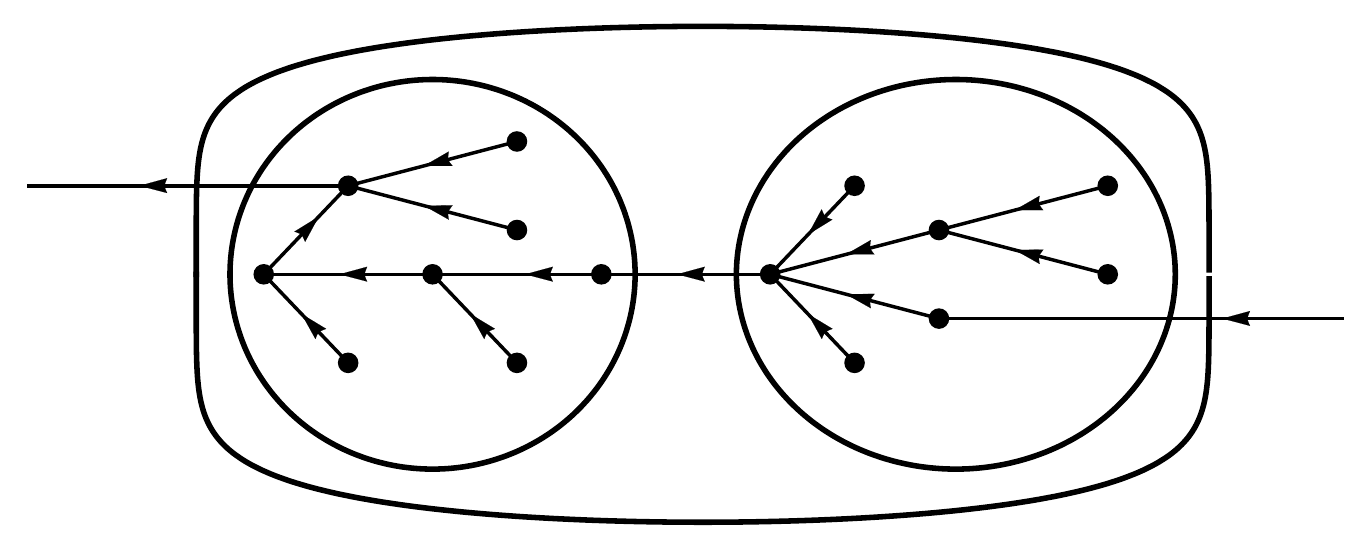} 
\caption{A RC $T'\in \calF_1(T)$, with $T$ as in Figure \ref{canc2}.}
\label{canc2bis} 
\end{figure} 
%%%%%%%%%%%%%%%%%%%%%%%%%%%%%%%%%%%%%%%%%%%%%%%%%%%%%%%%%%%%%%%%%%%%%%%%

%%%%%%%%%%%%%%%%%%%%%%%%%%%%%%%%%%%%%%%%%%%%%%%%%%%%%%%%%%%%%%%%%%%%%%%%%%
\begin{ex} \label{effe2}
\emph{
Similarly, an example of RC $T'\in \calF_2(T)$ is given in Figure \ref{canc2ter};
in such a case all the arrows along the path $\calP_{T'}$, including $\ell$, are reverted. 
Note that in Figure \ref{canc2ter} the line $\ell_{T'}$ exiting $T'$ is on the right, while the line entering $T'$ is on the left,
so that, in order to draw the RC with the arrow pointing to the right, according to the convention introduced in Section \ref{alberi},
one should consider the mirrored image of Figure \ref{canc2ter}, as represented in Figure \ref{canc2quater}.
}
\end{ex}
%%%%%%%%%%%%%%%%%%%%%%%%%%%%%%%%%%%%%%%%%%%%%%%%%%%%%%%%%%%%%%%%%%%%%%%%%%

Thanks to Corollary \ref{cancellazione} and Lemma \ref{reverse}, substituting in \eqref{taylor} we obtain
\begin{equation} \nonumber 
\matM_{n}^{(k)} (\om\cdot\nu_{\ell'_T}) =  (\om\cdot\nu_{\ell'_T})^2\matR^{(k)}_n (\om\cdot\nu_{\ell'_T}) .
\end{equation}
So we have to control the remainder $\matR^{(k)}_n (\om\cdot\nu_{\ell'_T})$.

%%%%%%%%%%%%%%%%%%%%%%%%%%%%%%%%%%%%%%%%%%%%%%%%%%%%%%%%%%%%%%%%%%%%%%%%%%
% Figure T' 2
%%%%%%%%%%%%%%%%%%%%%%%%%%%%%%%%%%%%%%%%%%%%%%%%%%%%%%%%%%%%%%%%%%%%%%%%%%
\begin{figure}[ht] 
\vspace{-.3cm}
\centering 
\ins{200pt}{-100pt}{$T'$}
\ins{160pt}{-088pt}{$T_2'$}
\ins{270pt}{-088pt}{$T_1'$}
\ins{086pt}{-083pt}{$\ell_{T'}'$}
\ins{350pt}{-045pt}{$\ell_{T'}$}
\ins{212pt}{-065pt}{$\ell$}
\includegraphics[width=4.2in]{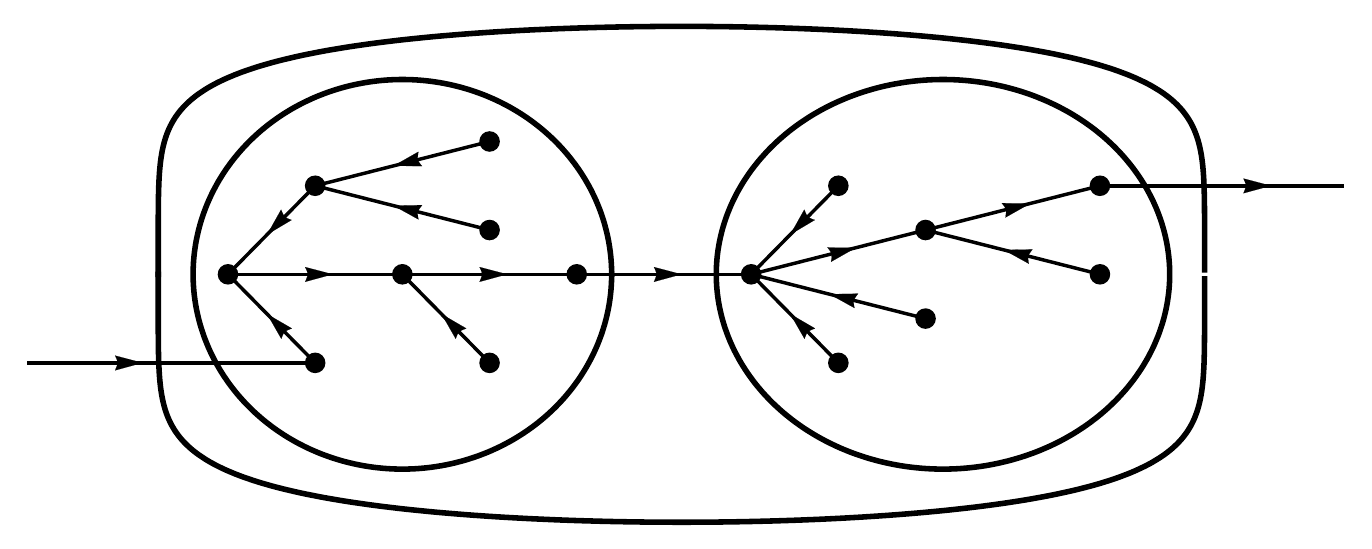} 
\caption{A RC $T'\in \calF_2(T)$, with $T$ as in Figure \ref{canc2}.}
\label{canc2ter} 
\end{figure} 
%%%%%%%%%%%%%%%%%%%%%%%%%%%%%%%%%%%%%%%%%%%%%%%%%%%%%%%%%%%%%%%%%%%%%%%%%%

%%%%%%%%%%%%%%%%%%%%%%%%%%%%%%%%%%%%%%%%%%%%%%%%%%%%%%%%%%%%%%%%%%%%%%%%%%
% Figure T' 2
%%%%%%%%%%%%%%%%%%%%%%%%%%%%%%%%%%%%%%%%%%%%%%%%%%%%%%%%%%%%%%%%%%%%%%%%%%
\begin{figure}[H] 
\vspace{-.2cm}
\centering 
\ins{200pt}{-100pt}{$T'$}
\ins{140pt}{-078pt}{$T_1'$}
\ins{250pt}{-088pt}{$T_2'$}
\ins{086pt}{-063pt}{$\ell_{T'}$}
\ins{354pt}{-064pt}{$\ell_{T'}'$}
\ins{215pt}{-065pt}{$\ell$}
\includegraphics[width=4.2in]{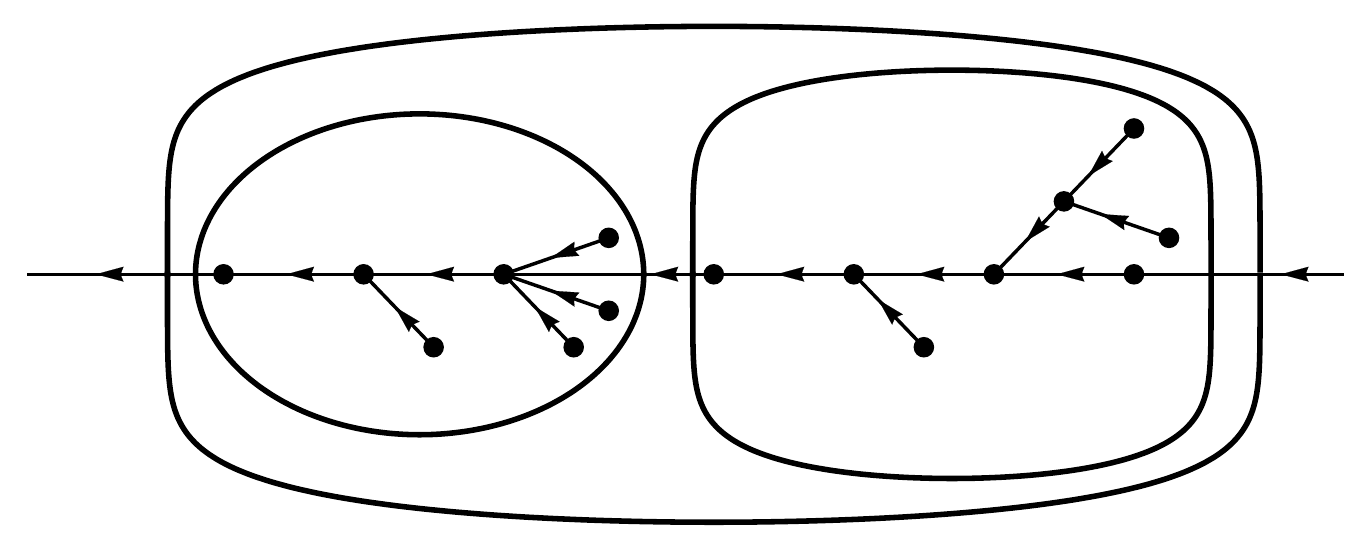} 
\caption{The RC $T'$ of Figure \ref{canc2ter} redrawn with the exiting line on the left and the entering line on the right.}
\label{canc2quater} 
\end{figure} 
%%%%%%%%%%%%%%%%%%%%%%%%%%%%%%%%%%%%%%%%%%%%%%%%%%%%%%%%%%%%%%%%%%%%%%%%%%

%%%%%%%%%%%%%%%%%%%%%%%%%%%%%%%%%%%%%%%%%%%%%%%%%%%%%%%%%%%%%%%%%%%%%%%%%%
%%%%%%%%%%%%%%%%%%%%%%%%%%%%%%%%%%%%%%%%%%%%%%%%%%%%%%%%%%%%%%%%%%%%%%%%%%
\zerarcounters 
\section{Renormalization and convergence} 
\label{R&C} 
%%%%%%%%%%%%%%%%%%%%%%%%%%%%%%%%%%%%%%%%%%%%%%%%%%%%%%%%%%%%%%%%%%%%%%%%%%
%%%%%%%%%%%%%%%%%%%%%%%%%%%%%%%%%%%%%%%%%%%%%%%%%%%%%%%%%%%%%%%%%%%%%%%%%%

To conclude the proof of the convergence of the Lindstedt series we have still to prove
that $\|\matR^{(k)}_n(x)\|_{\rm op} \le C^k$ for some constant $C>0$ independent of $n$. 

%%%%%%%%%%%%%%%%%%%%%%%%%%%%%%%%%%%%%%%%%%%%%%%%%%%%%%%%%%%%%%%%%%%%%%%%%%
\begin{ex} \label{SKn-example2-bis}
\emph{
The analysis of Section \ref{canc} allows us to solve the problem related to the chain of RCS appearing in the tree represented
in Figure \ref{catenabis}. Indeed, reasoning as in Example \ref{SKn-example2} and using the notation therein,
we see that $\calW_{T'}(x)$ depends on $x$ only when the RC $T'$ replacing $T_i$ is as on the left of Figure \ref{sec2}.
Thus, we can compute explicitly $\matR^{(2)}_n(\om\cdot\nu_{\ell_i})$ for all $i=1,\ldots,p$.
Indeed, setting $x=\om\cdot\nu$, first of all we note that we can confine to the case
$|x|<|\om\cdot\nu_0|/2$, since otherwise we can bound $|\om\cdot\nu|\le 2/|\om\cdot\nu_0|$.
Then, we obtain find
\begin{equation} \label{R2}
\begin{aligned}
x^2 \matR^{(2)}_n(x) & = 
\matM_{n}^{(k)} (x) - \matM_{n}^{(k)} (0) - x \partial_x \matM_{n}^{(k)} (0) \\
& =
\sum_{\substack{ \nu_0 \in \ZZZ^\ZZZ_f \\ n'<n}}
(\ii \nu_0) |\nu_0|^2 |f_{\nu_0}|^2 (-\ii\nu_0) \left(
\frac{\Psi_{n'} (\om\cdot\nu_0+x)}{(\om\cdot\nu_0 + x)^2} -  \frac{\Psi_{n'} (\om\cdot\nu_0)}{(\om\cdot\nu_0)^2}
+ \frac{2 x \, \Psi_{n'} (\om\cdot\nu_0)}{(\om\cdot\nu_0)^3} \right) .
\end{aligned}
\end{equation}
Assume first that there exists $n_1,n_2<n$ such that $\Psi_{n_1} (\om\cdot\nu_0)=\Psi_{n_2} (\om\cdot\nu_0+x)=1$.
Then, we have
\[
\begin{aligned}
\frac{1}{(\om\cdot\nu_0 + x)^2} -  \frac{1}{(\om\cdot\nu_0)^2} + \frac{2 x }{(\om\cdot\nu_0)^3} &
= \frac{(\om\cdot\nu_0)^3 - (\om\cdot\nu_0)(\om\cdot\nu_0+x)^2 + 2x(\om\cdot\nu_0+x)^2}{(\om\cdot\nu_0)^3(\om\cdot\nu_0+x)^2} \\
& \frac{3(\om\cdot\nu_0)\,x^2 + x^3}{(\om\cdot\nu_0)^3(\om\cdot\nu_0+x)^2} ,
\end{aligned}
\]
so that, using that $2|x|<|\om\cdot\nu_0|$, we obtain from \eqref{R2}
\[
\left| \matR^{(2)}_n(\om\cdot\nu_{\ell_i}) \right| \le {\sum_{\nu_0 \in \ZZZ^\ZZZ_f}}^{\!\!\!*,n} \;
\frac{4\| \nu_0 \| \| \nu_0 \|_1|^3 |f_{\nu_0}|^2}{(\om\cdot\nu_0)^4} , \qquad i=1,\ldots, p.
\]
Using the last bound we conclude that, even though the value of the tree $\vartheta$ considered in Example \ref{aiuto}
is only bounded as in \eqref{above}, nevertheless, when we consider all the tree which can be obtained from $\vartheta$
by replacing each RC $T_i$ with any other possible RC, we obtain an overall contribution, given by sums of terms each of which is
estimated, possibly with a different constant $C$ w.r.t.~\eqref{above}, as
\begin{equation} \nonumber 
\begin{aligned}
& \frac{4^{\frac{k-1}{2}}\|\nu_0\|_1^{2(k-1)} |f_{\nu_0}|^{k-1}}{2^{\frac{k-1}{2}} |\om\cdot\nu_0|^{2(k-1)}}
\frac{\|\nu\|_1|f_\nu|}{|\om\cdot\nu|^{2}} 
\sim C^k (\h_{\nu_0})^{k-1} \h_\nu
\g^{-(k+1)} \|\nu\|_1^{2\tau} e^{-2s\|\nu\|_1}
\end{aligned}
\end{equation}
where we can bound $ \|\nu\|_1^{2\tau} e^{-2s\|\nu\|_1} \le e^{-s\|\nu\|_1} s^{-2\tau} \lceil 2 \tau \rceil !$.
Finally, consider the case in which either $\Psi_{n'} (\om\cdot\nu_0)=1$ for some $n'<n$ and $\Psi_{n'} (\om\cdot\nu_0+x)=0$
for all $n'<n$ or, \emph{vice versa}, $\Psi_{n'} (\om\cdot\nu_0+x)=1$ for some $n'<n$ and $\Psi_{n'} (\om\cdot\nu_0)=0$ for all $n'<n$.
If this happens, then the RC does not represent a problem. To see this,
consider explicitly the first case (the second one can be discussed in the same way).
First of all note that, in such a case, $|\om\cdot\nu|$ cannot be smaller than $|\om\cdot\nu_0|/2$
because the condition $K(T_i) < 2^{m_{n}-1}$ (see item \ref{item3} in Definition \ref{RC}) implies that $|\om\cdot\nu_0|>\be^\star_{\om}(m_{n-1})$
and hence $|\om\cdot\nu_0+x|>\be_{\om}^\star (m_{n-1})/2$, so that there exists $n'<n$ such that $\Psi_{n'} (\om\cdot\nu_0+x)=1$.
On the other hand if $|\om\cdot\nu| \ge |\om\cdot\nu_0|/2$ then in \eqref{above}
we can bound $|\om\cdot\nu|^{-(k-1)}$ with $2^{k-1}|\om\cdot\nu_0|^{-(k-1)}$,
which leads once more to the bound \eqref{above} without exploiting any cancellations.
}
\end{ex}
%%%%%%%%%%%%%%%%%%%%%%%%%%%%%%%%%%%%%%%%%%%%%%%%%%%%%%%%%%%%%%%%%%%%%%%%%%

%%%%%%%%%%%%%%%%%%%%%%%%%%%%%%%%%%%%%%%%%%%%%%%%%%%%%%%%%%%%%%%%%%%%%%%%%%
\begin{rmk} \label{whysmooth}
\emph{
The argument in Example \ref{SKn-example2-bis} hints at the advantage in using smooth scale function:
this allows to estimate the remainder $\matR^{(k)}_n (\om\cdot\nu_{\ell'_T})$ in \eqref{taylor} through the integral form
\[
x^2 \matR^{(k)}_n(x) = x^2 \int_0^1 \der y \, \partial_x^2 \matM_n^{(k)}(yx) ,
\]
without worrying that derivatives can be applied to the scale functions. In particular \eqref{R2} becomes
\[
x^2 \matR^{(2)}_n(x) = \sum_{\substack{ \nu_0 \in \ZZZ^\ZZZ_f \\ n'<n}}
(\ii \nu_0) |\nu_0|^2 |f_{\nu_0}|^2 (-\ii\nu_0) \, \partial_x^2 \left( \frac{\Psi_{n'} (\om\cdot\nu_0+x)}{(\om\cdot\nu_0+x)^2} \right) .
\]
}
\end{rmk}
%%%%%%%%%%%%%%%%%%%%%%%%%%%%%%%%%%%%%%%%%%%%%%%%%%%%%%%%%%%%%%%%%%%%%%%%%%

However, extending the argument in Example \ref{SKn-example2-bis} to any possible RC is delicate,
because of the following issue: a RC may contain another RC, which in turn may contain another RC and so on.
The problem is better pinpointed using smooth scale functions, in terms of accumulation of derivatives.
Suppose that a tree contains a set of nested RCs $T_1\supset T_2 \ldots \supset T_h$ such that
$\calP_{T_h} \subset \ldots \subset \calP_{T_2} \subset \calP_{T_1}$.
When considering the  expansion \eqref{quello} for each of them,
it may happen that the derivatives act all on the propagator of a single line $\ell\in L(T_h)$,
so that $\calG_{n_\ell}$ is differentiated $2h$ times.
This leads to
\[
\del_x^{2h} \calG_{n_{\ell}}(\om\cdot\nu_{\ell}^0+x)
\sim \frac{(2h+1)!}{(\om\cdot\nu_{\ell}^0+x)^{2+2h}}
\]
Since $h$ can be $O(k)$, this may jeopardize the gain
$(\om\cdot\nu_{\ell}^0+x)^2$
obtained when summing together the values
of all RCs which may replace the RC $T_h$.

A crucial observation is that, as Example \ref{aiuto} shows, problems arise because the presence of arbitrarily long chains of RCs.
Indeed, if only chains of at most $\gotn_0$ RCs were possible, with $\gotn_0\in\NNN$ fixed once and for all,
then the number of resonant lines on scale $\ge n$ would be bounded by $\gotn_0\gotN_n(\vartheta)$,
because the line entering any chain with resonant lines on scale $n$ is non-resonant and hence contributes to $\gotN_n(\vartheta)$.
Then, we could estimate
\[
\prod_{\ell\in L(\vartheta)}|\calG_\ell| \leq\prod_{n\ge0}\left(\frac{4}{(\be^\star_\om(m_n))^2}\right)^{(1+\gotn_0)\gotN_n(\vartheta)} ,
\]
and the same argument as in Lemma \ref{convergerebbe} would be enough to obtain a bound like \eqref{magara}.

A byproduct of the argument above is that for any chain of RCs what we really need is to obtain a gain for all but a finite number $\gotn_0$ of them.
Therefore, we can reason as follows. Suppose that a RC $T$ is not contained in any other RC and contains one or more chains $\gotC$ of RCs
such that $\calP_{T'} \subset \calP_T$ for all $T'\in\gotC$,
and write $\matW(T)$ according to \eqref{quello}, with
\begin{equation} \label{RT}
R_T(x) = x^2 \int_0^1 \der y \, \partial_x^2 \matW_T(yx) 
\end{equation}
and
\begin{equation} \label{valfunderder}
\del_x^2 \matW_T(x) =
\frac{1}{p_{v_T}!}\ii \nu_{v_T} f_{\nu_{v_T}} \Biggl( \prod_{v\in N(\vartheta)\setminus \{v_T\}} 
\frac{1}{p_{v}!} (\ii \nu_{\pi(v)}\cdot \ii \nu_v ) f_{\nu_v} \Biggl) 
\del_x^2 \!\!\! 
\prod_{\ell\in L(T)} \calG_{n_\ell}(\x_\ell(x))) 
\ii \nu_{v'_T} .
\end{equation}
The second derivative of $\matW_T(yx)$ in \eqref{RT} produces several contributions: in each of them
the derivatives are applied either both to the propagator of the same line or to the propagators of two different lines.
\begin{itemize}[topsep=.5ex]
\itemsep0em
\item If the differentiated lines are outside all the RCs of the chains $\gotC$, then for each chain $\gotC$ and each $T'\in\gotC$ we can write
$\matW(T')$ as in \eqref{quello},
and ignore the contributions $\matW_{T'}(0)$ and $x \partial_x \matW_{T'}(0)$ since we know that,
when we sum the contributions of all possible RCs, they cancel out.
\item
If both differentiated lines are contained inside the
same RC $T_1$ of a chain $\gotC$, then we can expand again $\matW(T')$ as in \eqref{quello} for all $T'\in\gotC\setminus\{T_1\}$
and all $T'$ contained in the other chains,
and ignore the contributions $\matW_{T'}(0)$ and $x \partial_x \matW_{T'}(0)$.
\item
If the differentiated lines are contained in two different RCs $T_1$ and $T_2$, which may belong either to the same chain
or two different chains $\gotC$, we write $\matW(T')$ as in \eqref{quello} for all the $T'\notin\{T_1,T_2\}$.,
and, once more, ignore the contributions $\matW_{T'}(0)$ and $x \partial_x \matW_{T'}(0)$.
\end{itemize}
In such a way we obtain a gain $(\om\cdot\nu_{\ell_{T'}'})^2$ for all $T'\in\gotC$ except one or two of them.
Then, we can iterate the procedure, starting from each RC $T'$ which belongs to some chain $\gotC$ in $T$ and
which does not contain any differentiated lines produced by \eqref{RT}: we write $R_{T'}(x)$ as in \eqref{valfunderder},
with $T'$ instead of $T$, and repeat the construction above.
Eventually, we find that, in each chain of RCs, we obtain a gain factor $(\om\cdot\nu_{\ell_{T'}^0}+yx)^2$ for each RC $T'$ of the chain
except at most two of them.

\medskip 

%%%%%%%%%%%%%%%%%%%%%%%%%%%%%%%%%%%%%%%%%%%%%%%%%%%%%%%%%%%%%%%%%%%%%%%%%% 
\noindent\emph{Proof of Proposition \ref{esiste1}.}
The conclusion of the discussion above is that the bound in Lemma \ref{convergerebbe}
can be extended by replacing $\Val_{N\!R}(\vartheta)$ with $\Val(\vartheta)$,
provided the constant $C_0$ in \eqref{C0} be replaced with a new constant
\begin{equation} \label{C0'}
C_0' = c_0 \left( \frac{16 c_1}{s \be_\om^\star(m_{n_0})} \right)^4 \|f\|_{2s,\star,\RRR} ,
\end{equation}
for a suitable positive constant $c_0$. This proves the existence of a solution $u\in \matH_\star^{s}(\TTT^\ZZZ,\ell^\io(\RRR))$ to \eqref{lagrapde}
for all $\e\in(-\e_1,\e_1)$, with $\e_1=(C_0')^{-1}$ depending on $s$, $\om$ and $\|f\|_{s,\star,\RRR}$.
\EP
%%%%%%%%%%%%%%%%%%%%%%%%%%%%%%%%%%%%%%%%%%%%%%%%%%%%%%%%%%%%%%%%%%%%%%%%%% 

%%%%%%%%%%%%%%%%%%%%%%%%%%%%%%%%%%%%%%%%%%%%%%%%%%%%%%%%%%%%%%%%%%%%%%%%%%
%%%%%%%%%%%%%%%%%%%%%%%%%%%%%%%%%%%%%%%%%%%%%%%%%%%%%%%%%%%%%%%%%%%%%%%%%%
\zerarcounters 
\section{Uniform bounds}
\label{10} 
%%%%%%%%%%%%%%%%%%%%%%%%%%%%%%%%%%%%%%%%%%%%%%%%%%%%%%%%%%%%%%%%%%%%%%%%%%
%%%%%%%%%%%%%%%%%%%%%%%%%%%%%%%%%%%%%%%%%%%%%%%%%%%%%%%%%%%%%%%%%%%%%%%%%%

By using the explicit expression of $\e_1$ provided at the end of Section \ref{R&C} we can go through the last step we need to conclude the
proof of Theorem \ref{main}, as outlined at the end of Section \ref{functional}.

For all $\om\in \DD_{\g,\mu_1,\mu_2}$ with $\mu_1>1+q$ and $\mu>1$,
by Proposition \ref{brubru} we can bound  $\be_\om^\star(m_{n_0}) \ge \be_*(m_{n_0})$ in \eqref{C0'}.
More generally, for all $\om\in\gotB^\star$ such that $\be^\star_\om(m) \ge \be_*(m)$,
that is for all $\om$ in a subset of $\gotB^\star$ which contains the set of Diophantine vectors $\DD_{\g,\mu_1,\mu_2}$,
the value of $\e_1$ in Proposition \ref{esiste1} can be estimated uniformily. Indeed, if we set
\begin{equation} \label{estarol}
\ol{\e}_{1} = \left( c_0 \left( \frac{16 c_1}{s \be_*(m_{n_0})} \right)^4 \|f\|_{2s,\star,\RRR} \right)^{-1} \!\!\! ,
\end{equation}
then we obtain $\e_1 \ge \ol{\e}_1$. 

%%%%%%%%%%%%%%%%%%%%%%%%%%%%%%%%%%%%%%%%%%%%%%%%%%%%%%%%%%%%%%%%%%%%%%%%%%
\begin{rmk} \label{justrecall}
\emph{
By looking at \eqref{C0'}, we see that $\ol{\e}_1$ is proportional to $(s \be_*(m_{n_0})/\|f\|_{2s,\star,\RRR})^4$ and hence
it depends on $\om$ only through the sequence $\{\be_*(m)\}_{m\in\NNN}$, i.e.~only through the values if $\g$, $\mu_1$ and $\mu_2$.
If we recall the definition of $\be_{*}(m)$ in \eqref{betam} and define $m_{n_0}$ by requiring that
(compare with \eqref{bryunoo})
\begin{equation} \label{bryunooo}
8 \sum_{n>n_0 } \frac{1}{2^{m_n}}\log \left(\frac{1}{\be_*(m_n)} \right)< \frac{s}{2} ,
\end{equation}
we conclude that $\ol{\e}_1$ depend on $\g$, $\mu_1$, $\mu_2$, $s$ and $\s$. In particular, \eqref{bryunoo} yields that
$\ol{\e}_1$ is proportional to $\g^4$. In fact, such a bound may be improved by scaling arguments (see Remark \ref{finale}).
}
\end{rmk}
%%%%%%%%%%%%%%%%%%%%%%%%%%%%%%%%%%%%%%%%%%%%%%%%%%%%%%%%%%%%%%%%%%%%%%%%%%

%%%%%%%%%%%%%%%%%%%%%%%%%%%%%%%%%%%%%%%%%%%%%%%%%%%%%%%%%%%%%%%%%%%%%%%%%%
\noindent\emph{Proof of Proposition \ref{bieco}.}
Let $\e_1$ be such that $\e_1 C_0'=1$, with $C_0'$ in \eqref{C0'}. For all $\e\in(0,\e_1)$, we have
\[
\| u \|_{s,\star,\ell^\io(\RRR)}
\le \sum_{k\ge 1} \Bigl( \frac{\e}{\e_1} \Bigr)^k \le \frac{\e}{\e_1-\e} ,
\]
so that, if we require $\e_2\in(0,\e_1)$ to be such that $\| u \|_{s,\star,\ell^\io(\RRR)} \le (s-s_2)/2$ (see Proposition \ref{esiste2}), we need
\[
|\e| \le \e_2 , \qquad \e_2 := \frac{s-s_2}{s-s_2+2} \e_1 \ge \ol{\e}_2 := \frac{s-s_2}{s-s_2+2} \ol{\e}_1,
\]
Therefore, $\ol{\e}_2$ depends on $\s$, $s$, $s_2$, $\|f\|_{2s,\star,\RRR}$ and the sequence $\{\be_*(m)\}_{m\in\NNN}$.
This implies that, for all $\om\in\Omega_\e$, with $\Omega_\e$ defined as at the end of Section \ref{functional},
$\ol{\e}_2$ depends on $\g$, $\mu_1$, $\mu_2$, $\s$, $s$ and $s_2$ (that we can fix as $s_2=s/2$, as observed
at the end of Section \ref{functional}). In particular, $\ol{\e}_2$ depends on $f$ only through
$s$ and $\|f\|_{2s,\star,\RRR}$ and on $\om$ only through $\g$, $\mu_1$ and $\mu_2$.
Then, if we define $\e_*$ and $\e_\star$ as in \eqref{e*} and in \eqref{estar}, respectively,
one finds $\e_\star \ge \e_* \ge \ol{\e}_2$, and all the assertions in Proposition \ref{bieco} follow.
%%%%%%%%%%%%%%%%%%%%%%%%%%%%%%%%%%%%%%%%%%%%%%%%%%%%%%%%%%%%%%%%%%%%%%%%%%

%%%%%%%%%%%%%%%%%%%%%%%%%%%%%%%%%%%%%%%%%%%%%%%%%%%%%%%%%%%%%%%%%%%%%%%%%%
\begin{rmk} \label{finale}
\emph{
If, for any $\om\in \DD_{\g,\mu_1,\mu_2}$, we set $\om=\g \om_*$ since the beginning, then
we can work with $\om_*$ instead of $\om$, and all the analysis from Section \ref{brutale} on can be repeated word by word,
provided that $\be^\star_\om(m)$ and $\be_*(m)$ are replaced, respectively, with the $\g$-independent
\[
\ol{\be}^\star_{\om_*}(m):=\g^{-1}\be^\star_\om(m) , \qquad \ol{\be}_{*}(m):=K_1^{-1}e^{-K_2 2^m (\log 2^m)^{-(\s-1)}} ,
\]
and the first inequality in \eqref{ultimaformulanumerata} is replaced with
\[
\prod_{\ell\in L_{\!N\!R}(\vartheta)}|\calG_\ell| \le \g^{-2k} \prod_{n\ge0} \biggl(\frac{4}{(\ol{\be}^\star_\om(m_n))^2}\biggr)^{\gotN_n(\vartheta)} 
\le \g^{-2k} \prod_{n\ge0} \left(\frac{4}{(\ol{\be}_{*}(m_n))^2}\right)^{\gotN_n(\vartheta)} ,
\]
where we have used that any tree of order $k$ has $k$ lines.
Therefore, the only difference in the value of the constant $C_0'$ is that $\be_*(m_{n_0})$ is replaced with
$\ol{\be}_{*}(m_{n_0})$ and a factor $\g^2$ appears in front of $c_0$.
}
\end{rmk}
%%%%%%%%%%%%%%%%%%%%%%%%%%%%%%%%%%%%%%%%%%%%%%%%%%%%%%%%%%%%%%%%%%%%%%%%%%

%%%%%%%%%%%%%%%%%%%%%%%%%%%%%%%%%%%%%%%%%%%%%%%%%%%%%%%%%%%%%%%%%%%%%%%%%%
%%%%%%%%%%%%%%%%%%%%%%%%%%%%%%%%%%%%%%%%%%%%%%%%%%%%%%%%%%%%%%%%%%%%%%%%%%

\end{document}